\documentclass[11pt]{article}
\usepackage{url}
\usepackage[margin=1.0in]{geometry}
\usepackage[T1]{fontenc}
\usepackage[utf8]{inputenc}
\usepackage{lmodern}
\usepackage{microtype}
\usepackage{amsmath,amssymb,amsthm,mathtools}
\usepackage{enumitem}
\usepackage{titlesec}

\usepackage{comment}
\usepackage[colorlinks=true,linkcolor=blue,citecolor=blue,urlcolor=blue]{hyperref}
\usepackage{tikz}
\usetikzlibrary{arrows.meta,calc,decorations.pathreplacing}
\usetikzlibrary{arrows.meta,shapes.geometric,calc}
\DeclareMathOperator{\Span}{Span}
\usetikzlibrary{positioning,decorations.pathreplacing}
\setlist{nosep,leftmargin=*}
\allowdisplaybreaks
\titlespacing*{\section}{0pt}{2ex plus .5ex minus .2ex}{1ex plus .2ex}
\titlespacing*{\subsection}{0pt}{1.4ex plus .4ex minus .2ex}{.7ex plus .2ex}

\theoremstyle{plain}
\newtheorem{theorem}{Theorem}[section]
\newtheorem{lemma}[theorem]{Lemma}

\newtheorem{problem}[theorem]{Problem}
\newtheorem{claim}[theorem]{Claim}

\newtheorem{case}{Case}
\newtheorem{fact}{Fact}
\newtheorem{conjecture}[theorem]{Conjecture}
\newtheorem{observation}[theorem]{Observation}
\theoremstyle{definition}
\newtheorem{definition}[theorem]{Definition}

\newcounter{denseproperty}
\renewcommand{\thedenseproperty}{(U\arabic{denseproperty})}

\DeclareMathOperator{\Cay}{Cay}
\DeclareMathOperator{\ord}{ord}
\DeclareMathOperator{\Stab}{Stab}
\DeclareMathOperator{\aff}{aff}

\newcommand{\Prob}{\mathbb P}
\newcommand{\diss}{\dim_{\mathrm{diss}}}
\newcommand{\F}{\mathbb F}

\title{Graham's rearrangement conjecture in finite abelian groups}
\author{Peiru Kuang\footnote{School of Mathematical Sciences, Shanghai Jiao Tong University, Shanghai 200240, China. Supported by Shanghai Institute for Mathematics and Interdisciplinary Sciences, SIMIS-ID-26-AMS-003. Email: \url{peiru_k@sjtu.edu.cn}}
\and
Yan Wang\footnote{School of Mathematical Sciences, Shanghai Jiao Tong University, Shanghai 200240, China. Supported by National Key R\&D Program of China under Grant No. 2022YFA1006400 and National Natural Science Foundation of China under Grant No. 12571376. Email: \url{yan.w@sjtu.edu.cn}}}
\date{}

\begin{document}
\maketitle
\definecolor{samplingblue}{RGB}{43,79,122}
\definecolor{samplinggreen}{RGB}{36,116,88}
\definecolor{samplinggray}{RGB}{105,112,120}
\begin{abstract}
Graham conjectured in 1971 that every subset of the nonzero elements of a prime cyclic group admits an ordering with pairwise distinct nonempty partial sums. Recently, the work of Pham and Sauermann, together with earlier results, settled this conjecture for all sufficiently large primes. In this paper, we prove that there is an absolute integer $n_0$ such that every subset $S\subseteq G\setminus\{0\}$ of a finite abelian group with $|S|\ge n_0$ admits such an ordering, which also solves a conjecture of Alspach and Liversidge for all sufficiently large sets.

A key ingredient of our proof is a sufficient condition via additive dimension: If $|S|=n$ and $S$ contains a dissociated subset of size at least $\Omega((\log n)^2/\log\log n)$, then $S$ has an ordering in which no proper nonempty interval of consecutive terms sums to zero. This bound is within a factor of $O(\log n/\log\log n)$ of the smallest possible additive dimension of an $n$-element set.
\end{abstract}

\section{Introduction}
Let $G$ be a finite abelian group, written additively, and let $S\subseteq G\setminus\{0\}$ have size $n$. An ordering $s_1,\ldots,s_n$ of $S$ is \emph{valid} if its nonempty partial sums $s_1,s_1+s_2,\ldots,s_1+\cdots+s_n$ are pairwise distinct. We call $S$ \emph{rearrangeable} if it admits a valid ordering. In 1971, Graham~\cite{Graham} raised the following conjecture, which was later reiterated by Erd\H{o}s and Graham~\cite{ErdosGraham}.

\begin{conjecture}[Graham]\label{conj:graham}
If $p$ is prime, then every $S\subseteq\mathbb Z_p\setminus\{0\}$ is rearrangeable.
\end{conjecture}
In his original paper, Graham~\cite{Graham} posed the rearrangement problem alongside several questions about sumset structures, a central topic in additive combinatorics. Later, he suggested that the problem was also motivated by applications to juggling~\cite{BuhlerEtAlJuggling,ChungGrahamJuggling,GrahamJuggling}; also see~\cite{BucicEtAl}. 

Questions about distinct partial sums first arose in the construction of combinatorial designs. Motivated by constructions of complete Latin squares, Gordon~\cite{Gordon} asked in 1961 which finite groups admit an ordering of all their elements with distinct partial products, and characterized the abelian groups with this property. Such an ordering gives a complete Latin square by permuting the rows and columns of the multiplication table. A Latin square is \emph{complete} if every ordered pair of distinct symbols occurs exactly once in horizontally adjacent cells and exactly once in vertically adjacent cells. 
For $1\le i<j\le n$, the partial sums $s_1+\cdots+s_i$ and $s_1+\cdots+s_j$ coincide precisely when $s_{i+1}+\cdots+s_j=0$. Thus, rearrangeability asks for an ordering in which every consecutive interval starting after the first term has nonzero sum.
This viewpoint connects rearrangement problems with classical questions in additive combinatorics.
For example, the Erd\H{o}s--Ginzburg--Ziv theorem
~\cite{ErdosGinzburgZiv} states that every sequence of $2m-1$ elements of $\mathbb Z_m$ contains $m$ terms whose sum is zero.

Moving beyond prime orders, Archdeacon, Dinitz, Mattern and Stinson~\cite{ArchdeaconDinitzMatternStinson} conjectured the following.
\begin{conjecture}[Archdeacon, Dinitz, Mattern and Stinson~\cite{ArchdeaconDinitzMatternStinson}]\label{conj:Arch}
For every positive integer $n$, every subset $S \subseteq \mathbb{Z}_n \setminus \{0\}$ is rearrangeable.
\end{conjecture}
A slightly stronger conjecture is due to Alspach and Liversidge~\cite{AlspachLiversidge} (see also~\cite{BBKMM,BucicEtAl,CostaDellaFioreWeak,CostaDellaFioreFontanaVena,CostaDellaFioreOllisRovnerFrydman,CMPPPartialSums,CostaPellegrini,HicksOllisSchmitt,OllisSurvey,PasottiDinitzHeffterSurvey,PhamSauermann}):
\begin{conjecture}[Alspach and Liversidge]\label{conj:abelian}
Let $G$ be a finite abelian group. Every $S\subseteq G\setminus\{0\}$ is rearrangeable.
\end{conjecture}

Most previous work towards Conjecture \ref{conj:graham} focused primarily on the two extreme regimes where $|S|$ is either very small or very large. An approach based on the Combinatorial Nullstellensatz~\cite{AlonNullstellensatz} established the conjecture for $|S|\le 10$~\cite{HicksOllisSchmitt} and $|S|\le 12$~\cite{CostaPellegrini}. Direct constructive arguments settled the nearly full cases $|S|\ge p-3$ with nonzero sum~\cite{HicksOllisSchmitt}. Later, as a consequence of their resolution of the random Hall--Paige conjecture, M\"uyesser and Pokrovskiy~\cite{MuyesserPokrovskiy} established that every subset $S\subseteq G\setminus\{\mathrm{id}\}$ of size $|S|\ge (1-o(1))|G|$ is rearrangeable in any finite group $G$. This was subsequently extended to $|S|\ge |G|^{1-c}$ for an absolute constant $c>0$ by Bedert, Buci\'c, Kravitz, Montgomery and M\"uyesser~\cite{BBKMM}. For small sets, Kravitz~\cite{Kravitz} established the conjecture for $|S|\le\log p/\log\log p$ using rectification, also see Sawin~\cite{Sawin}. In a major breakthrough, Bedert and Kravitz~\cite{BedertKravitz} overcame this rectification barrier in \cite{Kravitz,Sawin}, extending the range to $|S|\le \exp(c(\log p)^{1/4})$ for every fixed $c>0$ and sufficiently large $p$. 

Very recently, Pham and Sauermann~\cite{PhamSauermann} bridged the intermediate regime by proving Conjecture~\ref{conj:graham} for $C_\alpha\le |S|\le p^{1-\alpha}$, where $\alpha\in(0,1)$ is arbitrary and $C_\alpha>0$ is a constant depending only on $\alpha$. Together with the dense bound of Bedert, Buci\'c, Kravitz, Montgomery and M\"uyesser~\cite{BBKMM} (by setting $\alpha=c$) and the small-set bounds for $|S|\le C_\alpha$, their result shows a complete resolution of Graham's rearrangement conjecture for all sufficiently large primes $p$.

Nevertheless, proving rearrangeability for general finite abelian groups is significantly more difficult. Our understanding beyond these remains very limited. Bedert, Buci\'c, Kravitz, Montgomery and M\"uyesser~\cite{BBKMM}
proved that every sufficiently large subset of
$\F_2^r\setminus\{0\}$ is rearrangeable. For further partial results in cyclic groups or other groups,
see~\cite{CostaDellaFioreZk,CostaEtAlKneser,CostaDellaFioreFontanaVena}.

\subsection{Rainbow subgraphs}
Finding rainbow subgraphs is an important theme in extremal graph
theory, with connections to additive combinatorics and combinatorial designs~\cite{SudakovRainbowSurvey}. In this paper, we study Conjecture~\ref{conj:abelian} via rainbow subgraphs. In fact, Conjecture~\ref{conj:abelian} can equivalently be stated in terms of rainbow paths in Cayley digraphs. The \emph{colored Cayley digraph} $\Cay(G,S)$ has vertex set $G$ and an edge $x\to x+s$ of color $s$ for each $x\in G$ and $s\in S$. A path or cycle is \emph{rainbow} if its edges have distinct colors. If a directed walk $x_0x_1\cdots x_k$ has successive edge colors $s_1,\ldots,s_k$, then $x_i=x_0+s_1+\cdots+s_i$ for $1\le i\le k$. Thus Conjecture~\ref{conj:abelian} is equivalent to the existence of a rainbow directed path of length $|S|-1$ in every such Cayley digraph (see \cite{BucicEtAl}, and also Lemma~\ref{lem:cayley-walk} in this paper).
Schrijver~\cite{Schrijver} asked the following problem for properly edge-colored regular graphs.

\begin{problem}[Schrijver~\cite{Schrijver}]\label{prob:schrijver}
Does every $d$-regular graph with a proper edge-coloring using $d$ colors contain a rainbow path of length $d-1$?
\end{problem}

Buci\'c, Frederickson, M\"uyesser, Pokrovskiy and Yepremyan~\cite{BucicEtAl} asked the following analogue in directed graphs. An edge-coloring of a digraph is \emph{proper} if edges with a common initial vertex, or a common terminal vertex, have distinct colors.

\begin{problem}[Buci\'c, Frederickson, M\"uyesser, Pokrovskiy and Yepremyan~\cite{BucicEtAl}]\label{prob:directed-schrijver}
Let $H$ be a digraph in which every vertex has in-degree and out-degree $d$, with a proper edge-coloring using $d$ colors. Does $H$ contain a rainbow directed path of length $d-1$?
\end{problem}

An affirmative answer to Problem~\ref{prob:directed-schrijver} would imply the rearrangement conjecture for all finite groups. Buci\'c, Frederickson, M\"uyesser, Pokrovskiy and Yepremyan~\cite{BucicEtAl} proved approximate results for rainbow paths and showed that every finite subset $S$ of a finite group has an ordering in which all but $o(|S|)$ partial products are distinct. Wilson~\cite{Wilson} subsequently proved that the digraphs in Problem~\ref{prob:directed-schrijver} contain rainbow directed paths of length $d-o(d)$, and deduced that every finite subset $S$ of a finite group contains a rearrangeable subset of size $|S|-o(|S|)$. 

These questions belong to a broader study of rainbow embeddings. Alon, Jiang, Miller and Pritikin~\cite{AlonJiangMillerPritikin} studied rainbow subgraphs under local coloring constraints. Montgomery~\cite{MontgomeryRyserBrualdiStein} proved that every Latin square of sufficiently large order $m$ contains a partial transversal of size $m-1$. Alon, Pokrovskiy and Sudakov~\cite{AlonPokrovskiySudakov} found rainbow cycles of length $m-O(m^{3/4})$ in properly edge-colored $K_m$, and recently, Bowtell, Montgomery, M\"uyesser and Pokrovskiy~\cite{BowtellEtAl} proved Andersen's conjecture for large $m$ by finding a rainbow path on $m-1$ vertices. Benzing, Pokrovskiy and Sudakov~\cite{BenzingPokrovskiySudakov} proved that every properly edge-colored complete symmetric digraph on $m$ vertices contains a rainbow directed cycle of length $m-O(m^{4/5})$. Related embedding and absorption methods also underlie the work of Montgomery, Pokrovskiy and Sudakov on rainbow trees and Ringel's conjecture~\cite{MontgomeryPokrovskiySudakovTrees,MontgomeryPokrovskiySudakovRingel}.

\subsection{Main results via dissociated sets}\label{subsec:main-result}

In this paper, we solve Conjecture~\ref{conj:abelian} (also Conjectures \ref{conj:graham} and \ref{conj:Arch}) for sufficiently large sets.

\begin{theorem}\label{thm:asymptotic}
There exists an absolute integer $n_0$ such that the following holds. Let $G$ be a finite abelian group and let $S\subseteq G\setminus\{0\}$ have size $n\ge n_0$. Then $S$ is rearrangeable.
\end{theorem}

By Lemma~\ref{lem:cayley-walk}, this also answers Problem~\ref{prob:directed-schrijver} for the canonical colorings of Cayley digraphs over finite abelian groups whenever $d$ is sufficiently large.

We first give a sufficient condition for rearrangeability in terms of dissociated sets, an additive analogue of linearly independent sets. 
\begin{definition}\label{def:dissociated}
A set $D$ in an abelian group is \emph{dissociated} if every relation $\sum_{d\in D}\varepsilon_d d=0$, with $\varepsilon_d\in\{-1,0,1\}$, has $\varepsilon_d=0$ for every $d\in D$. Equivalently, the subset sums $\sum_{d\in U}d$, where $U\subseteq D$, are pairwise distinct. For $S\subseteq G$, define $\diss(S)=\max\{|D|:D\subseteq S\text{ is dissociated}\}$.
\end{definition}

Dissociated sets play a role analogous to that of linearly independent sets in the study of additive structure~\cite{Rudin,TaoVu}. The maximum size of a dissociated subset of a set is called its \emph{additive dimension}, and controlling this parameter is an important topic in additive combinatorics~\cite{SchoenShkredov,Shkredov}. For example, Chang's theorem bounds the additive dimension of large Fourier spectra~\cite{Chang}, which are used in arguments based on density increments, including Sanders's work on Roth's theorem~\cite{SandersRoth}. For integer intervals, Erd\H{o}s~\cite[Section~6]{ErdosDistinctSubsetSums}
asked whether $\diss(\{1,\ldots,N\})\le\log_2 N+O(1)$. Dubroff, Fox and Xu~\cite{DubroffFoxXu} presented two short proofs yielding the weaker upper bound $\diss(\{1,\ldots,N\})\le\log_2 N+\frac12\log_2\log_2 N+O(1)$.

Dissociated sets also give properly edge-colored graphs with no rainbow cycle~\cite{AlonEtAlRainbowCycles,KuangWangRainbowSubdivisions}. Given $D\subseteq G$, let $X$ and $Y$ be two disjoint copies of $G$. Join $x\in X$ to $y\in Y$ whenever $y-x\in D$, and color the edge $xy$ with $y-x$. This defines a properly edge-colored bipartite graph. A rainbow cycle with successive edge colors $d_1,\ldots,d_{2\ell}$ yields $d_1-d_2+\cdots+d_{2\ell-1}-d_{2\ell}=0$, a contradiction. A set $D\subseteq G$ is called a \emph{Sidon set} if every equality $a+b=c+d$ with $a,b,c,d\in D$ implies $\{a,b\}=\{c,d\}$~\cite{ErdosTuranSidon}. When $D$ is Sidon, the same construction contains no rainbow $4$-cycle, since such a cycle would yield $d_1+d_3=d_2+d_4$ with four distinct elements of $D$.

For rearrangement problems, dissociativity is useful because it yields anticoncentration estimates for random partial sums. Let $D$ be a dissociated set of size $r$, and let $j$ be an integer with $0\le j\le r$. If $U$ is a uniformly random $j$-element subset of $D$, then $\mathbb P(\sum_{u\in U}u=z)\le\binom{r}{j}^{-1}$ for every $z\in G$. We use this estimate to insert suitably chosen blocks of $S\setminus D$ into a random ordering of $D$ while avoiding zero-sum intervals. This gives the following sufficient condition.

\begin{theorem}\label{thm:main}
There is an absolute constant $A>0$ such that the following holds. Let $G$ be a finite abelian group and let $S\subseteq G\setminus\{0\}$ have size $n\ge3$. If
\begin{equation}\label{dis}
\diss(S)\ge A\frac{(\log n)^2}{\log\log n},
\end{equation}
then $\Cay(G,S)$ contains a rainbow directed path of length $n$ using every color in $S$ when $\sum_{s\in S}s\ne0$, and a rainbow directed cycle of length $n$ using every color in $S$ when $\sum_{s\in S}s=0$. In particular, $S$ is rearrangeable.
\end{theorem}

Every $n$-element set $S\subseteq G\setminus\{0\}$ satisfies $\diss(S)\ge\log_3(n+1)$. The dimension threshold in~\eqref{dis} is within a factor $O(\log n/\log\log n)$ of this bound. Indeed, a maximal dissociated set $D\subseteq S$ satisfies $S\subseteq\{\sum_{d\in D}\varepsilon_dd:\varepsilon_d\in\{-1,0,1\}\}$, so $n+1\le3^{|D|}$. Equality holds for $S=\F_3^r\setminus\{0\}$, where dissociativity is equivalent to linear independence.

To handle sets of small additive dimension, we study the intersections of $S$ with cosets of subgroups. We show that each such intersection can be recovered from a small subset of its elements, yielding a bound on the number of distinct intersections in terms of $|S|$ and $\diss(S)$. This leads to the following structural result. The definitions can be seen in Section \ref{sec:prelim}.

\begin{theorem}[Algebraic coset containers]\label{lem:coset-traces}
Let $S\subseteq G\setminus\{0\}$ have size $n\ge3$, and let $d=\diss(S)$ and $M=\lceil4d\log_2(4d+1)\rceil$. For nonempty $F\subseteq S$, define $\mathcal C_S(F)=S\cap\aff(F)$, and let $\mathcal C_S(\varnothing)=\varnothing$. Then the following properties hold.
\begin{enumerate}[label=\textup{(\roman*)}]
\item For every $a\in S$, we have $\diss(S-a)\le M$.
\item Every coset trace $T\in\mathcal T(S)$ has a \emph{fingerprint} $F\subseteq T$ with $|F|\le M+1$ and $T=\mathcal C_S(F)$. If $T\ne\varnothing$, we can moreover choose $F$ with $\aff(F)=\aff(T)$. Thus,
\begin{equation}\label{eq:algebraic-container-count}
|\mathcal T(S)|\le\sum_{j=0}^{\min\{M+1,n\}}\binom nj\le n^{M+2}.
\end{equation}
The set system $\mathcal T(S)$ has VC-dimension at most $M+1$.
\item For every $0<\kappa<n/(M+1)$, there is a nonempty $F\subseteq S$ of size at most $M+1$ such that, with $A=\mathcal C_S(F)$ and $H=\langle F-F\rangle$,
\begin{equation}\label{eq:algebraic-stable-container}
|S\setminus A|\le M\kappa\quad\text{and}\quad
|A\setminus(x+K)|\ge\kappa
\end{equation}
for every proper subgroup $K<H$ and every $x\in G$.
\end{enumerate}
\end{theorem}
For graphs, the use of small fingerprints to describe independent sets goes back to the Kleitman--Winston algorithm; see Samotij~\cite[Section~2]{SamotijContainers}. The counting argument shares the idea of encoding sets by small fingerprints with the hypergraph container method, developed independently by Balogh, Morris and Samotij~\cite{BMSContainers} and Saxton and Thomason~\cite{STContainers}. Under suitable hypotheses, every independent set $I$ is contained in a container $C(F)$ determined by a small fingerprint $F\subseteq I$. Counting the possible fingerprints controls the number of containers; together with bounds on the sets inside each container, this yields enumeration results, including bounds on the number of Sidon sets~\cite{SaxtonThomasonOnline}. Balogh, Das, Delcourt, Liu and Sharifzadeh~\cite{BaloghDasEtAl} similarly bound the number of maximal intersecting families by counting small generating subfamilies. In our setting, each coset trace $T$ is recovered exactly from a small subset $F\subseteq T$ through the identity $T=S\cap\aff(F)$. Counting these fingerprints bounds $|\mathcal T(S)|$ in terms of $|S|$ and $\diss(S)$, independently of the number of subgroups of $G$.

\subsection{Proof overview}
\subsubsection{Theorem \ref{thm:main}}
Let $D\subseteq S$ be a dissociated set of size $r$, and write $R=S\setminus D$. If $R=\varnothing$, then we are done. Otherwise, consider a uniformly random ordering $\pi=(\pi_1,\ldots,\pi_r)$ of $D$. For every $J\subseteq[r]$ with $|J|=j$ and every $z\in G$, dissociativity implies $$\mathbb P\bigl(\sum_{i\in J}\pi_i=z\bigr)\le\binom{r}{j}^{-1}.$$

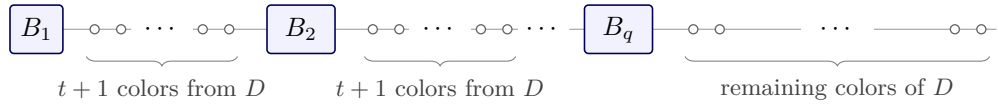
\begin{figure}[!htb]
\centering
\begin{tikzpicture}
[
    x=1cm,
    y=1cm,
    font=\small,
    block/.style={
        draw=blue!35!black,
        fill=blue!5,
        rounded corners=1.5pt,
        line width=.6pt,
        minimum height=.64cm,
        minimum width=.9cm,
        inner sep=0pt
    },
    colorpoint/.style={
        circle,
        draw=black!60,
        fill=white,
        line width=.5pt,
        minimum size=3.4pt,
        inner sep=0pt
    },
    brace/.style={
        decorate,
        decoration={brace,mirror,amplitude=3.5pt},
        draw=black!50,
        line width=.4pt
    },
    note/.style={
        font=\footnotesize,
        text=black!75
    }
]
\path[use as bounding box] (-.1,-1.14) rectangle (13.25,.44);
\draw[black!35,line width=.45pt] (.45,0)--(6.72,0);
\draw[black!35,line width=.45pt] (7.40,0)--(13.05,0);
\node[block,minimum width=.7cm] at (.35,0) {$B_1$};
\node[block] at (3.85,0) {$B_2$};
\node[block] at (8.05,0) {$B_q$};
\foreach \x in {1.12,1.47,2.55,2.90,4.80,5.15,6.23,6.58,9.04,9.39,12.50,12.85}
    \node[colorpoint] at (\x,0) {};
\node[fill=white,inner xsep=5pt,inner ysep=2pt] at (2.01,0) {$\cdots$};
\node[fill=white,inner xsep=5pt,inner ysep=2pt] at (5.69,0) {$\cdots$};
\node at (7.06,0) {$\cdots$};
\node[fill=white,inner xsep=7pt,inner ysep=2pt] at (10.95,0) {$\cdots$};
\draw[brace] (1.01,-.36)--(3.01,-.36);
\node[note] at (2.01,-.79) {$t+1$ colors from $D$};
\draw[brace] (4.69,-.36)--(6.69,-.36);
\node[note] at (5.69,-.79) {$t+1$ colors from $D$};
\draw[brace] (8.93,-.36)--(12.96,-.36);
\node[note] at (10.945,-.79) {remaining colors of $D$};
\end{tikzpicture}
\caption{The blocks $B_1,\ldots,B_q$ of $R$ are inserted into a uniformly random ordering of $D$.}
\label{fig:main-splicing}
\end{figure}

By Lemma~\ref{lem:path-decomposition}, we divide $R$ into ordered blocks $B_1,\ldots,B_q$, where $q=O(\log(en))$ and no block contains a nonempty zero-sum interval. We insert these blocks into $\pi$, with $t+1$ elements of $D$ between consecutive blocks (Figure~\ref{fig:main-splicing}).

Let $Q_i$ be the set of sums of the nonempty suffixes of $B_i$. If a suffix of $B_i$ together with the next $j<t$ elements of $D$ had sum zero, the sum of these $j$ elements would belong to $-Q_i$. By symmetry of the uniformly random ordering, the sum of these $j$ elements has the same distribution as $\pi_1+\cdots+\pi_j$. For a fixed small constant $\tau>0$, by Lemma~\ref{lem:path-decomposition}, we have
$
\sum_{j=1}^{t-1}e^{\tau j}\Prob\bigl(\pi_1+\cdots+\pi_j\in-Q_i\bigr)\le\frac1{100}.
$
As $e^{\tau j}\ge1$, the union bound yields
\[
\Prob\bigl(\exists\,1\le j<t:\ \pi_1+\cdots+\pi_j\in-Q_i\bigr)\le\frac1{100}.
\]
The analogous estimate holds for prefixes of $B_i$ together with preceding elements of $D$. Thus, for each $i<q$, the two weighted bounds at the ends of the segment of $D$ between $B_i$ and $B_{i+1}$ sum to at most $1/50$. By Lemma~\ref{lem:conditioning}, these local zero-sum intervals can therefore be excluded simultaneously, with a controlled loss in the remaining subset-sum estimates.

Every interval containing fewer than $t$ elements of $D$ meets at most one block. After excluding the intervals above, a short zero-sum interval must contain an entire block and elements of $D$ on both sides. For each $j<t$, there are at most $q(j-1)$ such intervals, and the conditional subset-sum estimates control their probabilities. When $t\le j\le r-t$, the bound $\binom rj^{-1}$ suffices to handle all possible endpoints. When $r-t<j<r$, we instead consider the few omitted elements of $D$. Finally, the first block is chosen as a singleton $B_1=(x)$ with $x\ne\sum_{s\in S}s$, which excludes proper zero-sum intervals containing all of $D$.

With suitable choices of $t=\Theta(\log n/\log\log n)$ and $r=\Theta(t\log n)=\Theta((\log n)^2/\log\log n)$, the probability that some proper nonempty interval has sum zero is less than one; see Theorem~\ref{thm:splicing}.

\subsubsection{Theorem \ref{thm:asymptotic}}
By Theorem~\ref{thm:main}, it suffices to consider $\diss(S)=O((\log n)^2)$. By Theorem~\ref{lem:coset-traces}, almost all elements of $S$ lie in a coset $a+H$, while $S\cap(a+H)$ has many elements outside every coset of a proper subgroup of $H$. We distinguish three cases according to $|H|$ and $\ord_{G/H}(a+H)$. The following informal statements correspond to these cases; their precise statements appear in Lemma~\ref{lem:large-subgroup-core}, Theorem~\ref{thm:large-order-absorption} and Theorem~\ref{lem:dense-core-sparse-outside}, respectively.

\begin{lemma}[A large direction subgroup, informal]\label{lem:informal-large-subgroup}
Suppose that $\diss(S)=O((\log n)^2)$ and all but a sufficiently small power of $n$ elements of $S$ lie in a coset $a+H$. If $|H|\ge n^{1+\delta}$ for some fixed $\delta>0$ and $S\cap(a+H)$ has at least a fixed positive power of $n$ elements outside every proper subcoset, then $S$ is rearrangeable.
\end{lemma}

\begin{lemma}[A coset of large order, informal]\label{lem:informal-large-order}
Suppose that all but a small power of $n$ elements of $S$ lie in a coset $a+H$. If the order of $a+H$ in $G/H$ is at least a larger fixed power of $n$, then $S$ is rearrangeable.
\end{lemma}

\begin{lemma}[A dense subgroup core, informal]\label{lem:informal-dense-core}
Let $Q\le G$ be a sufficiently large subgroup. If $S\cap Q$ is sufficiently dense in $Q$ and $|S\setminus Q|$ is sufficiently small relative to $|Q|$, then $S$ is rearrangeable.
\end{lemma}

The first two cases, covered by Lemmas~\ref{lem:informal-large-subgroup} and~\ref{lem:informal-large-order}, use local repair, adapting the methods of Pham and Sauermann~\cite{PhamSauermann} and Costa, Della Fiore, Feng and Liu~\cite{CostaEtAlKneser}. In the first case, the bounds outside proper subcosets yield anticoncentration estimates for random subset sums. In the second case, the large order of $a+H$ restricts the possible lengths of zero-sum intervals in an ordering of $S\cap(a+H)$. These provide two different ways to control the events arising in local repair. In the third case, both $|H|$ and $q=\ord_{G/H}(a+H)$ are small, so $Q=H+\langle a\rangle$ has order $q|H|$ and contains almost all of $S$. The set $S\cap Q$ is then sufficiently dense for Lemma~\ref{lem:informal-dense-core}. Here we combine the expansion approach of Buci\'c, Frederickson, M\"uyesser, Pokrovskiy and Yepremyan~\cite{BucicEtAl} with the absorption method of Bedert, Buci\'c, Kravitz, Montgomery and M\"uyesser~\cite{BBKMM}. The full deduction, including the choice of constants, appears in Section~\ref{sec:dense-core}.

\paragraph{Organization.}
Section~\ref{sec:prelim} introduces notation and some basic results. In Section~\ref{sec:mixing}, we show estimates for subset sums arising from random orderings of dissociated sets. We use these estimates in Section~\ref{sec:blocks} to prove Theorem~\ref{thm:main}. Section~\ref{sec:concentrated-core} develops the coset concentration and slice anticoncentration estimates. Sections~\ref{sec:large-subgroup} and~\ref{sec:large-order-coset} handle the cases of a large direction subgroup and a coset of large order in the quotient, respectively. In Section~\ref{sec:dense-core}, we settle the remaining case using expansion and absorption and complete the proof of Theorem~\ref{thm:asymptotic}. The proof of Lemma~\ref{lem:nested-to-local-remote} appears in Appendix~\ref{app:local-repair}.

\section{Preliminaries}\label{sec:prelim}
Throughout, $G$ is a finite abelian group, written additively, with identity $0$. All logarithms are natural unless a base is specified. For integers $a,b$, write $[a,b]=\{j\in\mathbb Z:a\le j\le b\}$, and let $[m]=[1,m]$ for $m\ge0$. For a finite set $X$, let $\binom{X}{j}$ denote the family of its $j$-element subsets. We write $X\mathbin{\dot\cup}Y$ for a disjoint union and $X\triangle Y=(X\setminus Y)\cup(Y\setminus X)$ for the \emph{symmetric difference}.

For $X,Y\subseteq G$ and $a\in G$, write $X\pm Y=\{x\pm y:x\in X,\ y\in Y\}$ and $a\pm X=\{a\pm x:x\in X\}$. For an integer $r\ge1$, let $rX=\{x_1+\cdots+x_r:x_1,\ldots,x_r\in X\}$, where repetitions are allowed. We write $\sigma(X)=\sum_{x\in X}x$, with $\sigma(\varnothing)=0$, and let $\langle X\rangle$ denote the \emph{subgroup generated by $X$}. The \emph{span} of $X$ is $\Span(X)=\{\sum_{x\in X}\varepsilon_xx:\varepsilon_x\in\{-1,0,1\}\}$.

For a nonempty set $X\subseteq G$, let $\aff(X)$ denote the smallest coset containing $X$. Thus $\aff(X)=a+\langle X-X\rangle=a+\langle X-a\rangle$ for any $a\in X$, and we call $\langle X-X\rangle$ its \emph{direction subgroup}. A proper subcoset of $a+H$ is a coset $b+K$ with $b\in a+H$ and $K<H$. We write $\mathcal T(X)=\{X\cap(a+H):a\in G,\ H\le G\}$ for the family of \emph{coset traces} of $X$. For $H\le G$ and $a\in G$, let $\ord(a+H)=\ord_{G/H}(a+H)=\min\{q\ge1:qa\in H\}$.

The Vapnik--Chervonenkis dimension (\emph{VC-dimension} for short), introduced by Vapnik and Chervonenkis~\cite{VapnikChervonenkis}, measures the complexity of a set system. By the Sauer and Shelah lemma~\cite{Sauer,Shelah}, set systems of bounded VC-dimension have polynomially many traces on finite subsets of their ground sets. This parameter also underlies the $\varepsilon$-net theorem of Haussler and Welzl~\cite{HausslerWelzl}.

A \emph{set system} on a finite set $X$ is a family $\mathcal F\subseteq2^X$. For $Y\subseteq X$, its \emph{trace} on $Y$ is $\mathcal F|_Y=\{F\cap Y:F\in\mathcal F\}$. We say that $\mathcal F$ \emph{shatters} $Y$ if $\mathcal F|_Y=2^Y$. For $\mathcal F\ne\varnothing$, its VC-dimension is $\operatorname{VC}(\mathcal F)=\max\{|Y|:Y\subseteq X,\ \mathcal F|_Y=2^Y\}$.

For a digraph $D$, let $V(D)$ and $E(D)$ denote its vertex and edge sets. Oppositely directed edges between the same two vertices are allowed. For an edge-colored digraph, let $\operatorname{col}(D)$ denote its set of edge colors. The coloring is \emph{proper} if edges with a common initial vertex, or a common terminal vertex, have distinct colors. A subgraph is \emph{rainbow} if its edges have pairwise distinct colors. All paths and cycles in digraphs are directed, and their length is their number of edges. For $S\subseteq G\setminus\{0\}$, the \emph{Cayley digraph} $\Cay(G,S)$ has vertex set $G$ and an edge $x\to x+s$ of color $s$ for each $x\in G$ and $s\in S$. We write $x\xrightarrow{s}y$ when $y=x+s$. The \emph{underlying undirected graph of \(\operatorname{Cay}(G,B)\)} is the graph with vertex set \(G\) in which \(x\) and \(y\) are adjacent if \(y-x\in B\cup(-B)\).

For $X,Y\subseteq V(D)$, let $e_D(X,Y)$ denote the number of edges directed from $X$ to $Y$. We omit the subscript when the digraph is clear. We say that $D$ has \emph{no $\eta$-sparse cut} if $e_D(X,V(D)\setminus X)\ge\eta|X||V(D)\setminus X|$ for every $X\subseteq V(D)$.

\subsection{Orderings and Cayley walks}
For an ordered tuple $B=(b_1,\ldots,b_\ell)$, write $|B|=\ell$, $\sigma(B)=\sigma(b_1,\ldots,b_\ell)=b_1+\cdots+b_\ell$, $V(B)=\{\sigma(b_1,\ldots,b_i):i\in[\ell]\}$, and $\overleftarrow B=(b_\ell,\ldots,b_1)$. An \emph{interval} of an ordered tuple $B=(b_1,\ldots,b_\ell)$ is a subtuple $(b_i,\ldots,b_j)$, where $1\le i\le j\le\ell$. It is \emph{proper} if $(i,j)\ne(1,\ell)$, and \emph{zero-sum} if its entries sum to $0$.

An ordering $B$ of distinct nonzero elements is \emph{valid} if $|V(B)|=\ell$, and is a \emph{path ordering} if, in addition, $0\notin V(B)$. Thus, for a path ordering, $0,b_1,b_1+b_2,\ldots,b_1+\cdots+b_\ell$ are pairwise distinct.

The next lemma shows the connection between rainbow subgraphs in Cayley graphs and rearrangeability.

\begin{lemma}\label{lem:cayley-walk}
Let \(S\subseteq G\setminus\{0\}\) and let \(m\ge1\). Let
\(P= v_0\xrightarrow{b_1}v_1\xrightarrow{b_2}\cdots \xrightarrow{b_m}v_m\)
be a directed walk in \(\Cay(G,S)\), where
\(b_1,\ldots,b_m\) are distinct.
Then the following hold.
\begin{enumerate}[label=\textup{(\roman*)}]
    \item For \(0\le i<j\le m\), we have \(v_i=v_j\) if and only if \(\sigma(b_{i+1},\ldots,b_j)=0\). Thus, \(P\) is a
    rainbow directed path if and only if no nonempty interval of
    \((b_1,\ldots,b_m)\) has sum zero.

    \item If \(\sigma(b_1,\ldots,b_m)=0\), then \(P\) is a rainbow
    directed cycle if and only if no proper nonempty interval of
    \((b_1,\ldots,b_m)\) has sum zero.

    \item A nonempty finite set \(S\subseteq G\setminus\{0\}\) is
    rearrangeable if and only if \(\Cay(G,S)\) contains a rainbow
    directed path of length \(|S|-1\).
\end{enumerate}
\end{lemma}
\begin{proof}
For $0\le i<j\le m$, we have $v_j-v_i=\sigma(b_{i+1},\ldots,b_j)$. This proves~\textup{(i)} and~\textup{(ii)}, since a directed cycle allows only the coincidence $v_0=v_m$.

For~\textup{(iii)}, the nonempty partial sums of a valid ordering $s_1,\ldots,s_n$ form a rainbow path with colors $s_2,\ldots,s_n$. Conversely, given a rainbow path $v_0\xrightarrow{b_1}\cdots\xrightarrow{b_{n-1}}v_{n-1}$, let $a$ be its unique unused color. The partial sums of $a,b_1,\ldots,b_{n-1}$ are the distinct elements $a+v_j-v_0$, $0\le j<n$.
\end{proof}

\subsection{Fourier analysis}
We use standard Fourier notation on finite abelian groups; see, for example, \cite{GreenFiniteField,TaoVu}. Let $\widehat G$ be the character group of $G$, so each $\chi\in\widehat G$ satisfies $\chi(x+y)=\chi(x)\chi(y)$ and $|\chi(x)|=1$. For $f:G\to\mathbb C$, write $\widehat f(\chi)=\sum_{x\in G}f(x)\overline{\chi(x)}$ and $\|f\|_2^2=\sum_{x\in G}|f(x)|^2$. With this normalization, Fourier inversion takes the form $f(x)=|G|^{-1}\sum_{\chi\in\widehat G}\widehat f(\chi)\chi(x)$, and Parseval's identity states that $\|f\|_2^2=|G|^{-1}\sum_{\chi\in\widehat G}|\widehat f(\chi)|^2$.

\subsection{Probabilistic lemmas}\label{sec:probabilistic-inequalities}
We use $\Prob$ and $\mathbb E$ for probability and expectation, respectively. And we write $\mathbf1_P$ for the indicator function of a set $P$,
and $\mathbf1_E$ for the indicator random variable of an event $E$. For $0\le p\le1$, a $p$-random subset of a finite set is obtained by including each element independently with probability $p$. We shall use the following standard inequalities below.

\begin{lemma}[Chernoff's inequality~\cite{AlonSpencer}]\label{lem:chernoff}
Let $X$ be a sum of independent Bernoulli random variables, with $\mu=\mathbb EX>0$. Then, for every $t>0$,
\begin{enumerate}[label=\textup{(\roman*)}]
\item $\Prob(X\le\mu-t)\le\exp(-t^2/(2\mu))$;
\item $\Prob(X\ge\mu+t)\le\exp(-t^2/(2\mu+t))$.
\end{enumerate}
\end{lemma}

The following lemma was proved by Hoeffding \cite{Hoeffding}.
\begin{lemma}[Chernoff's inequality for the hypergeometric distribution \cite{Hoeffding}]
\label{lem:hypergeometric}
Let $\Omega$ be a set of size $N\ge1$, let $A\subseteq\Omega$ have size $K$, and let $R$ be a uniformly random $m$-subset of $\Omega$. Write $X=|R\cap A|$ and $\mu=\mathbb EX=mK/N>0$. Then, for $0<\delta<1$,
\begin{enumerate}[label=\textup{(\roman*)}]
\item $\Prob(X\le(1-\delta)\mu)\le\exp(-\delta^2\mu/2)$;
\item $\Prob(X\ge(1+\delta)\mu)\le\exp(-\delta^2\mu/3)$.
\end{enumerate}
\end{lemma}

\section{Random orderings of dissociated sets}\label{sec:mixing}
In this section, we establish the estimates needed in Section~\ref{sec:blocks}. We first analyze the subset sums of a dissociated set, and subsequently prove a conditioning lemma and a key mixing estimate.

\subsection{Subset sums of dissociated sets}
For a finite set $D\subseteq G$ and $0\le j\le|D|$, write $\Sigma_j(D)=\{\sigma(U):U\in\binom Dj\}$. 
\begin{lemma}\label{lem:boolean-cube}\label{3.1}
Let $D\subseteq G$ be dissociated with $|D|=r$. Then the sets $\Sigma_0(D),\ldots,\Sigma_r(D)$ are pairwise disjoint, and $|\Sigma_j(D)|=\binom rj$ for every $0\le j\le r$.
\end{lemma}
\begin{proof}
Since $D$ is dissociated, the sums $\sigma(U)$, $U\subseteq D$, are pairwise distinct. Both assertions follow.
\end{proof}

Fix a dissociated set \(D\) of size \(r\). For \(A\subseteq G\)
and \(0\le j\le r\), define
\(\mu_j^D(A)=|A\cap\Sigma_j(D)|/\binom rj\).
If \(\pi=(\pi_1,\ldots,\pi_r)\) is a uniformly random ordering
of \(D\), then
\(\mu_j^D(A)=\Prob(\pi_1+\cdots+\pi_j\in A)\), since its
first \(j\) elements form a uniformly random \(j\)-subset of \(D\).

For \(1\le t\le r\), \(\tau\ge0\), \(z\in G\), and \(P\subseteq G\),
define
\[
   \Lambda_{t,\tau}^D(z,P)
   =\sum_{j=1}^{t-1}e^{\tau j}\mu_j^D(z-P).
\]

For a uniformly random ordering $\pi$ of $D$, we have
$\Lambda_{t,\tau}^D(z,\{x\})=\sum_{j=1}^{t-1}e^{\tau j}\Prob(\pi_1+\cdots+\pi_j=z-x)$. Since $D$ is dissociated, the walk starting at $x$ with successive increments $\pi_1,\ldots,\pi_r$ visits each vertex at most once. Thus, when $\tau=0$, this quantity is exactly the probability that the walk reaches $z$ exactly at one of the times $1,\ldots,t-1$. For $\tau>0$, it is a weighted sum of these hitting probabilities, with a visit at time $j$ receiving weight $e^{\tau j}$. Note that the weight $e^{\tau j}$ indeed absorbs the loss from conditioning. When \(D,t,\tau\) are fixed, we write \(\Lambda(z,P)=\Lambda_{t,\tau}^D(z,P)\). We have the following easy fact which will be used many times.
\begin{fact}
\label{fact:subadditivity}
For \(P_1,P_2\subseteq G\), we have \(\mu_j^D(z-(P_1\cup P_2)) \le\mu_j^D(z-P_1)+\mu_j^D(z-P_2)\).
\end{fact}

We will use the following lemma to control the effect of adding a single vertex to $P$.
\begin{lemma}\label{lem:one-vertex}
Suppose that \(r\ge 4t\) and \(0\le \tau\le \ln 2\). Then, for all \(x,z\in G\), $\Lambda_{t,\tau}^D(z,\{x\})\le2/r$.
\end{lemma}

\begin{proof}
If $t=1$, then $\Lambda_{t,\tau}^D(z,\{x\})=0$. Assume that $t\ge2$. By Lemma \ref{3.1}, there is at most one $h\in\{1,\ldots,t-1\}$ such that $z-x\in\Sigma_h(D)$. If no such $h$ exists, then $\Lambda_{t,\tau}^D(z,\{x\})=0$. Otherwise, $\Lambda_{t,\tau}^D(z,\{x\})=e^{\tau h}\binom rh^{-1}$. Since $e^{\tau h}\binom rh^{-1}$ decreases with $h$ for $1\le h<t$, we have $\Lambda_{t,\tau}^D(z,\{x\})\le e^\tau/r\le2/r$.
\end{proof}

\subsection{Subset sums under local conditioning}
The following lemma controls subset sums under local restrictions on a random ordering of $D$.

\begin{lemma}\label{lem:conditioning}
Let $\pi=(\pi_1,\ldots,\pi_r)$ be a uniformly random ordering of a dissociated set $D\subseteq G$. Let $\mathcal U$ be a family of pairwise disjoint subsets of $[r]$, and write $M=|\bigcup_{U\in\mathcal U}U|\le r/100$.

For each $U\in\mathcal U$, let $\mathcal C(U)$ be a finite family of pairs $(J,A)$ with $\varnothing\ne J\subseteq U$ and $A\subseteq G$. Define
\[
\mathcal G_U=\left\{\sum_{i\in J}\pi_i\notin A\text{ for every }(J,A)\in\mathcal C(U)\right\},
\qquad
\mathcal G=\bigcap_{U\in\mathcal U}\mathcal G_U.
\]
Suppose that $0<\rho\le1/2$, $\tau\ge4M/r$, and
\begin{equation}\label{e1}
\sum_{(J,A)\in\mathcal C(U)}e^{\tau|J|}\mu_{|J|}^D(A)\le\rho
\qquad\text{for every }U\in\mathcal U.
\end{equation}
Then the following hold.
\begin{enumerate}[label=\textup{(\roman*)}]
\item $\Prob(\mathcal G)\ge e^{-4\rho|\mathcal U|}$.
\item For every $I\subseteq[r]$ with $|I|=s\le r/20$ and every $B\subseteq G$,
\[
\Prob\left(\sum_{i\in I}\pi_i\in B\,\middle|\,\mathcal G\right)
\le e^{4Ms/r+4\rho m_{\mathcal U}(I)}\mu_s^D(B),
\]
where $m_{\mathcal U}(I)$ is the number of members of $\mathcal U$ meeting $I$.
\item For every $I\subseteq[r]$ with $|I|=s$ and every $y\in G$,
\[
\Prob\left(\sum_{i\in I}\pi_i=y\,\middle|\,\mathcal G\right)
\le e^{4\rho|\mathcal U|}\binom rs^{-1}.
\]
\end{enumerate}
\end{lemma}
\begin{proof}
For integers $0\le m\le M$ and $0\le j\le r/20$, we have
\[
\frac{\binom rj}{\binom{r-m}{j}}
=\prod_{k=0}^{j-1}\left(1+\frac{m}{r-m-k}\right)
\le\exp\left(\frac{mj}{r-m-j+1}\right)
\le e^{4Mj/r}.
\]

Let $W$ be a union of members of $\mathcal U$. After exposing $(\pi_i)_{i\in W}$, the remaining entries form a uniformly random ordering of $D\setminus\{\pi_i:i\in W\}$. Since $|J|\le M\le r/100$ for every $(J,A)\in\mathcal C(U)$, by the union bound and \eqref{e1}, we have
\[
\begin{aligned}
\Prob(\mathcal G_U^c\mid\pi|_W)
&\le\sum_{(J,A)\in\mathcal C(U)}
\frac{\binom r{|J|}}{\binom{r-|W|}{|J|}}\mu_{|J|}^D(A)\le\sum_{(J,A)\in\mathcal C(U)}
e^{\tau|J|}\mu_{|J|}^D(A)
\le\rho.
\end{aligned}
\]
Since each $\mathcal G_U$ holds with conditional probability at least $1-\rho$, we have $\Prob(\mathcal G)\ge(1-\rho)^{|\mathcal U|}\ge e^{-4\rho|\mathcal U|}$, which proves \textup{(i)}.

For \textup{(ii)}, fix $I\subseteq[r]$ with $|I|=s\le r/20$ and $B\subseteq G$. Set $\mathcal U_0=\{U\in\mathcal U:U\cap I=\varnothing\}$, $W=\bigcup_{U\in\mathcal U_0}U$, $\mathcal G_0=\bigcap_{U\in\mathcal U_0}\mathcal G_U$, and $E=\{\sum_{i\in I}\pi_i\in B\}$. Since $W\cap I=\varnothing$ and $\mathcal G_0$ is determined by $\pi|_W$, we have
\[
\begin{aligned}
\Prob(E\mid\mathcal G_0)
&=\mathbb E\bigl[\Prob(E\mid\pi|_W)\mid\mathcal G_0\bigr]\le\frac{\binom rs}{\binom{r-|W|}{s}}\mu_s^D(B)
\le e^{4Ms/r}\mu_s^D(B).
\end{aligned}
\]
Since each of the $m_{\mathcal U}(I)$ remaining events $\mathcal G_U$ holds with conditional probability at least $1-\rho$, we have $\Prob(\mathcal G\mid\mathcal G_0)\ge(1-\rho)^{m_{\mathcal U}(I)}\ge e^{-4\rho m_{\mathcal U}(I)}$. Since $\mathcal G\subseteq\mathcal G_0$, we obtain
\[
\Prob(E\mid\mathcal G)
\le\frac{\Prob(E\mid\mathcal G_0)}{\Prob(\mathcal G\mid\mathcal G_0)}
\le e^{4Ms/r+4\rho m_{\mathcal U}(I)}\mu_s^D(B).
\]

For \textup{(iii)}, since $D$ is dissociated, at most one $s$-subset of $D$ has sum $y$. Thus, by \textup{(i)},
\[
\Prob\left(\sum_{i\in I}\pi_i=y\,\middle|\,\mathcal G\right)
\le\frac{\Prob(\sum_{i\in I}\pi_i=y)}{\Prob(\mathcal G)}
\le e^{4\rho|\mathcal U|}\binom rs^{-1}.
\]
\end{proof}

The next lemma shows that \(\Lambda_{t,\tau}^D(z,P)\) can be large for only a small number of \(z\in G\).

\begin{lemma}\label{lem:cube-mixing}
Let \(D\subseteq G\) be dissociated with \(|D|=r\), and suppose that
\(t\ge3\), \(r\ge10^{12}t\), \(0\le\tau\le1/8000\), and
\(|P|\le e^{r/100}\). Then, for every \(\theta>0\),
$
\left|\{z\in G:\Lambda_{t,\tau}^D(z,P)>\theta\}\right|
\le
10^{12}\theta^{-2}|P|.
$
\end{lemma}
\begin{proof}
It suffices to prove $\sum_{z\in G}\Lambda_{t,\tau}^D(z,P)^2\le10^{12}|P|$. We may assume $P\ne\varnothing$.

Define $T_jf(z)=\binom rj^{-1}\sum_{u\in\Sigma_j(D)}f(z-u)$ and $T=\sum_{j=1}^{t-1}e^{\tau j}T_j$, so $(T\mathbf1_P)(z)=\Lambda_{t,\tau}^D(z,P)$. For $\chi\in\widehat G$, we have $T_j\chi=\lambda_j(\chi)\chi$, where $\lambda_j(\chi)=e_j((\overline{\chi(d)})_{d\in D})/\binom rj$ and $e_j$ is the $j$th elementary symmetric polynomial.

Let $a=999/1000$ and $\Omega=\{\chi:|\sum_{d\in D}\chi(d)|>ar\}$. We first prove that $|\lambda_j(\chi)|\le100e^{-j/4000}$ whenever $\chi\notin\Omega$ and $1\le j<t$. For such $\chi$ and $j$, Cauchy's coefficient estimate~\cite[Chapter~4]{AhlforsComplex} and the arithmetic-geometric mean inequality imply $|e_j((\overline{\chi(d)}))|\le\varrho^{-j}(1+\varrho^2+2a\varrho)^{r/2}\le(ar/j)^j e^{j+j^2/(2a^2r)}$, where $\varrho=j/(ar)$. Since $j!\le3\sqrt j(j/e)^j$ and $\log(1-u)\ge-2u$ for $0\le u\le1/2$, we have $\binom rj=\frac{r^j}{j!}\prod_{k=0}^{j-1}(1-k/r)\ge\frac1{3\sqrt j}(er/j)^j e^{-j^2/r}$. Consequently, $|\lambda_j(\chi)|\le3\sqrt j\,a^j e^{(1+1/(2a^2))j^2/r}\le3\sqrt j\,e^{-(1/1000-2/10^{12})j}\le100e^{-j/4000}$.

For $\omega\in\{1,-1,i,-i\}$, consider the Riesz product $Q_\omega(\chi)=\prod_{d\in D}(1+\tfrac12\Re(\omega\chi(d)))$; see~\cite{Rudin,TaoVu}. Since $D$ is dissociated, we have $|G|^{-1}\sum_{\chi\in\widehat G}Q_\omega(\chi)=1$ by character orthogonality. Each $\chi\in\Omega$ satisfies $\Re(\omega\sum_d\chi(d))>ar/\sqrt2>7r/10$ for some $\omega$. Since $\log(1+x/2)\ge\tfrac{1+x}{2}\log3-\log2$ for $-1\le x\le1$, we have $\log Q_\omega(\chi)>(\tfrac{17}{20}\log3-\log2)r>6r/25$. Hence $|\Omega|/|G|\le4e^{-6r/25}\le e^{-r/10}$.

The eigenvalue $\lambda(\chi)=\sum_{j=1}^{t-1}e^{\tau j}\lambda_j(\chi)$ of $T$ satisfies $|\lambda(\chi)|\le100\sum_{j\ge1}e^{-j/8000}\le8\cdot10^5$ outside $\Omega$, and $|\lambda(\chi)|\le te^{\tau t}$ everywhere. By Parseval's identity and $|\widehat{\mathbf1_P}(\chi)|\le|P|$, we obtain $\|T\mathbf1_P\|_2^2\le6.4\cdot10^{11}|P|+t^2e^{2\tau t}|P|^2e^{-r/10}$.
The second term is at most $|P|$, since
$2\log t+2\tau t+\log|P|-r/10
\le(1+1/4000)r/10^{12}-9r/100<0$.
Thus $\|T\mathbf1_P\|_2^2<10^{12}|P|$.
\end{proof}

\section{Rearrangements via dissociated sets}\label{sec:blocks}
\subsection{Rainbow path decomposition}

This section decomposes \(R\) into a small number of path orderings
and controls the corresponding vertex sets in both directions.

\begin{lemma}\label{lem:path-decomposition}
Let \(0<\eta<1\), \(0\le\tau\le1/8000\), and
\(S=R\mathbin{\dot\cup}D\subseteq G\setminus\{0\}\), where \(D\) is
dissociated, \(|D|=r\), and \(R\ne\varnothing\). Suppose that
\(t\ge3\), \(r\ge10^{14}\max\{t,\eta^{-1}\}\), and
\(|R|\le e^{r/100}\). Let \(T=\sigma(S)\). Then \(R\) can be
partitioned into path orderings \(B_1,\ldots,B_q\) such that
\(q\le10^{14}\eta^{-2}\log(e|R|)\) and
\begin{enumerate}[label=\textup{(\roman*)}]
\item \(|B_1|=1\) and \(T\notin V(B_1)\);
\item \(\Lambda(0,V(B_i))\le\eta\) and
\(\Lambda(0,V(\overleftarrow B_i))\le\eta\) for every \(i\in[q]\).
\end{enumerate}
\end{lemma}

\begin{proof}
We first prove the following claim.
\begin{claim}\label{c1}
Let \(3\le\ell\le e^{r/100}\), and let
\(X\subseteq G\setminus\{0\}\) satisfy
\(|X|\ge10^{13}\eta^{-2}\ell\). There is a path ordering \(B\) of
length \(\ell\), using colors from \(X\), such that
\(\Lambda(0,V(B))\le\eta\) and
\(\Lambda(0,V(\overleftarrow B))\le\eta\).
\end{claim}

Choose a path ordering \(C=(c_1,\ldots,c_{\ell-2})\) from \(X\) greedily. After \(k-1\) steps with partial sums \(p_0=0,p_1,\ldots,p_{k-1}\), avoid used colors and \(\{p_i-p_{k-1}:0\le i<k\}\), excluding fewer than \(2k\) colors.

Let \(P=\{0\}\cup V(C)\), \(Q=\{0\}\cup V(\overleftarrow C)\), and
\(X'=X\setminus\{c_1,\ldots,c_{\ell-2}\}\). We have
\(|P|=|Q|=\ell-1\le e^{r/100}\). Since
\(\Lambda(0,x+P)=\Lambda(-x,P)\), we obtain from Lemma~\ref{lem:cube-mixing} that
\[
 |\{x\in G:\Lambda(0,x+P)>\eta/2\}|
 \le4\cdot10^{12}\eta^{-2}|P|.
\]
At most \(|P|\) further choices have \(0\in x+P\). Hence there is
\(x\in X'\) such that \(\Lambda(0,x+P)\le\eta/2\) and
\(0\notin x+P\).
Similarly, choose \(y\in X'\setminus\{x\}\) such that
\(\Lambda(0,y+Q)\le\eta/2\), \(0\notin y+Q\), and \(x+\sigma(C)+y\ne0\). These requirements, including \(y\ne x\), exclude at most \((4\cdot10^{12}\eta^{-2}+1)|Q|+2<5\cdot10^{12}\eta^{-2}\ell\)
elements of \(X'\), whereas \(|X'|>9\cdot10^{12}\eta^{-2}\ell\).
Let $B=(x,c_1,\ldots,c_{\ell-2},y)$. By the choice of $x$ and $y$, $B$ is a path ordering. Moreover, \(V(B)=(x+P)\cup\{x+\sigma(C)+y\}\) and
\(V(\overleftarrow B)=(y+Q)\cup\{x+\sigma(C)+y\}\).
By Fact \ref{fact:subadditivity} and Lemma~\ref{lem:one-vertex}, both
\(\Lambda(0,V(B))\) and \(\Lambda(0,V(\overleftarrow B))\) are at most \(\eta/2+2/r\le\eta\), proving the claim. \hfill\(\blacksquare\)
\medskip

Choose \(x_0\in R\setminus\{T\}\) and let \(B_1=(x_0)\).
Such an \(x_0\) exists, since otherwise \(R=\{T\}\) and
\(\sigma(D)=T-\sigma(R)=0\), a contradiction to dissociativity.
Let \(m=\lceil10^{13}\eta^{-2}\rceil\) and \(Y=R\setminus\{x_0\}\).
While \(|Y|\ge3m\), set \(\ell=\lfloor|Y|/m\rfloor\) and apply
Claim~\ref{c1} to \(X=Y\), then remove the chosen colors from \(Y\).
Indeed, \(3\le\ell\le|R|\le e^{r/100}\) and
\(|Y|\ge m\ell\ge10^{13}\eta^{-2}\ell\).
When \(|Y|<3m\), let every remaining color form a one-term path
ordering. Lemma~\ref{lem:one-vertex} bounds both potentials of each
one-term ordering by \(2/r\le\eta\), so \textup{(i)} and
\textup{(ii)} hold.

Let $k$ be the number of blocks chosen while $|Y|\ge3m$, and let $Y_i$ be the set remaining after the first $i$ such blocks, with $Y_0=R\setminus\{x_0\}$. For $0\le i<k$, each chosen block has size at least $|Y_i|/(2m)$, so
$
|Y_i|\le |R|e^{-i/(2m)}.
$
If $k\ge1$, then $|Y_{k-1}|\ge3m$, which implies $k\le1+2m\log|R|\le2m\log(e|R|)$. This bound also holds when $k=0$. Since fewer than $3m$ colors remain, the initial block and the remaining singletons contribute at most $3m$ blocks. Consequently,
$
q\le k+3m\le5m\log(e|R|)
\le10^{14}\eta^{-2}\log(e|R|).
$
\end{proof}

\subsection{Path orderings along the cube path}\label{sec:splicing}

We insert the path orderings from Lemma~\ref{lem:path-decomposition} at regularly spaced positions in a random ordering of $D$, and exclude proper zero-sum intervals by their number of colors from $D$.

\begin{theorem}\label{thm:splicing}
There are absolute constants \(c,C>0\) such that the following holds.
Let
\(S=R\mathbin{\dot\cup}D\subseteq G\setminus\{0\},\)
where \(D\) is dissociated and \(|D|=r\), and let \(n=|S|\). Suppose
that \(t\in\mathbb N\) and \(t\ge C\),
\(r\ge Ct\log(en),\)
and
\begin{equation}\label{R}
|R|\le
c(en)^{-C}
\sqrt{\frac1r\binom r{t+1}}.
\end{equation}
\begin{itemize}
\item[(i)] If \(\sigma(S)\ne0\), then \(\Cay(G,S)\) contains a rainbow
directed path of length \(n\) using every color.
\item[(ii)] If \(\sigma(S)=0\), then \(\Cay(G,S)\) contains a rainbow
directed cycle of length \(n\) using every color.
\end{itemize}
\end{theorem}
\begin{proof}
Choose $0<c\le1/10$ and $C\ge10^{25}$, and let $\tau=1/16000$. If $R=\varnothing$, then $\sigma(S)=\sigma(D)\ne0$ by dissociativity, and any ordering of $D$ determines the required path by Lemma~\ref{lem:boolean-cube}. Hence assume that $R\ne\varnothing$. Since $\log(e|R|)\le\log(en)\le r/(Ct)\le r/C^2$, we have $|R|\le e^{r/100}$. The hypotheses also ensure $t\ge3$, $r\ge10^{14}\max\{t,100\}$, and $t\le r/20$.

Let $T=\sigma(S)$. By Lemma~\ref{lem:path-decomposition} with $\eta=1/100$, we obtain path orderings $B_1,\ldots,B_q$ that partition $R$, where $q\le10^{18}\log(en)$, $B_1=(x)$ with $x\ne T$, and $\Lambda(0,V(B_i)),\Lambda(0,V(\overleftarrow B_i))\le1/100$ for every $i\in[q]$.

For $i\in[q]$, let $U_i=\{(i-1)(t+1)+1,\ldots,i(t+1)\}$ and $\mathcal U=\{U_1,\ldots,U_q\}$. Let $\pi=(\pi_1,\ldots,\pi_r)$ be a uniformly random ordering of $D$. Insert $B_i$ immediately before $\pi_{(i-1)(t+1)+1}$, and denote the resulting rainbow walk by $W(\pi)$; see Figure \ref{fig:main-splicing}.

Let $\mathcal G$ be the event that, for every $i\in[q]$ and $1\le j<t$,
\[
 \pi_{(i-1)(t+1)+1}+\cdots+\pi_{(i-1)(t+1)+j}
 \notin -V(\overleftarrow B_i),
\]
and, for $i<q$,
\[
 \pi_{i(t+1)-j+1}+\cdots+\pi_{i(t+1)}\notin -V(B_{i+1}).
\]
For $i\in[q]$ and $1\le j<t$, let $L_{i,j}$ and $R_{i,j}$ denote the first and last $j$ positions of
$U_i$, respectively. Define
\[
\begin{aligned}
\mathcal C(U_i)
={}\{(L_{i,j},-V(\overleftarrow B_i)):1\le j<t\}\cup\{(R_{i,j},-V(B_{i+1})):1\le j<t,\ i<q\}.
\end{aligned}
\]
For each $i\in[q]$, the left-hand side of \eqref{e1} equals $\Lambda(0,V(\overleftarrow B_i))+\mathbf1_{\{i<q\}}\Lambda(0,V(B_{i+1}))$, which is at most $1/50$. Hence, by Lemma~\ref{lem:conditioning} with $M=q(t+1)$, we have $\Prob(\mathcal G)\ge e^{-2q/25}>0$.

A nonempty consecutive subwalk is \emph{closed} if and only if its colors sum to zero. We show that, conditional on $\mathcal G$, with positive probability $W(\pi)$ contains no proper nonempty closed subwalk. For a proper nonempty consecutive subwalk of $W(\pi)$, let $h$ be the number of its colors from $D$. If $h=0$, the subwalk lies in one $B_i$ and is not closed. If $h=r$, the subwalk contains both $\pi_1$ and $\pi_r$, and hence every color except possibly the initial color $x$. Thus a proper subwalk must exclude $x$ and has sum $T-x\ne0$. Now we consider the remaining three cases.

\smallskip
\begin{case}
$1\le h<t$.
\end{case}
In this case, the subwalk meets at most one $B_i$. If it meets no block, it is not closed by dissociativity. If it meets a block and all its $D$-colors lie on one side of that block, it is excluded by $\mathcal G$. Thus it contains an entire $B_i$ and at least one $D$-color on each side. Thus $h\ge2$, and for each $h$ there are at most $q(h-1)$ such subwalks.

Fix such a subwalk containing $B_i$, and let $I\subseteq[r]$ be its set of $D$-positions. The subwalk is closed if and only if $\sum_{j\in I}\pi_j=-\sigma(B_i)$. Since $|I|=h<t\le r/20$ and $D$ is dissociated, we have $\mu_h^D(\{-\sigma(B_i)\})\le\binom rh^{-1}$.
Thus, by Lemma~\ref{lem:conditioning}\textup{(ii)}, we obtain
\[
\begin{aligned}
\Prob\left(\sum_{j\in I}\pi_j=-\sigma(B_i)\,\middle|\,\mathcal G\right)
&\le \exp\left(\frac{4q(t+1)h}{r}+\frac4{25}\right)\binom rh^{-1}\le e^{\tau h+4/25}\binom rh^{-1},
\end{aligned}
\]
where the last inequality uses $4q(t+1)/r\le\tau$. By the union bound,
\[
\begin{aligned}
&\Prob(\exists\text{ a closed subwalk with }2\le h<t\mid\mathcal G)\le \sum_{h=2}^{t-1}q(h-1)e^{\tau h+4/25}\binom rh^{-1}\\
&\le q e^{4/25+2\tau}\binom r2^{-1}
\sum_{h=2}^{\infty}\left(\frac{2e^\tau t}{r-t}\right)^{h-2}\le 2q e^{4/25+2\tau}\binom r2^{-1}
\le\frac{10q}{r^2}.
\end{aligned}
\]

\smallskip
\begin{case}
$t\le h\le r-t$.
\end{case}
Fix consecutive positions $a,\ldots,b$ in $\pi$, where $b-a+1=h$. A subwalk whose $D$-colors are exactly $\pi_a,\ldots,\pi_b$ contains the entire segment from $\pi_a$ to $\pi_b$. It may additionally contain a suffix of an $R$-block immediately before $\pi_a$ and a prefix of an $R$-block immediately after $\pi_b$, either of which may be empty. Thus there are at most $(|R|+1)^2$ such subwalks. Each choice fixes all colors from $R$ in the subwalk, so a closed walk requires the $h$ colors from $D$ to have a fixed sum. There are $r-h+1$ intervals of length $h$ in $\pi$. Since $\binom rh\ge\binom rt$ for $t\le h\le r-t$, $|R|+1\le2|R|$, and $e^{2q/25}\le(en)^C$ by the bound on $q$, Lemma~\ref{lem:conditioning}\textup{(iii)} and the union bound show that
\[
\begin{aligned}
&\Prob\bigl(\exists\text{ a closed subwalk with }t\le h\le r-t\mid\mathcal G\bigr)\\
&\quad\le 
\sum_{h=t}^{r-t}(r-h+1)(|R|+1)^2e^{2q/25}\binom rh^{-1}\\
&\quad\le 4r^2|R|^2(en)^C\binom rt^{-1}\le 4c^2r(en)^{-C}\frac{\binom r{t+1}}{\binom rt}\\
&\quad=4c^2\frac{r(r-t)}{t+1}(en)^{-C}
\le4c^2(en)^{2-C}
\le4c^2.
\end{aligned}
\]
Here the bound involving $c^2$ follows from \eqref{R}, the equality from $\binom r{t+1}/\binom rt=(r-t)/(t+1)$, and the final two bounds from $r\le n$ and $C\ge2$.

\smallskip
\begin{case}
$h=r-s$, where $1\le s<t$.
\end{case}
For $a\in[s+1]$, let $I_a=\{1,\ldots,a-1\}\cup\{a+r-s,\ldots,r\}$ be the missing $D$-positions. Since $a\le t$ and $a+r-s-1\ge r-s>q(t+1)$, the subwalk contains $B_2,\ldots,B_q$. If $a>1$, the subwalk omits $x$. If $a=1$, it starts either with $x$, including $B_1$, or with $\pi_1$, omitting $B_1$. Thus a closed subwalk must satisfy one of the $s+2$ equalities
$
\sum_{j\in I_1}\pi_j=T$,
or
$\sum_{j\in I_a}\pi_j=T-x$ for some
$a\in[s+1]$.
We have $|I_a|=s\le r/20$ and $I_a\cap(\bigcup\mathcal U)\subseteq U_1$, so $m_{\mathcal U}(I_a)\le1$. By Lemma~\ref{lem:conditioning}\textup{(ii)} and dissociativity, each of these events has conditional probability at most $\exp(4q(t+1)s/r+2/25)\binom rs^{-1}\le e^{\tau s+2/25}\binom rs^{-1}$. Again by the union bound,
\[
\begin{aligned}
&\Prob\bigl(\exists\text{ a closed subwalk with }r-t<h<r\mid\mathcal G\bigr)\\
&\quad\le e^{2/25}\sum_{s=1}^{t-1}(s+2)e^{\tau s}\binom rs^{-1}\\
&\quad\le \frac{3e^{2/25+\tau}}{r}
\sum_{s=1}^{\infty}\left(\frac{4e^\tau t}{3(r-t)}\right)^{s-1}\\
&\quad\le \frac{6e^{2/25+\tau}}{r}
\le\frac7r.
\end{aligned}
\]
\medskip

Combining the three cases, we obtain that
\[
\Prob\bigl(\exists\text{ a proper nonempty closed subwalk of }W(\pi)
\mid\mathcal G\bigr)
\le \frac{10q}{r^2}+4c^2+\frac7r<1.
\]
Thus some $W(\pi)$ has no proper nonempty closed subwalk. If $T\ne0$, the whole walk is not closed, so all its vertices are distinct and it is a rainbow directed path of length $n$. If $T=0$, only its initial and final vertices coincide, so it is a rainbow directed cycle of length $n$.
\end{proof}

\subsection{Proof of Theorem~\ref{thm:main}}\label{sec:dimension-proof}

\begin{proof}[Proof of Theorem~\ref{thm:main}]
Let $c,C$ be the constants in Theorem~\ref{thm:splicing}. We first show that there are constants $A_0>0$ and $n_0$ such that the conclusion holds for every $n\ge n_0$ under the hypothesis $\diss(S)\ge A_0(\log n)^2/\log\log n$. Let $m=\log n$, choose a constant $\alpha>2(C+2)$, and then choose $A_0$ sufficiently large in terms of $C$ and $\alpha$. Let $r=\lfloor A_0m^2/\log m\rfloor$ and $t=\lfloor\alpha m/\log m\rfloor$. Under this hypothesis, $S$ contains a dissociated set $D$ of size $r$; let $R=S\setminus D$.

For large $n$, one has $t\ge C$ and $r\ge Ct\log(en)$, since $r/(t\log(en))=A_0/\alpha+o(1)$. It remains to verify the bound on $|R|$. Since $\binom r{t+1}\ge(r/(t+1))^{t+1}$ and $r/(t+1)=(A_0/\alpha+o(1))m$,
$
   \frac12\log\binom r{t+1}
   \ge \frac{t+1}{2}\log\frac r{t+1}
   =\left(\frac\alpha2+o(1)\right)m.
$
Also $\log r=o(m)$. Therefore the logarithm of the right-hand side of \eqref{R} is at least $(-C+\alpha/2+o(1))m$. Our choice of $\alpha$ makes this at least $(1+\delta)m$ for some fixed $\delta>0$ and all sufficiently large $n$. The right-hand side therefore exceeds $n$, and hence exceeds $|R|$. By Theorem~\ref{thm:splicing}, we obtain the required rainbow path or cycle.

Finally, choose $A\ge A_0$ sufficiently large that
$A(\log n)^2/\log\log n>n$ for $3\le n<n_0$.
Since $\diss(S)\le n$, this makes the hypothesis impossible for these remaining values of $n$.
\end{proof}

\section{Coset concentration and slice anticoncentration}\label{sec:concentrated-core}
We first show that when $S$ has small additive dimension, every coset trace is determined by a small subset of its elements. This allows us to bound the number of coset traces and to find a coset containing almost all of $S$ whose intersection with $S$ is not concentrated in any proper subcoset.

\subsection{A concentrated coset core}
\begin{lemma}\label{lem:concentrated-core}
Fix $C_0>0$ and $0<\alpha<1/4$. Let $S\subseteq G\setminus\{0\}$ have size $n$, where $n$ is sufficiently large in terms of $C_0$ and $\alpha$, and suppose that $\diss(S)\le C_0(\log n)^2$. Then there exist $a\in G$ and $H\le G$ such that, for $A=S\cap(a+H)$,
\begin{equation}\label{eq:core-properties}
|S\setminus A|\le n^{2\alpha}\quad\text{and}\quad |A\setminus(x+K)|\ge n^\alpha
\end{equation}
for every proper subgroup $K<H$ and every $x\in G$.
\end{lemma}

Before proving Lemma \ref{lem:concentrated-core}, we first give a proof of Theorem \ref{lem:coset-traces}.

\begin{proof}[Proof of Theorem \ref{lem:coset-traces}]
Choose a maximum dissociated set $D_0\subseteq S$. Since $S\subseteq\Span(D_0)$, every element of $S-a$ has a representation over $D_0$ with coefficients in $\{-2,-1,0,1,2\}$. If $E\subseteq S-a$ is dissociated and $|E|=h$, its $2^h$ distinct subset sums imply $2^h\le(4h+1)^d$. Thus $v=h/d\le\log_2(4dv+1)$. If $v\ge4\log_2(4d+1)$, then $v\ge4$ and $\log_2(4dv+1)\le\log_2(4d+1)+\log_2v\le v/4+v/2<v$, a contradiction. This proves~\textup{(i)}.

Now we prove (ii). For a nonempty trace $T\in\mathcal T(S)$, choose $x\in G$ and $K\le G$ such that $T=S\cap(x+K)$. Choose $b\in T$ and an inclusion-maximal dissociated set $E\subseteq T-b$. By~\textup{(i)}, $|E|\le M$. Since $E\subseteq T-b\subseteq\Span(E)$, $\langle T-b\rangle=\langle E\rangle\le K$. Consequently, $T\subseteq S\cap(b+\langle E\rangle)\subseteq S\cap(b+K)=T$, and hence $T=S\cap(b+\langle E\rangle)$.
With $F=\{b\}\cup(b+E)$, we have $|F|\le M+1$, $\aff(F)=b+\langle E\rangle=\aff(T)$, and $T=\mathcal C_S(F)$. For $T=\varnothing$, take $F=\varnothing$. Since every trace has the form $\mathcal C_S(F)$ for some $F\subseteq S$ with $|F|\le M+1$, there are at most $\sum_{j=0}^{\min\{M+1,n\}}\binom nj$ traces, giving \eqref{eq:algebraic-container-count}.

To bound the VC-dimension, let $\varnothing\ne Y\subseteq S$ be shattered by $\mathcal T(S)$. Choose $b\in Y$ and an inclusion-maximal dissociated set $E\subseteq Y-b$, and let $F=\{b\}\cup(b+E)$. By~\textup{(i)}, $|E|\le M$, so $|F|\le M+1$. Since $Y-b\subseteq\Span(E)$, $Y\subseteq b+\langle E\rangle=\aff(F)$. Since $Y$ is shattered, there is a coset $C$ such that $Y\cap C=F$. Since $F\subseteq C$, we have $\aff(F)\subseteq C$, so $Y\subseteq C$. Thus $Y=Y\cap C=F$ and $|Y|\le M+1$. This proves (ii).

Now we prove (iii). Fix $0<\kappa<n/(M+1)$. Start with $A_0=S$ and $\mathcal C_0=\aff(S)$. Whenever a proper subcoset $\mathcal D_i\subsetneq\mathcal C_i$ satisfies $|A_i\setminus\mathcal D_i|<\kappa$, let $A_{i+1}=A_i\cap\mathcal D_i$ and $\mathcal C_{i+1}=\aff(A_{i+1})$. Throughout the first $M+1$ steps the sets $A_i$ are nonempty, since $|A_i|\ge n-i\kappa>0$. Also $A_i=S\cap\mathcal C_i$, and $\mathcal C_{i+1}\subsetneq\mathcal C_i$.

If the procedure continued at least $M+1$ steps, choose $b\in A_{M+1}$, write $\mathcal C_i=b+H_i$, and choose $s_i\in A_i\setminus A_{i+1}$ for $0\le i\le M$. Then $H_0>\cdots>H_{M+1}$ and $s_i-b\in H_i\setminus H_{i+1}$. We claim that $s_0-b,\ldots,s_M-b$ are dissociated: in any nontrivial $\{-1,0,1\}$-relation, the term of smallest index lies outside $H_{i+1}$, whereas all later terms belong to $H_{i+1}$. This contradicts~\textup{(i)}.

Thus the procedure stops after at most $M$ steps. For its final trace $A=S\cap\mathcal C_h$, we have $|S\setminus A|\le M\kappa$ and $|A\setminus\mathcal D|\ge\kappa$ for every proper subcoset $\mathcal D\subsetneq\mathcal C_h$. By~\textup{(ii)}, choose $F\subseteq A$ of size at most $M+1$ with $A=\mathcal C_S(F)$ and $\aff(F)=\aff(A)=\mathcal C_h$. If $K<\langle F-F\rangle$, then every coset $x+K$ is either a proper subcoset of $\mathcal C_h$ or is disjoint from it. In the latter case, $|A|\ge n-M\kappa>\kappa$. This proves~\textup{(iii)}.
\end{proof}
Now we are ready to prove Lemma~\ref{lem:concentrated-core}.
\begin{proof}[Proof of Lemma~\ref{lem:concentrated-core}]
Let $M$ be as in Theorem~\ref{lem:coset-traces}, and set $\kappa=n^\alpha$. Since $M=O_{C_0}((\log n)^2\log\log n)$, we have $(M+1)\kappa<n$ and $M\kappa\le n^{2\alpha}$ for sufficiently large $n$.
By Theorem~\ref{lem:coset-traces}\textup{(iii)}, choose a nonempty $F\subseteq S$ such that $A=S\cap\aff(F)$ satisfies \eqref{eq:algebraic-stable-container}. Writing $\aff(F)=a+H$ gives \eqref{eq:core-properties}.
\end{proof}

\subsection{Slice anticoncentration}\label{sec:slice-anticoncentration}
Throughout this subsection, $H$ is a finite abelian group. The connection between concentration of random sums and additive structure is a central theme of inverse Littlewood--Offord theory~\cite{TaoVuInverse,NguyenVu}. Campos, Mattos, Morris and Morrison~\cite{CamposMattosMorrisMorrison} developed a container approach to this theory over prime fields. Here we consider sums over uniformly chosen subsets of prescribed size. The next lemma bounds the probability that a uniform subset of prescribed size has a given sum.

\begin{lemma}\label{lem:one-slice-traces}
There is an absolute constant $C\ge1$ such that the following holds. Let $G$ be a finite abelian group, let $H\le G$, and let $U\subseteq G$ have size $u$. Let $0<\rho\le1/4$ and $T\ge1$. Suppose that there is $a\in G$ such that, writing $B=U\cap(a+H)$, the following conditions hold:
\begin{enumerate}[label=\textup{(\roman*)}]
\item $|B|\ge3u/4$;
\item $|(B-a)\setminus(x+K)|\ge\rho|B|$ for every proper subgroup $K<H$ and every $x\in H$;
\item the cosets of proper subgroups of $H$ induce at most $T$ distinct traces on $B-a$, and $\rho|B|\ge C\log(2|B|T)$.
\end{enumerate}
If $R$ is a uniformly random $m$-subset of $U$, where $1\le m\le3u/4$, then
\begin{equation}\label{eq:one-slice-bound}
\sup_{z\in G}\Prob(\sigma(R)=z)
\le\frac1{|H|}+\frac{C}{u\sqrt{\rho m}}.
\end{equation}
\end{lemma}

Informally speaking, conditions~\textup{(i)} and~\textup{(ii)} say that most of $U$ lies in one coset of $H$, and this large part is not concentrated in any smaller coset. This prevents a natural obstruction to anticoncentration: If $B$ lay entirely in a coset of a proper subgroup $K<H$, then sums of a fixed number of elements of $B$ would also be confined to a single coset of $K$. In the proof, we need the random remaining sets to retain sufficiently many elements outside every coset of a proper subgroup of $H$. Two cosets with the same intersection with $B-a$ give the same event, so for each random set we only need to consider at most $T$ distinct intersections. The lower bound in~\textup{(iii)} ensures that the total failure probability is small after taking a union bound.

For a probability distribution $f$ on $H$, write $\mathsf T_gf(x)=f(x-g)$ and $\|f\|_2^2=\sum_{x\in H}f(x)^2$. The quantity $\|f\|_2^2$ is the probability that two independent vertices sampled from $f$ coincide. For $Q\subseteq H$ and an integer $r\ge1$, let $rQ$ be the set of sums of $r$ elements of $Q$, allowing repetitions. 

\begin{lemma}\label{lem:slice-neighbourhood-packing}
Let $Q=-Q\subseteq H$ contain $0$ and generate $H$. If $r\ge1$ is an integer and $rQ\ne H$, then $|rQ|\ge r|Q|/3$.
\end{lemma}
\begin{proof}
Let $\Gamma$ be the graph on $H$ in which two distinct vertices $x,y$ are adjacent if $y-x\in Q$. Since $Q=-Q$ and $0\in Q$, the set $rQ$ consists precisely of the vertices at distance at most $r$ from $0$ in $\Gamma$. As $Q$ generates $H$ and $rQ\ne H$, there is a shortest path $v_0v_1 \ldots v_{r+1}$ with $v_0=0$ and some $v_{r+1} \not\in rQ$. For $i=0,3,6,\ldots$ with $i<r$, the sets $v_i+Q$ are contained in $rQ$ and are pairwise disjoint, since any two of their centers are at distance at least three. There are $\lceil r/3\rceil$ such sets, each of size $|Q|$. Hence $|rQ|\ge\lceil r/3\rceil|Q|\ge r|Q|/3$.
\end{proof}

Let $X$ have distribution $f$. Adding $a$ or $b$ to $X$, each with probability $1/2$ and independently of $X$, gives the distribution $(\mathsf T_af+\mathsf T_bf)/2$. By the parallelogram identity, we have
\begin{equation}\label{eq:pair-averaging-identity}
\left\|\frac{\mathsf T_af+\mathsf T_bf}{2}\right\|_2^2
=\|f\|_2^2-\frac14\|\mathsf T_af-\mathsf T_bf\|_2^2.
\end{equation}
Thus averaging two translates cannot increase the squared $\ell_2$-norm. We next give a lower bound on the average decrease over all pairs $a,b\in P$.

\begin{lemma}\label{lem:pair-averaging}
Let $P\subseteq H$ have size $w$, let $h=|H|$, and let $0<\rho\le1/4$. Suppose that $|P\setminus(a+K)|\ge\rho w$ for every coset of every proper subgroup $K<H$. If $f$ is a probability distribution on $H$ with $\|f\|_2^2\le16/w$, then
\begin{equation}\label{eq:pair-averaging}
\frac1{w^2}\sum_{a,b\in P}\|\mathsf T_af-\mathsf T_bf\|_2^2
\ge c_0\rho w^2\bigl(\|f\|_2^2-h^{-1}\bigr)^3,
\end{equation}
where $c_0=2^{-20}$.
\end{lemma}

\begin{proof}
Write $v=\|f\|_2^2-h^{-1}$ and assume $v>0$. For $g\in H$, define $d(g)=\|f-\mathsf T_gf\|_2$, and let $r=\lceil16/(wv)\rceil\le32/(wv)$ and $Q=\{g\in H:d(g)^2\le v/r^2\}$. 
Since $\mathsf T_0f=f$ and translations preserve the $\ell_2$-norm, we have $d(0)=0$ and $d(-g)=\|\mathsf T_gf-f\|_2=d(g)$ for every $g\in H$. Thus $0\in Q$ and $Q=-Q$. Moreover, the triangle inequality and translation invariance give $d(g+h)\le d(g)+d(h)$ for all $g,h\in H$. For any $g\in rQ$, write $g=g_1+\cdots+g_r$ with $g_i\in Q$. By the definition of $Q$, we have $d(g_i)\le\sqrt v/r$, and hence $d(g)\le\sum_{i=1}^r d(g_i)\le\sqrt v$. Therefore $d(g)^2\le v$.


Since $h^{-1}\sum_{g\in H}d(g)^2=2v$, we have $rQ\ne H$. For $g\in rQ$, we also have $\sum_{x\in H}f(x)f(x-g)=\|f\|_2^2-d(g)^2/2\ge v/2$. As $\sum_{g,x\in H}f(x)f(x-g)=1$, it follows that $|rQ|\le2/v$.

If $\langle Q\rangle=H$, then by Lemma~\ref{lem:slice-neighbourhood-packing}, $|Q|\le3|rQ|/r\le6/(rv)\le3w/8$, so $|P\setminus(b+Q)|\ge5w/8\ge\rho w$. If $\langle Q\rangle<H$, the hypothesis ensures $|P\setminus(b+Q)|\ge\rho w$. Thus, for every $b\in P$, at least $\rho w$ elements $a\in P$ satisfy $d(a-b)^2>v/r^2$. Consequently,
\[
\frac1{w^2}\sum_{a,b\in P}\|\mathsf T_af-\mathsf T_bf\|_2^2
\ge\frac{\rho v}{r^2}
\ge2^{-10}\rho w^2v^3
\ge c_0\rho w^2v^3.\qedhere
\]
\end{proof}
\begin{proof}[Proof of Lemma~\ref{lem:one-slice-traces}]
We first prove the following claim.
\begin{claim}
Let $P\subseteq H$ have size $u$, and let $0<\kappa\le u/4$. Suppose that every coset $C'$ of a proper subgroup of $H$ satisfies $|P\setminus C'|\ge\kappa$, that these cosets induce at most $T$ distinct traces on $P$, and that $\kappa\ge512\log(2uT)$. If $R$ is a uniformly random $m$-subset of $P$, where $1\le m\le3u/4$, then
\begin{equation}\label{eq:pure-slice-estimate}
\sup_{z\in H}\Prob(\sigma(R)=z)\le\frac1{|H|}+\frac{C_0}{\sqrt{u\kappa m}},
\end{equation}
where $C_0$ is an absolute constant.
\end{claim}

If $m<16$, we sample $R$ by first sampling a uniformly random $(m-1)$-subset $R'$ of $P$, and then choosing a uniformly random element of $P\setminus R'$. Conditional on any outcome of $R'$, at most one choice gives $\sigma(R)=z$. Thus $\Prob(\sigma(R)=z)\le1/(u-m+1)\le2/u\le4/\sqrt{u\kappa m}$.

We may therefore assume that $m\ge16$. Let $s=\lfloor m/16\rfloor$ and $\ell=\lceil u/16\rceil$. We sample $R$ in the following way. First, divide an initial segment of a uniformly random ordering of $P$ into sets $F,A_0,A_1$ of sizes $m-2s-2,\ell,\ell$, followed by $2s$ ordered pairs $(a_i,b_i)$. This is possible since $m+2s+2\ell-2\le31u/32$. We then take all elements of $F$ and independently choose one uniformly random element from each of $A_0,A_1$ and the $2s$ pairs. The resulting set $R$ is a uniformly random $m$-subset of $P$, since it has size $m$ and its distribution is invariant under permutations of $P$. The sampling procedure is illustrated in Figure~\ref{fig:pair-sampling}.

\begin{figure}[t]
\centering
\resizebox{0.85\linewidth}{!}{%
\begin{tikzpicture}[x=1cm,y=1cm,font=\small,>=Latex,line width=.5pt]
\definecolor{samplingblue}{RGB}{42,85,125}
\tikzset{selected/.style={draw=samplingblue,fill=samplingblue!12,line width=.65pt},choice/.style={selected,circle,minimum size=8mm,inner sep=0pt},arrow/.style={->,draw=black!65,line width=.6pt}}
\node[anchor=west] at (0,3.35) {A uniformly random ordering of $P$};
\draw[draw=black!70,fill=black!3] (0,1.55) rectangle (2.9,2.8);
\node at (1.45,2.35) {$F$};
\node[font=\footnotesize] at (1.45,1.92) {$|F|=m-2s-2$};
\draw[draw=black!70,fill=black!3] (3.2,1.55) rectangle (5,2.8);
\node at (4.1,2.35) {$A_0$};
\node[font=\footnotesize] at (4.1,1.92) {$|A_0|=\ell$};
\draw[draw=black!70,fill=black!3] (5.3,1.55) rectangle (7.1,2.8);
\node at (6.2,2.35) {$A_1$};
\node[font=\footnotesize] at (6.2,1.92) {$|A_1|=\ell$};
\draw[draw=black!70,fill=black!3] (7.4,1.55) rectangle (11.7,2.8);
\draw[draw=black!35] (8.85,1.55)--(8.85,2.8) (10.25,1.55)--(10.25,2.8);
\node at (8.125,2.42) {$a_1$};
\node at (8.125,1.91) {$b_1$};
\node at (9.55,2.175) {$\cdots$};
\node at (10.975,2.42) {$a_{2s}$};
\node at (10.975,1.91) {$b_{2s}$};
\node[font=\footnotesize] at (9.55,3.02) {$2s$ pairs};
\draw[draw=black!40,dashed,fill=black!2] (12,1.55) rectangle (13.5,2.8);
\node[align=center,font=\footnotesize,text=black!60] at (12.75,2.175) {unused\\elements};
\draw[arrow] (1.45,1.55)--(1.45,.24) node[midway,fill=white,inner sep=3pt,font=\footnotesize] {take all};
\draw[arrow] (4.1,1.55)--(4.1,.24) node[midway,fill=white,inner sep=3pt,font=\footnotesize] {choose one};
\draw[arrow] (6.2,1.55)--(6.2,.24) node[midway,fill=white,inner sep=3pt,font=\footnotesize] {choose one};
\draw[arrow] (8.125,1.55)--(8.125,.24);
\draw[arrow] (10.975,1.55)--(10.975,.24);
\node[font=\footnotesize,align=center] at (9.55,.9) {choose one\\from each pair};
\node[anchor=east] at (-.15,-.2) {$R:$};
\draw[selected] (0,-.62) rectangle (2.9,.22);
\node at (1.45,-.2) {$F$};
\node[choice] at (4.1,-.2) {$\alpha_0$};
\node[choice] at (6.2,-.2) {$\alpha_1$};
\node[choice] at (8.125,-.2) {$c_1$};
\node at (9.55,-.2) {$\cdots$};
\node[choice] at (10.975,-.2) {$c_{2s}$};
\draw[decorate,decoration={brace,mirror,amplitude=5pt},draw=black!65] (0,-.9)--(11.7,-.9) node[midway,below=10pt] {$|R|=(m-2s-2)+1+1+2s=m$};
\end{tikzpicture}%
}
\caption{Sampling $R$. Here $\alpha_j$ is the element chosen from $A_j$
for $j\in\{0,1\}$, and $c_i$ is the element chosen from $\{a_i,b_i\}$
for $i\in[2s]$. Thus $R=F\cup\{\alpha_0,\alpha_1,c_1,\ldots,c_{2s}\}$.}
\label{fig:pair-sampling}
\end{figure}
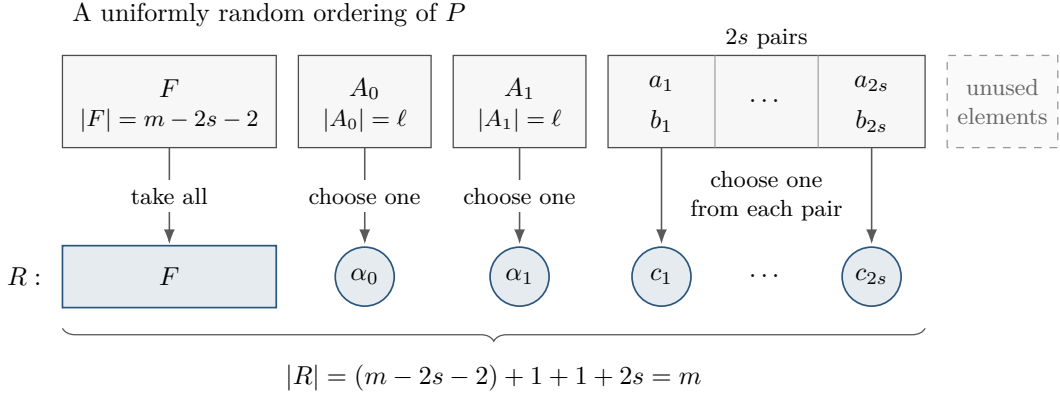

We expose $F,A_0,A_1$, and then the first $s$ pairs one at a time. For $0\le j\le s$, let $\mathcal F_j$ denote the $\sigma$-algebra generated by $F,A_0,A_1$ and the first $j$ pairs. Write $P_j=P\setminus(F\cup A_0\cup A_1\cup\bigcup_{i=1}^j\{a_i,b_i\})$ and $w_j=|P_j|$. Then $w_j\ge u-m+2-2\ell\ge u/8$. We call $P_j$ \emph{good} if $|P_j\setminus C'|\ge\kappa w_j/(2u)$ for every coset $C'$ of a proper subgroup of $H$. By Lemma~\ref{lem:hypergeometric} and a union bound over $0\le j\le s$ and the at most $T$ traces, we obtain
\begin{equation}\label{eq:pair-pool-failure}
\Prob(\text{some }P_j\text{ is not good})\le uT e^{-\kappa/64}\le u^{-7}.
\end{equation}

Let $f_j$ be the conditional distribution, given $\mathcal F_j$, of the sum of one uniformly chosen element from $A_0$ and one from each of the first $j$ pairs, with all choices independent. Write $v_j=\|f_j\|_2^2-|H|^{-1}$. By \eqref{eq:pair-averaging-identity}, $\|f_j\|_2^2\le\|f_0\|_2^2=1/\ell\le16/u$. Now fix $0\le j<s$ and condition on $\mathcal F_j$. The next pair is a uniformly random ordered pair of distinct elements of $P_j$. Whenever $P_j$ is good, by Lemma~\ref{lem:pair-averaging} with $P=P_j$ and $\rho=\kappa/(2u)$, together with \eqref{eq:pair-averaging-identity}, we obtain
\begin{equation}\label{eq:pair-step}
\mathbb E(v_{j+1}\mid\mathcal F_j)
\le v_j-\frac{c_0\kappa w_j^2}{8u}v_j^3
\le v_j-\lambda v_j^3,
\end{equation}
where $\lambda=c_0\kappa u/512$.

Let $W_j=v_j\mathbf1_{\{P_0,\ldots,P_j\text{ are good}\}}$ and $y_j=\mathbb EW_j$. By \eqref{eq:pair-step}, $\mathbb E(W_{j+1}\mid\mathcal F_j)\le W_j-\lambda W_j^3$. Taking expectations and using Jensen's inequality, we obtain $y_{j+1}\le y_j-\lambda y_j^3$. If $y_{j+1}>0$, it follows that $y_{j+1}^{-2}\ge y_j^{-2}+2\lambda$. Thus $y_s\le(2\lambda s)^{-1/2}$; if some $y_j=0$, this bound holds trivially. Recalling that $v_s\le16/u$ and $s\ge m/32$, we deduce from \eqref{eq:pair-pool-failure} that
\begin{equation}\label{eq:pair-endpoint-bound}
\mathbb E\|f_s\|_2^2
\le\frac1{|H|}+(2\lambda s)^{-1/2}+16u^{-8}
\le\frac1{|H|}+\frac{C_0}{\sqrt{u\kappa m}}.
\end{equation}
Let $g_s$ be the corresponding distribution obtained from $A_1$ and the last $s$ pairs. Interchanging $A_0,A_1$ and the first and last $s$ pairs does not change the distribution of the random ordering, so the same bound holds for $\mathbb E\|g_s\|_2^2$. Conditional on any outcome of the complete ordering, the two sums are independent and $\sigma(F)$ is fixed. Applying Cauchy--Schwarz first to their distributions and then to the expectation, we obtain $\Prob(\sigma(R)=z)\le(\mathbb E\|f_s\|_2^2\,\mathbb E\|g_s\|_2^2)^{1/2}$. This proves the claim. \hfill $\blacksquare$
\medskip

We now continue the proof of the lemma, first assuming that $m\le u/2$. Write $b=|B|$, $E=U\setminus B$, and $J=|R\cap B|$. Conditional on any outcome of $R\setminus B$, the value of $J$ is determined and $(R\cap B)-a$ is a uniformly random $J$-subset of $B-a$. Whenever $J\ge1$, we can apply \eqref{eq:pure-slice-estimate} with $P=B-a$ and $\kappa=\rho b$, since $J\le m\le2b/3$. Averaging over the outcomes of $R\setminus B$, we obtain
\[
\Prob(\sigma(R)=z,\ J\ge1)
\le\frac1{|H|}+\frac{C_0}{b\sqrt\rho}
\sum_{j\ge1}\frac{\Prob(J=j)}{\sqrt j}.
\]
Since $J$ is hypergeometric with mean at least $3m/4$, by Lemma~\ref{lem:hypergeometric} we have $\Prob(J<m/2)\le e^{-m/24}$. It follows that $\sum_{j\ge1}\Prob(J=j)/\sqrt j\le\sqrt{2/m}+e^{-m/24}\le4/\sqrt m$.

It remains to bound the probability that $\sigma(R)=z$ and $J=0$. Every $(m-1)$-subset of $E$ is contained in at most one $m$-subset of $E$ with sum $z$, and every such $m$-subset contains exactly $m$ subsets of size $m-1$. Hence, recalling that $|E|\le u/4$, we obtain
\[
\Prob(\sigma(R)=z,\ J=0)
\le\frac{\binom{|E|}{m-1}}{m\binom um}
\le\frac{2}{u}\,4^{-(m-1)}
\le\frac{2}{u\sqrt m}.
\]
Combining these estimates with $b\ge3u/4$ proves \eqref{eq:one-slice-bound} for $m\le u/2$, provided $C$ is sufficiently large in terms of $C_0$.

Finally, suppose that $u/2<m\le3u/4$. The set $U\setminus R$ is a uniformly random $(u-m)$-subset of $U$, and $\sigma(R)=z$ if and only if $\sigma(U\setminus R)=\sigma(U)-z$. Since $m/(u-m)\le3$, we apply the bound just proved to $U\setminus R$ and increase $C$ by a factor of $\sqrt3$, this implies \eqref{eq:one-slice-bound} and thus finishes the proof.
\end{proof}

\section{A large direction subgroup}\label{sec:large-subgroup}
In this section, we use the local-repair argument of Pham and Sauermann~\cite[Section~5]{PhamSauermann}, also see Costa, Della Fiore, Feng and Liu~\cite{CostaEtAlKneser}. Starting from a random ordering, the argument removes zero-sum intervals by disjoint exchanges of nearby entries, processing their right endpoints from right to left. We combine this criterion with the slice anticoncentration estimate from Section~\ref{sec:concentrated-core} to handle the case in which the direction subgroup of the concentrated coset is large. Related switching arguments have been used in graph packing~\cite{SauerSpencerPacking} and in the generation of random regular graphs~\cite{McKayWormaldRegular} and Latin rectangles~\cite{McKayWormaldLatin}. For random permutations, Harris and Srinivasan~\cite{HarrisSrinivasanSODA,HarrisSrinivasanPermutations} developed a general swapping algorithm.

The following lemma reduces the local-repair argument to anticoncentration estimates for uniformly random nested chains; its proof is deferred to Appendix~\ref{app:local-repair}.
\begin{lemma}\label{lem:nested-to-local-remote}
Fix $\beta,B>0$ and an integer $D>3/\beta$. There exists $n_0=n_0(\beta,D,B)$ such that the following holds. Let $G$ be a finite abelian group and $S\subseteq G\setminus\{0\}$ have size $n\ge n_0$. Let $(X_1,\ldots,X_n)$ be a uniformly random ordering of $S$. For $F\subseteq[n]$, let $\mathcal F_F=\sigma(X_j:j\in F)$ and write $S_F=S\setminus\{X_j:j\in F\}$.
Suppose that for every $F\subseteq[n]$ with $|F|\le100D^2$ and every $1\le h\le D$, almost surely
\begin{equation}\label{eq:nested-chain-sum}
\sum_{1\le m_1<\cdots<m_h<|S_F|}\sup_{z_1,\ldots,z_h\in G}
\Prob\bigl(\sigma(R_i)=z_i\text{ for all }i\in[h]\mid\mathcal F_F\bigr)
\le Bn^{-\beta h},
\end{equation}
where, for each $(m_1,\ldots,m_h)$, the sets $R_1\subsetneq\cdots\subsetneq R_h\subsetneq S_F$ form a uniformly random nested chain with $|R_i|=m_i$, conditional on $\mathcal F_F$. Then $S$ is rearrangeable.
\end{lemma}

To apply Lemma~\ref{lem:nested-to-local-remote}, we use the following estimate for uniformly random nested chains.
\begin{lemma}\label{lem:block-exposure}
Let $X\subseteq G$ have size $N$, let $1\le h<N$, and let
$b:[N]\to[0,\infty)$ and $\varepsilon\ge0$.
Call a subset $U\subseteq X$ \emph{good} if, for every integer
$1\le m\le |U|/2$, a uniformly random $m$-subset $R$ of $U$
satisfies
$\sup_{z\in G}\Prob(\sigma(R)=z\mid U)\le b(m)$.
Suppose that, for every integer $u$ with
$N/(h+1)\le u\le N$, a uniformly random $u$-subset of $X$
is good with probability at least $1-\varepsilon$.

Fix $1\le m_1<\cdots<m_h<N$, and let
$R_1\subsetneq\cdots\subsetneq R_h\subsetneq X$
be a uniformly random chain with $|R_i|=m_i$.
Set $m_0=0$, $m_{h+1}=N$, and
$r_i=m_i-m_{i-1}$ for $i\in[h+1]$.
Choose $j$ with $r_j=\max_i r_i$. Then
\begin{equation}\label{eq:block-exposure-product}
\sup_{z_1,\ldots,z_h\in G}
\Prob\bigl(\sigma(R_i)=z_i\text{ for all }i\in[h]\bigr)
\le
\prod_{\substack{i\in[h+1]\\i\ne j}}b(r_i)
+h\varepsilon.
\end{equation}
Thus, we have
\begin{equation}\label{eq:block-exposure-sum}
\sum_{1\le m_1<\cdots<m_h<N}
\sup_{z_1,\ldots,z_h\in G}
\Prob\bigl(\sigma(R_i)=z_i\text{ for all }i\in[h]\bigr)
\le
(h+1)\left(\sum_{m=1}^N b(m)\right)^h
+h\binom{N-1}{h}\varepsilon,
\end{equation}
where each chain in the sum has $|R_i|=m_i$.
\end{lemma}

\begin{proof}
Set $R_0=\varnothing$ and $R_{h+1}=X$, and write
$C_i=R_i\setminus R_{i-1}$ for $i\in[h+1]$.
Then $(C_1,\ldots,C_{h+1})$ is a uniformly random
ordered partition of $X$ with $|C_i|=r_i$.
Fix $z_1,\ldots,z_h\in G$, and set $z_0=0$ and
$z_{h+1}=\sigma(X)$. Since
$\sigma(R_i)=\sum_{s=1}^i\sigma(C_s)$, the conditions
$\sigma(R_i)=z_i$ for $i\in[h]$ are equivalent to
$
\sigma(C_i)=z_i-z_{i-1}$, $i\in[h+1]\setminus\{j\}$.

We expose the sets $C_i$, $i\ne j$, one by one,
leaving $C_j$ unexposed. At each step, conditional
on the previously exposed sets, the next set
$C_i$ is a uniformly random $r_i$-subset of the
remaining set $U$.

Since $C_j$ remains unexposed and $r_j$ is maximal,
we have
$
|U|\ge r_i+r_j\ge2r_i$,
and
$|U|\ge r_j\ge\frac{N}{h+1}$.
Consequently, whenever $U$ is good, the
conditional probability that
$\sigma(C_i)=z_i-z_{i-1}$ is at most $b(r_i)$ by assumption.

For each step, the remaining set $U$ is uniformly
distributed among subsets of $X$ of its size.
By our assumption, the probability that $U$
is not good is at most $\varepsilon$.
As there are $h$ steps, the probability that
some remaining set is not good is at most
$h\varepsilon$. Thus, with probability
at least $1-h\varepsilon$, all remaining sets are
good. Successive conditioning then gives
\[
\Prob\bigl(\sigma(R_i)=z_i
\text{ for all }i\in[h]\bigr)
\le\prod_{i\ne j}b(r_i)+h\varepsilon.
\]
This proves \eqref{eq:block-exposure-product}.

To obtain \eqref{eq:block-exposure-sum}, we sum
over all $1\le m_1<\cdots<m_h<N$.
For each fixed $j\in[h+1]$, dropping the
restrictions on the block sizes gives
\[
\sum_{\substack{r_1,\ldots,r_{h+1}\ge1\\
\sum_i r_i=N,\ r_j=\max_i r_i}}
\prod_{i\ne j}b(r_i)
\le\left(\sum_{m=1}^N b(m)\right)^h.
\]
Since there are $h+1$ choices for $j$ and
$\binom{N-1}{h}$ choices for the chain sizes,
we obtain \eqref{eq:block-exposure-sum}.
\end{proof}

\begin{lemma}\label{lem:large-subgroup-core}
Fix $C_0>0$, $0<\alpha<1/4$, and $\delta>0$. Let $n$ be sufficiently large in terms of these constants, and let $S\subseteq G\setminus\{0\}$ have size $n$ with $\diss(S)\le C_0(\log n)^2$. Suppose that there exist $a\in G$ and a subgroup $H\le G$ such that, writing $A=S\cap(a+H)$, the following hold:
\begin{enumerate}[label=\textup{(\roman*)}]
    \item $|S\setminus A|\le n^{2\alpha}$, and
$|A\setminus(x+K)|\ge n^\alpha$
for every proper subgroup $K<H$ and every $x\in G$;
    \item $|H|\ge n^{1+\delta}$;
\end{enumerate}
Then $S$ is rearrangeable.
\end{lemma}
\begin{proof}
Fix an integer $D>\max\{3/\delta,6/\alpha\}$ and choose $\beta$ with $3/D<\beta<\min\{\delta,\alpha/2\}$. Set $k=\lceil n^\alpha\rceil$. By Lemma~\ref{lem:nested-to-local-remote},
it suffices to establish \eqref{eq:nested-chain-sum}
with $B=D+2$.

Fix $F\subseteq[n]$ with $|F|\le100D^2$, condition on $\mathcal F_F$, and write $N=|S_F|=n-|F|$. Then $N\ge n/2$ and $|(A\cap S_F)\setminus(x+J)|\ge k-|F|\ge k/2$ for every $J<H$ and $x\in G$. Let $T=|\mathcal T(S)|$. By Theorem~\ref{lem:coset-traces}\textup{(ii)} and our choice of $k$, we have $\log(2nT)=o(k)$.

Let $n/(2(D+1))\le u\le N$, let $U$ be a uniform $u$-subset of $S_F$, and write $B=U\cap A$. Call $U$ \emph{good} if $|B\setminus(x+J)|\ge uk/(4N)$ for every $J<H$ and $x\in G$. We first prove the following claim.

\medskip
\begin{claim}\label{clm:good-pool}
There are constants $c_D,C_D>0$, depending only on $D$, such that, writing
$b(m)=|H|^{-1}+C_D/\sqrt{nkm}$ for $1\le m\le n$, the following hold:
\begin{enumerate}[label=\textup{(\roman*)}]
\item $\Prob(U\text{ is not good})\le e^{-c_Dk}$;
\item If $U$ is good and $R$ is a uniform $m$-subset of $U$, where $1\le m\le u/2$, then
\[
\sup_{z\in G}\Prob(\sigma(R)=z\mid U)
\le b(m).
\]
\end{enumerate}
\end{claim}
\medskip

Since $|S\setminus A|\le n^{2\alpha}$, we have $|B|\ge u-n^{2\alpha}\ge3u/4$. For every $J<H$ and $x\in G$,
$
\mathbb E|B\setminus(x+J)|
=\frac{u}{N}|(A\cap S_F)\setminus(x+J)|
\ge\frac{uk}{2N}.
$
By Lemma~\ref{lem:hypergeometric} and a union bound over the at most $T$ coset traces, we have
\[
\Prob(U\text{ is not good})
\le T\exp\left(-\frac{uk}{16N}\right)
\le e^{-c_Dk},
\]
where $c_D=1/(64(D+1))$. This proves \textup{(i)}.

If $U$ is good, set $\rho=k/(8n)\le1/4$ and $\kappa=\rho|B|$. Then $|(B-a)\setminus(x+J)|\ge uk/(4N)\ge\kappa$ for every $J<H$ and $x\in H$. Moreover, $|\mathcal T(B-a)|\le T$ and $\kappa\ge3k/(64(D+1))\ge C\log(2|B|T)$ for sufficiently large $n$. By Lemma~\ref{lem:one-slice-traces},
\[
\sup_{z\in G}\Prob(\sigma(R)=z\mid U)
\le\frac1{|H|}+\frac{C}{u\sqrt{\rho m}}
\le b(m),
\]
proving \textup{(ii)}. 
\hfill $\blacksquare$ \medskip

Fix $1\le h\le D$.
Since $N/(h+1)\ge n/(2(D+1))$, by
Claim~\ref{clm:good-pool} we may apply
Lemma~\ref{lem:block-exposure}, conditionally on
$\mathcal F_F$, with $X=S_F$,
$\varepsilon=e^{-c_Dk}$, and
$b(m)=|H|^{-1}+C_D/\sqrt{nkm}$.
We obtain
\begin{equation}\label{eq:summed-chain}
\begin{aligned}
&\sum_{1\le m_1<\cdots<m_h<N}
\sup_{z_1,\ldots,z_h\in G}
\Prob\bigl(
  \sigma(R_i)=z_i\text{ for all }i\in[h]
  \mid\mathcal F_F
\bigr)\\
&\qquad\le
(h+1)\left(\sum_{m=1}^N b(m)\right)^h
+h\binom{N-1}{h}e^{-c_Dk}\\
&\qquad\le
(h+1)\left(
  \frac n{|H|}+\frac{2C_D}{\sqrt k}
\right)^h
+Dn^he^{-c_Dk},
\end{aligned}
\end{equation}
where we used
$\sum_{m=1}^N m^{-1/2}\le2\sqrt n$
and $h\binom{N-1}{h}\le Dn^h$.
Since $|H|\ge n^{1+\delta}$ and $k\ge n^\alpha$,
$
\frac n{|H|}+\frac{2C_D}{\sqrt k}
\le n^{-\delta}+2C_Dn^{-\alpha/2}
=o(n^{-\beta})$,
and
$Dn^De^{-c_Dk}=o(n^{-\beta D})$.
Thus \eqref{eq:summed-chain} is at most $(D+2)n^{-\beta h}$ for every $h\in[D]$, proving \eqref{eq:nested-chain-sum}.
\end{proof}

\section{A large-order coset}\label{sec:large-order-coset}
We next consider a coset $a+H$ of large order in $G/H$. Our main result in this section is the following.
\begin{theorem}\label{thm:large-order-absorption}
Fix $0<\varepsilon<\delta\le1$. There exists $n_0(\varepsilon,\delta)$ such that if $|S\setminus(a+H)|\le n^\varepsilon$ and $\ord(a+H)\ge n^\delta$, then every $n$-element set $S\subseteq G\setminus\{0\}$ with $n\ge n_0$ is rearrangeable.
\end{theorem}

In any ordering of $X\subseteq a+H$, an interval of length $m$ has sum in $ma+H$. Thus an interval with sum $z$ must satisfy $m(a+H)=z+H$, restricting the possible lengths to at most one residue class modulo $\ord(a+H)$. We use this restriction to construct valid orderings whose partial sums avoid a prescribed set, adapting the local adjustment argument of Pham and Sauermann~\cite{PhamSauermann}.

\begin{lemma}\label{lem:large-order-forbidden}
For every integer $D\ge2$, there are constants $C_D,N_D$ with the following property. Let $G$ be a finite abelian group, $H\le G$, and $X\subseteq(a+H)\setminus\{0\}$ with $|X|=N\ge N_D$. Write $q=\ord(a+H)$ in $G/H$ and $\theta=q^{-1}+N^{-1}$. Let $F\subseteq G$ be finite with $\sigma(X)\notin F$. If
\begin{equation}\label{eq:large-order-criterion}
C_D\bigl(\theta+N\theta^D\bigr)+2^{10D+1}|F|\theta<1,
\end{equation}
then $X$ has a valid ordering $(x_1, \dots, x_N)$ such that $x_1 + \ldots + x_k \not\in F$ for $1 \le k \le N$.
\end{lemma}

\begin{proof}
Our strategy is to repair repeated partial sums via disjoint swaps between nearby elements, ensuring that the updated partial sums avoid $F$.
Fix $L = 4D$, and let $K_D$ and $N_D \ge 2L+2$ be sufficiently large constants. Let $\pi \colon [N] \to X$ be a uniformly random bijection. For an ordering $\rho$ of $X$ and an index subset $J \subseteq [N]$, we denote $\Sigma_\rho(J) = \sum_{j \in J} \rho(j)$. We have the following key observation.

\begin{observation}
Every zero-sum interval has length divisible by $q$.
\end{observation}

We may assume that $q > L + 1$; indeed, if $q \le L + 1$, then $C_D \theta \ge C_D/q \ge 1$, which contradicts~\eqref{eq:large-order-criterion}.
If $A$ is a uniformly random $m$-subset of $v$ distinct group elements, where $1\le m<v$, exposing all but one element of $A$ or its complement gives
\begin{equation}\label{eq:large-order-atom}
\sup_{z\in G}\Prob(\sigma(A)=z)\le\frac{1}{\max\{m,v-m\}+1}\le\frac2v.
\end{equation}
Indeed, conditioning on all but one element of $A$ (resp.\ its complement) leaves the final element chosen uniformly from at least $v-m+1$ (resp.\ $m+1$) elements, which takes the unique required value with probability at most $1/(\max\{m,v-m\}+1)\le 2/v$.

Fix $W\subseteq[N]$ with $|W|\le L+1$, and condition on $\pi|_W$. Write $U=[N]\setminus W$ and $u=|U|\ge N/2$. Let $\varnothing\ne R_1\subsetneq\cdots\subsetneq R_h\subsetneq U$, where $1\le h\le D$, and let $Z_1,\ldots,Z_h\in G$ be arbitrary values depending on $\pi(j)$, $j\in W$. We bound the conditional probability that $\Sigma_\pi(R_i)=Z_i$ for every $i\in[h]$.
By \eqref{eq:large-order-atom}, every subset of $X\setminus\pi(W)$ of size at least $u/(h+1)$ is good with $b(m)=2(h+1)/u$. Applying Lemma~\ref{lem:block-exposure} conditionally on $\pi|_W$, with $\varepsilon=0$, we obtain
\begin{equation}\label{eq:large-order-fixed-chain}
\Prob\bigl(\Sigma_\pi(R_i)=Z_i\text{ for all }i\in[h]\mid\pi|_W\bigr)
\le\left(\frac{2(h+1)}u\right)^h
\le\left(\frac{4(D+1)}N\right)^h.
\end{equation}

For an ordering $\rho$ and $1\le i<j\le N$, write $\rho^{ij}$ for the ordering obtained by exchanging the entries at positions $i$ and $j$. A collection $\mathcal Q$ of exchanges is \emph{good} if each position is used at most once and $|i-j|\le L$ for every $(i,j)\in\mathcal Q$. Write $\rho_{\mathcal Q}$ for the ordering after these exchanges, and let $\mathcal Q_b$ be the set of good collections using only positions in $[b,N]$. Define
\[
\mathcal B(\rho)=\{b\in[N]:\Sigma_\rho([s,b])=0\text{ for some }2\le s<b\}.
\]
In fact, $\mathcal B(\rho)$ is the set of right endpoints of zero-sum intervals starting at position $2$ or later.
Since $0\notin X$, the ordering $\rho$ is valid precisely when $\mathcal B(\rho)=\varnothing$. 
Exchanging the entries at $b$ and $z$ changes only the sums of intervals containing exactly one of these positions. We call $z$ a \emph{bad choice} if the exchange creates a new zero-sum interval $[s,t]$ with $2\le s<t\le N$. Thus the set of bad choices is
\[
\mathcal Z_b(\rho)=\left\{z\in[b+1,\min\{N,b+L\}]:
\begin{array}{l}
\exists\,2\le s<t\le N,\ \Sigma_{\rho^{bz}}([s,t])=0,\\
|[s,t]\cap\{b,z\}|=1
\end{array}\right\}.
\]

If $z\notin\mathcal Z_b(\rho)$, every zero-sum interval $[s,t]$ of $\rho^{bz}$ with $2\le s<t$ satisfies $|[s,t]\cap\{b,z\}|\in\{0,2\}$. Hence $\Sigma_\rho([s,t])=0$ and $t\ne b$, giving
\begin{equation}\label{eq:exchange-removes-endpoint}
\mathcal B(\rho^{bz})\subseteq\mathcal B(\rho)\setminus\{b\}.
\end{equation}

We consider three events:
\[
\begin{aligned}
\mathcal E_1&=\{\mathcal B(\pi)\cap[N-L+1,N]\ne\varnothing\},\\
\mathcal E_2&=\{\exists\,b\in[N-L]:|\mathcal B(\pi)\cap[b,b+L]|>D\},\\
\mathcal E_3&=\{\exists\,b\in\mathcal B(\pi)\cap [N-L],\ \exists\,\mathcal Q\in \mathcal Q_b:
|\mathcal Z_b(\pi_{\mathcal Q})|\ge2D\}.
\end{aligned}
\]
The three events are illustrated in Figure~\ref{fig:local-repair-events}.

\begin{figure}[t]
\centering
\resizebox{\linewidth}{!}{%
\begin{tikzpicture}[
  x=1cm,y=1cm,line cap=round,line join=round,font=\footnotesize,
  panel/.style={draw=black!22,fill=black!1,rounded corners=2pt,line width=.45pt},
  axis/.style={draw=black!65,line width=.65pt},
  tick/.style={draw=black!65,line width=.6pt},
  bad/.style={circle,fill=red!75!black,inner sep=1.8pt},
  witness/.style={draw=red!75!black,line width=.85pt},
  blocked/.style={draw=violet!80!black,line width=1pt},
  brace/.style={decorate,decoration={brace,amplitude=3pt},
    draw=black!55,line width=.55pt}
]
\begin{scope}
  \path[panel] (0,0) rectangle (5.15,3.30);
  \node[font=\small] at (2.575,3.00) {$\mathcal E_1$};
  \node at (2.575,2.64) {Bad endpoint near the end};
  \fill[red!8] (2.65,.99) rectangle (4.85,1.31);
  \draw[axis] (.30,1.15)--(4.85,1.15);
  \foreach \a in {.30,.95,2.65,4.85}
    \draw[tick] (\a,1.08)--(\a,1.22);
  \node[bad] at (4.08,1.15) {};
  \node at (.30,.80) {$1$};
  \node at (.95,.80) {$s$};
  \node at (2.65,.80) {$N-L+1$};
  \node at (4.08,.80) {$b$};
  \node at (4.85,.80) {$N$};
  \draw[witness] (.95,1.81)--(.95,2.00)--(4.08,2.00)--(4.08,1.81);
  \node[text=red!75!black] at (2.515,2.29) {$\Sigma_\pi([s,b])=0$};
  \draw[brace] (4.85,.49)--(2.65,.49);
  \node at (3.75,.21) {$L$ positions};
\end{scope}
\begin{scope}[xshift=5.60cm]
  \path[panel] (0,0) rectangle (5.15,3.30);
  \node[font=\small] at (2.575,3.00) {$\mathcal E_2$};
  \node at (2.575,2.64) {Too many bad endpoints};
  \fill[orange!12] (1.05,.99) rectangle (4.10,1.31);
  \draw[axis] (.30,1.15)--(4.85,1.15);
  \foreach \a in {.30,1.05,4.10,4.85}
    \draw[tick] (\a,1.08)--(\a,1.22);
  \foreach \a in {1.45,2.15,3.50}
    \node[bad] at (\a,1.15) {};
  \node[fill=orange!12,inner sep=1pt,text=red!75!black]
    at (2.85,1.15) {$\cdots$};
  \node at (.30,.80) {$1$};
  \node at (1.05,.80) {$b$};
  \node at (4.10,.80) {$b+L$};
  \node at (4.85,.80) {$N$};
  \draw[witness] (1.45,1.81)--(1.45,2.00)--(3.50,2.00)--(3.50,1.81);
  \node[text=red!75!black] at (2.475,2.29) {$\ge D+1$ bad endpoints};
  \draw[brace] (4.10,.49)--(1.05,.49);
  \node at (2.575,.21) {$L+1$ positions};
\end{scope}
\begin{scope}[xshift=11.20cm]
  \path[panel] (0,0) rectangle (5.15,3.30);
  \node[font=\small] at (2.575,3.00) {$\mathcal E_3$};
  \node at (2.575,2.64) {Too many excluded partners};
  \fill[violet!8] (.95,.99) rectangle (3.80,1.31);
  \draw[axis] (.30,1.15)--(4.85,1.15);
  \foreach \a in {.30,1.25,3.80,4.40,4.85}
    \draw[tick] (\a,1.08)--(\a,1.22);
  \fill[black!75] (.65,1.15) circle (1.6pt);
  \foreach \a in {1.70,2.40,3.40}{
    \draw[blocked] (\a-.065,1.085)--(\a+.065,1.215);
    \draw[blocked] (\a-.065,1.215)--(\a+.065,1.085);
  }
  \node[fill=violet!8,inner sep=.5pt,text=violet!80!black]
    at (2.92,1.15) {$\cdots$};
  \node at (.30,.80) {$1$};
  \node at (.73,.80) {$b$};
  \node at (1.25,.80) {$u$};
  \node[text=violet!80!black] at (2.40,.80) {$z$};
  \node at (3.80,.80) {$b+L$};
  \node at (4.40,.80) {$x$};
  \node at (4.85,.80) {$N$};
  \draw[blocked] (1.25,1.81)--(1.25,2.00)--(4.40,2.00)--(4.40,1.81);
  \node[text=violet!80!black] at (2.825,2.29)
    {$\Sigma_{\eta^{bz}}([u,x])=0$};
  \draw[<->,>=Stealth,draw=violet!80!black,densely dashed,line width=.65pt]
    (.65,1.21) .. controls (.90,1.92) and (2.15,1.92) .. (2.40,1.21);
  \node[text=violet!80!black] at (1.525,1.47) {$(b,z)$};
  \draw[brace] (3.80,.49)--(.95,.49);
  \node[text=violet!80!black] at (2.375,.21) {$\ge2D$ excluded partners};
\end{scope}
\end{tikzpicture}%
}
\caption{The three bad events.}
\label{fig:local-repair-events}
\end{figure}
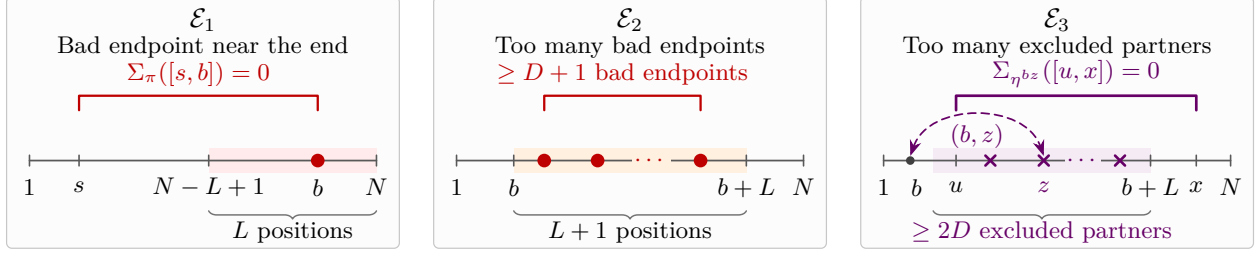

\begin{claim}\label{clm:large-order-event1}
We have $\Prob(\mathcal E_1)\le2L\theta.$
\end{claim}

For fixed $b$, a zero-sum interval $[s,b]$ requires $s\equiv b+1\pmod q$, so there are at most $N/q+1$ choices of $s$. By the union bound and \eqref{eq:large-order-atom},
$
\Prob(\mathcal E_1)
\le L\left(\frac Nq+1\right)\frac2N
=2L\theta.
$ \hfill $\blacksquare$
\medskip

\begin{claim}\label{clm:large-order-event2}
We have $\Prob(\mathcal E_2)\le K_DN\theta^D.$
\end{claim}

If $\mathcal E_2$ occurs, choose $b\in[N-L]$ such that $|\mathcal B(\pi)\cap[b,b+L]|>D$. Then there are $b<b_1<\cdots<b_D\le b+L$ with $b_i\in\mathcal B(\pi)$. For each $i\in[D]$, choose $s_i\ge2$ such that $\Sigma_\pi([s_i,b_i])=0$. Since $q\mid b_i-s_i+1$ and $q>L+1$, we have $s_i\le b_i-q+1<b$. Moreover, the $s_i$ are distinct: otherwise $q\mid b_i-b_j$ for some $i\ne j$, a contradiction to $0<|b_i-b_j|\le L<q$.

Fix these endpoints and condition on $\pi|_{[b,b+L]}$. The intervals $[s_i,b-1]$ can be ordered to form a strict chain of nonempty proper subsets of the unexposed positions. The zero-sum conditions become $\Sigma_\pi([s_i,b-1])=-\Sigma_\pi([b,b_i])$, where the right-hand sides are fixed by $\pi|_{[b,b+L]}$.
Thus \eqref{eq:large-order-fixed-chain} bounds the conditional probability by $(4(D+1)/N)^D$. There are at most $N$ choices of $b$, at most $\binom LD$ choices of $b_1,\ldots,b_D$, and at most $N/q+1$ choices of each $s_i$. Hence
$
\Prob(\mathcal E_2)
\le N\binom LD\left(\frac Nq+1\right)^D
\left(\frac{4(D+1)}N\right)^D
\le K_DN\theta^D.
$
\hfill $\blacksquare$ 
\medskip

Fix $3\le b\le N-L$, $1\le h\le D$, distinct $z_1,\ldots,z_h\in[b+1,b+L]$, and $u_i\in[b+1,z_i]$. We first bound the probability that, after some collection of exchanges from $\mathcal Q_b$, each individual exchange $(b,z_i)$ produces a zero-sum interval $[u_i,x_i]$, where $b+L<x_1<\cdots<x_h\le N$. More precisely, we claim that
\begin{equation}\label{eq:large-order-repair-chain}
\Prob\left(
\begin{array}{l}
\exists\,b+L<x_1<\cdots<x_h\le N,\ \exists\,\mathcal Q\in\mathcal Q_b:\\
\Sigma_{(\pi_{\mathcal Q})^{bz_i}}([u_i,x_i])=0
\text{ for all }i\in[h]
\end{array}
\right)
\le K_D\theta^h.
\end{equation}

Fix $x_1,\ldots,x_h$. The $i$th equation is equivalent to $\Sigma_{\pi_{\mathcal Q}}(([u_i,x_i]\setminus\{z_i\})\cup\{b\})=0$. An exchange affects this sum only if exactly one of its endpoints belongs to the indicated set. Since every exchange has length at most $L$, all exchanges affecting at least one equation have both endpoints in $[b,b+2L]\cup\bigcup_{i=1}^h[x_i-L,x_i+L]$, intersected with $[N]$. This union contains at most $(h+1)(2L+1)$ positions. We may therefore discard all other exchanges and take a union bound over at most $K_D$ fixed collections.

For each fixed collection, let $\eta$ be the ordering obtained by performing its exchanges on $\pi$. Then $\eta$ is uniformly random. Conditional on $\eta|_{[b,b+L]}$, the required equations become $\Sigma_\eta([b+L+1,x_i])=-\Sigma_\eta([u_i,b+L])-\eta(b)+\eta(z_i)$ for $i\in[h]$, with fixed right-hand sides. The intervals on the left form a strict chain of nonempty proper subsets of the unexposed positions; in particular, position $1$ remains outside all of them. By \eqref{eq:large-order-fixed-chain}, their conditional probability is at most $(4(D+1)/N)^h$.

Finally, a zero-sum interval $[u_i,x_i]$ requires $q\mid x_i-u_i+1$, so each $x_i$ has at most $N/q+1$ choices. Taking the union bound over the endpoints and the fixed collections gives an upper bound of $K_D(N/q+1)^h(4(D+1)/N)^h$. Since $(N/q+1)/N=\theta$, this proves \eqref{eq:large-order-repair-chain} after enlarging $K_D$.

\begin{claim}\label{clm:large-order-event3}
We have $\Prob(\mathcal E_3)\le K_DN\theta^D.$
\end{claim}

Suppose that $\mathcal E_3$ occurs. Choose $b\in\mathcal B(\pi)\cap[N-L]$ and $\mathcal Q\in\mathcal Q_b$ with $|\mathcal Z_b(\pi_{\mathcal Q})|\ge2D$, and write $\eta=\pi_{\mathcal Q}$. Since $q>L+1$, every zero-sum witness extends to the left of $b$ or to the right of $b+L$. Thus at least $D$ distinct excluded positions $z_i$ have witnesses of the same type: either $[u_i,x_i]$ with $b<u_i\le z_i\le b+L<x_i$, or $[x_i,v_i]$ with $2\le x_i<b\le v_i<z_i\le b+L$.

In either case the remote endpoints $x_i$ are distinct. In the first case, $x_i=x_j$ implies $q\mid u_i-u_j$, hence $u_i=u_j$ since $|u_i-u_j|<q$. Subtracting the two zero-sum equations then yields $\eta(z_i)=\eta(z_j)$, a contradiction. In the second case, $x_i=x_j$ similarly implies $v_i=v_j$ and the same contradiction.

Order the witnesses by their remote endpoints. By \eqref{eq:large-order-repair-chain} and its left-going analogue, summing over the two directions, at most $N$ choices of $b$, and at most $L^{2D}$ choices of local positions yields
$
\Prob(\mathcal E_3)\le2NL^{2D}K_D\theta^D
\le K_DN\theta^D.
$
\hfill $\blacksquare$
\medskip

By Claims~\ref{clm:large-order-event1}--\ref{clm:large-order-event3},
\begin{equation}\label{eq:repair-bad-events}
\Prob(\mathcal E_1\cup\mathcal E_2\cup\mathcal E_3)
\le K_D(\theta+N\theta^D).
\end{equation}

We now show that any ordering outside these three events can be repaired. Write $\mathcal B(\pi)=\{b_1>\cdots>b_m\}$ and set $\mathcal Q_0=\varnothing$. We construct good collections $\mathcal Q_i$ such that
$\mathcal B(\pi_{\mathcal Q_i})\subseteq\{b_{i+1},\ldots,b_m\}$.

Suppose that $\mathcal Q_{i-1}$ has been constructed. Since every previous exchange partner lies outside $\mathcal B(\pi)$, the position $b_i$ has not been used. Set $C_i=[b_i+1,b_i+L]\subseteq[N]$, where the inclusion follows from $\mathcal E_1^c$. Since $\mathcal E_2$ does not occur,
$|C_i\cap\mathcal B(\pi)|\le D-1$.
Also $|C_i\cap\{z_j:j<i\}|\le D-1$, since $b_i<b_j<z_j\le b_i+L$ whenever $z_j\in C_i$. The previous exchanges belong to $\mathcal Q_{b_i}$, so $\mathcal E_3^c$ implies
$|\mathcal Z_{b_i}(\pi_{\mathcal Q_{i-1}})|\le2D-1$.
Hence at most $4D-3<L$ positions of $C_i$ are excluded.

Choose $z_i\in C_i$ outside these three sets, and let $\mathcal Q_i=\mathcal Q_{i-1}\cup\{(b_i,z_i)\}$. Then $\mathcal Q_i$ is good, and \eqref{eq:exchange-removes-endpoint} implies
$
\mathcal B(\pi_{\mathcal Q_i})
\subseteq\mathcal B(\pi_{\mathcal Q_{i-1}})\setminus\{b_i\}
\subseteq\{b_{i+1},\ldots,b_m\}.
$
Thus $\mathcal B(\pi_{\mathcal Q_m})=\varnothing$, so the resulting ordering is valid.

It remains to ensure that the partial sums avoid $F$ after the repair. Let $\mathcal E_F=\{\exists\,\mathcal Q\in\mathcal Q_1,\ \exists\,j\in[N-1]:\Sigma_{\pi_{\mathcal Q}}([j])\in F\}$. Fix $1\le j<N$. Only exchanges crossing the boundary between $j$ and $j+1$ can change the sum on $[j]$. Since the exchanges are disjoint and have length at most $L$, every modified prefix sum equals $\Sigma_\pi(J)$ for some $J$ satisfying $|J|=j$ and $J\triangle[j]\subseteq[j-L+1,j+L]\cap[N]$. There are at most $2^{2L}$ such sets $J$. Moreover, $\Sigma_\pi(J)=f$ requires $j(a+H)=f+H$, which restricts $j$ to at most $N/q+1$ values for each $f\in F$. By \eqref{eq:large-order-atom} and the union bound,
\begin{equation}\label{eq:forbidden-prefixes}
\Prob(\mathcal E_F)
\le |F|\left(\frac Nq+1\right)2^{2L}\frac2N
=2^{8D+1}|F|\theta.
\end{equation}

Choose $C_D\ge\max\{L+1,K_D\}$ sufficiently large to cover the preceding estimates. By \eqref{eq:large-order-criterion}, \eqref{eq:repair-bad-events} and \eqref{eq:forbidden-prefixes}, we have $\Prob(\mathcal E_1\cup\mathcal E_2\cup\mathcal E_3\cup\mathcal E_F)\le C_D(\theta+N\theta^D)+2^{8D+1}|F|\theta<1$. Hence some initial ordering avoids all four events. The repair procedure produces a valid ordering, and $\mathcal E_F^c$ ensures that all its proper nonempty partial sums avoid $F$. Its full sum is $\sigma(X)\notin F$, completing the proof.
\end{proof}

We next use a greedy absorption argument to incorporate the exceptional colors into a short rainbow path. By inserting an auxiliary color before each exceptional color, we can extend the path while avoiding previously visited vertices and a prescribed forbidden set. The following lemma formalizes this robust path extension.

\begin{lemma}\label{lem:greedy-path-extension}
Let $G$ be a finite abelian group, let $P$ be a rainbow directed path in $\Cay(G,G\setminus\{0\})$, and let $A,E\subseteq G\setminus\{0\}$ be disjoint sets of colors unused by $P$. Suppose that $Z\subseteq G$ and $|A|\ge2(|V(P)|+|Z|)+5|E|$. Then, for any prescribed ordering of $E$, there exists a rainbow directed path $P'$ having $P$ as an initial subpath such that
\begin{enumerate}[label=\textup{(\roman*)}]
\item every color of $E$ occurs on a new edge, and these colors appear along $P'$ in the prescribed order;
\item all other new edges have colors in $A$, and there are between $1$ and $|E|+1$ such edges;
\item $(V(P')\setminus V(P))\cap Z=\varnothing$.
\end{enumerate}
\end{lemma}
\begin{proof}
If $E=\varnothing$, let $v$ be the endpoint of $P$. Since $|A|>|V(P)|+|Z|$, choose $a\in A$ with $v+a\notin V(P)\cup Z$ and append the edge of color $a$.

Otherwise, write $E=\{e_1,\ldots,e_m\}$ in the prescribed order and extend $P$ successively by color pairs $a_i,e_i$, choosing distinct $a_i\in A$. At step $i$, let $v$ be the current endpoint and $W$ the union of the current vertex set and $Z$. Choose $a_i\notin\{a_1,\ldots,a_{i-1}\}\cup(W-v)\cup(W-v-e_i)$. This excludes at most $2|W|+i-1\le2(|V(P)|+|Z|)+5(i-1)<|A|$ colors. Both $v+a_i$ and $v+a_i+e_i$ avoid $W$, and they are distinct since $e_i\ne0$. After $m$ steps, the required extension uses exactly $m$ colors of $A$.
\end{proof}

Now we prove Theorem \ref{thm:large-order-absorption}.
\begin{proof}[Proof of Theorem \ref{thm:large-order-absorption}]
Let $A=S\cap(a+H)$, $E=S\setminus A$, $m=|E|$, and $T=\sigma(S)$. If $m=0$, apply Lemma~\ref{lem:large-order-forbidden} with $F=\varnothing$. Otherwise, $|A|\ge5m+4$ for sufficiently large $n$. By Lemma~\ref{lem:greedy-path-extension}, applied with $P=(0)$ and $Z=\{T\}$, choose a prefix using every color of $E$ and at most $m+1$ colors of $A$, whose nonempty partial sums are distinct and avoid $0$ and $T$.

Let $W$ be its nonempty partial sums, $v$ its endpoint, and $X$ the unused part of $A$. Set $F=W-v$ and $N=|X|$. Then $N\ge n-2m-1\ge n/2$, $|F|\le2m+1\le3m$, and $\sigma(X)=T-v\notin F$. Also, $\ord(a+H)\ge N^\delta$ and $|F|\le3\cdot2^\varepsilon N^\varepsilon$. By Lemma~\ref{lem:large-order-forbidden}, $X$ has a valid ordering with nonempty partial sums $q_1,\ldots,q_N\notin F$. Appending it produces the pairwise distinct partial sums $W\cup\{v+q_1,\ldots,v+q_N\}$.
\end{proof}

\section{Rainbow paths via expansion and absorption}\label{sec:dense-core}
We now consider the remaining case in which both $|H|$ and $q=\ord_{G/H}(a+H)$ are small. The subgroup $Q=H+\langle a\rangle$ then has order $q|H|$ and contains almost all of $S$, so $S\cap Q$ is dense in $Q$. We address this case by proving Theorem~\ref{lem:dense-core-sparse-outside}.

\begin{theorem}\label{lem:dense-core-sparse-outside}
Let $\gamma=2^{-20}$. There is an absolute integer $h_0$ such that the following holds. Let $G$ be a finite abelian group, let $Q\le G$ have order $h\ge h_0$, and let $A\subseteq Q\setminus\{0\}$ and $E\subseteq G\setminus Q$. If $|A|\ge h^{1-\gamma}$ and $|E|\le h^{1-2\gamma}$, then $A\cup E$ is rearrangeable. 
\end{theorem}

Following the expansion approach of~\cite{BucicEtAl}, we first use a Cayley graph diameter bound to find a subgroup of $Q$ containing almost all colors of $A$, such that the Cayley graph generated by these colors has no sparse cut. We then construct a short rainbow path through the colors outside this subgroup and adapt the absorption argument of Bedert, Buci\'c, Kravitz, Montgomery and M\"uyesser~\cite[Section~9]{BBKMM} to obtain a rainbow path using all but one color of $A\cup E$.

\subsection{Passing to a subgroup with no sparse cut}
We first bound the diameter of a connected Cayley graph in terms of its order and the size of its generating set.
\begin{lemma}\label{lem:cayley-diameter}
Let $K$ be a finite abelian group of order $M$, and let $B\subseteq K\setminus\{0\}$ generate $K$. Then the underlying undirected graph of \(\Cay(K,B)\) has diameter less than \(3M/(|B|+1)\).
\end{lemma}

\begin{proof}
Take a shortest path $v_0\ldots v_D$ whose length is the diameter $D$. The closed neighborhoods of $v_0,v_3,\ldots,v_{3\lfloor D/3\rfloor}$ are pairwise disjoint, each of size at least $|B|+1$. Hence $(\lfloor D/3\rfloor+1)(|B|+1)\le M$, and the result follows from $D<3(\lfloor D/3\rfloor+1)$.
\end{proof}

We next find a subgroup containing almost all of $A$ such that the Cayley digraph generated by these elements has no sparse cut.
\begin{lemma}\label{lem:dense-expanding-subgroup}
Let $Q$ be a finite abelian group of order $h$, let $A\subseteq Q\setminus\{0\}$ have size $u\ge1$, and let $0<\varepsilon\le1/2$. Then there is a subgroup $H\le Q$ such that $|A\cap H|\ge(1-\varepsilon)u$ and $\Cay(H,A\cap H)$ has no $\zeta$-sparse cut, where $\zeta=\varepsilon u^2/(128h^2)$.
\end{lemma}

\begin{proof}
Set $K_0=Q$. We construct a descending sequence of subgroups as follows. Given $K_i$, write $M_i=|K_i|$, $A_i=A\cap K_i$, and $d_i=|A_i|$. If $d_i<u/2$ or $\Cay(K_i,A_i)$ has no $\zeta$-sparse cut, we stop. Otherwise, since the digraph has equal indegrees and outdegrees, there is a set $X\subsetneq K_i$ with $0<|X|\le M_i/2$ and $e(X,K_i\setminus X)<\zeta|X|(M_i-|X|)$.

Let $B_i=\{a\in A_i:|X\triangle(X+a)|\le |X|d_i/(12M_i)\}$. Since $\sum_{a\in A_i}|X\triangle(X+a)|=2e(X,K_i\setminus X)<2\zeta|X|M_i$, we have $|A_i\setminus B_i|<24\zeta M_i^2/d_i\le 3d_i/8$, where the last inequality follows from $\zeta\le u^2/(256h^2)$, $M_i\le h$, and $d_i\ge u/2$. In particular, $|B_i|>d_i/2$.

Suppose that $\langle B_i\rangle=K_i$. By Lemma~\ref{lem:cayley-diameter}, every $g\in K_i$ can be expressed as a sum of fewer than $6M_i/d_i$ elements of $B_i\cup(-B_i)$. Since $|X\triangle(X-b)|=|X\triangle(X+b)|$ for every $b\in K_i$, the triangle inequality implies $|X\triangle(X+g)|<(6M_i/d_i)\cdot(|X|d_i/(12M_i))=|X|/2$ for every $g\in K_i$. On the other hand, $\sum_{g\in K_i}|X\cap(X+g)|=|X|^2$, since each ordered pair $(x,y)\in X^2$ contributes once, with $g=x-y$. Using $|X\triangle(X+g)|=2|X|-2|X\cap(X+g)|$, we obtain
\[
\frac1{M_i}\sum_{g\in K_i}|X\triangle(X+g)|
=2|X|-\frac{2|X|^2}{M_i}
=2|X|\left(1-\frac{|X|}{M_i}\right)
\ge |X|,
\]
where the last inequality follows from $|X|\le M_i/2$. This contradicts $|X\triangle(X+g)|<|X|/2$ for every $g\in K_i$.
Thus $K_{i+1}=\langle B_i\rangle$ is a proper subgroup of $K_i$, so $M_{i+1}\le M_i/2$. Moreover, $d_i-d_{i+1}\le|A_i\setminus B_i|<48\zeta M_i^2/u$.

Since the subgroup orders strictly decrease, the process terminates, say at $K_s$. As $M_i\le h/2^i$, the total number of discarded colors satisfies
$
u-d_s<\frac{48\zeta}{u}\sum_{i<s}M_i^2\le\frac{64\zeta h^2}{u}=\frac{\varepsilon u}{2}.
$
In particular, $d_s>(1-\varepsilon/2)u\ge u/2$, so the process cannot terminate as $d_s<u/2$. Consequently, $H=K_s$ satisfies $|A\cap H|\ge(1-\varepsilon)u$, and $\Cay(H,A\cap H)$ has no $\zeta$-sparse cut.
\end{proof}

\subsection{Switching and absorption}

For a directed path $P$, let $\operatorname{col}(P)$ denote its set of edge colors. A $g$-pair is an ordered pair $(a,b)$ of distinct nonzero colors with $a+b=g$. A family of such pairs is color-disjoint if no color belongs to two pairs. For a family $\mathcal P$, write $\operatorname{col}(\mathcal P)=\bigcup_{(a,b)\in\mathcal P}\{a,b\}$. If $g\ne0$, the paths $v,v+a,v+g$ corresponding to these pairs have the same endpoints.

We shall use the following matching result of Montgomery~\cite{MontgomerySpanning}, in the form given by Bedert, Buci\'c, Kravitz, Montgomery and M\"uyesser~\cite[Proposition~9.3]{BBKMM}.

\begin{lemma}[{\cite[Proposition~9.3]{BBKMM}}]\label{lem:robust-pair-matching}
Let $1\le\ell\le t$ be integers. There is a bipartite graph $B$ with bipartition $(X\cup Y,Z)$, where $X$ and $Y$ are disjoint, $|X|=t+\ell$, $|Y|=2t$, and $|Z|=3t$, such that $\Delta(B)\le40$ and, for every $X'\subseteq X$ with $|X'|=t$, there is a perfect matching between $X'\cup Y$ and $Z$.
\end{lemma}

For each starting vertex $x$, we first choose an endpoint $w$ and a color set $F_R$. After removing from $F_R$ the colors in any $\ell$ pairs of $\mathcal P_{\mathrm{flex}}$, the remaining colors are exactly those of a rainbow $x$--$w$ path of length $9t$, with all vertices other than $x$ in $V$. The following lemma ensures that such a choice is possible.
\begin{lemma}\label{lem:switching-chain}
Let $G$ be a finite abelian group, let $H\le G$ have order $N$,
and let $U$ be a coset of $H$. Fix $F\subseteq H\setminus\{0\}$,
$g\in H\setminus\{0\}$, and integers $1\le\ell\le t$.
Let $\mathcal P$ be a family of $3t+\ell$ pairwise color-disjoint
$g$-pairs in $F$, and let
$\mathcal P_{\mathrm{flex}}\subseteq\mathcal P$ have size $t+\ell$.

Let $V$ be a $p$-random subset of $U$, where $0<p\le1$.
Suppose that $N$ is sufficiently large and
$|F|p^{42}\ge2^{20}\max\{3t,\log N\}$.

Then, with probability at least $1-1/N$, for every $x\in U$
there are a vertex $w\in V$ and a color set $F_R$ of size
$9t+2\ell$ satisfying
$\operatorname{col}(\mathcal P)\subseteq F_R\subseteq F$,
with the following property.
For every $\ell$-element subfamily
$\mathcal P'\subseteq\mathcal P_{\mathrm{flex}}$,
there is a rainbow directed $x$--$w$ path using exactly the colors
in $F_R\setminus\operatorname{col}(\mathcal P')$,
with all vertices other than $x$ in $V$.
\end{lemma}

\begin{proof}
Let $B$ be the bipartite graph from Lemma~\ref{lem:robust-pair-matching}, with bipartition $(X\cup Y,Z)$. Label $X$ and $Y$ by the pairs in $\mathcal P_{\mathrm{flex}}$ and $\mathcal P\setminus\mathcal P_{\mathrm{flex}}$, respectively. Write $Z=\{z_1,\ldots,z_{3t}\}$, and let $\mathcal P_j$ consist of the pairs labelling $N_B(z_j)$. Then $|\mathcal P_j|\le40$, and for every $\ell$-element $\mathcal P'\subseteq\mathcal P_{\mathrm{flex}}$, there is a bijection $j\mapsto p_j$ from $[3t]$ to $\mathcal P\setminus\mathcal P'$ with $p_j\in\mathcal P_j$.

Fix $z\in U$ and $j\in[3t]$, and define
$
D_j=\{0,g\}\cup\{a:(a,b)\in\mathcal P_j\}.
$
For $c\in F\setminus\operatorname{col}(\mathcal P)$, consider the walks
$
z\xrightarrow{c}z+c\xrightarrow{a}z+c+a
\xrightarrow{b}z+c+g$,
$(a,b)\in\mathcal P_j.
$
Define the gadget $Z_j(z,c):=\{z+c,z+c+g\}\cup\{z+c+a:(a,b)\in\mathcal P_j\}=z+c+D_j$. We have $|D_j|\le42$. If $z\notin Z_j(z,c)$, each walk is a rainbow path. These walks form the gadget illustrated in Figure~\ref{fig:switching-gadget}.

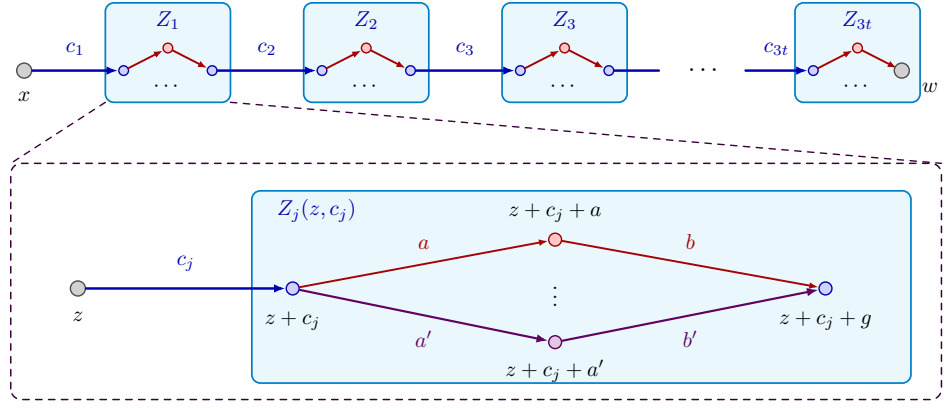
\begin{figure}[t]
\centering
\resizebox{0.75\linewidth}{!}{%
\begin{tikzpicture}[
  x=1cm,y=1cm,font=\small,line cap=round,line join=round,
  vertex/.style={circle,draw=blue!65!black,fill=blue!18,
    minimum size=4.5pt,inner sep=0pt,line width=.65pt},
  endpoint/.style={vertex,draw=black!70,fill=black!18,minimum size=7pt},
  pink/.style={vertex,draw=red!60!black,fill=red!22},
  purple/.style={vertex,draw=violet!70!black,fill=violet!22},
  connector/.style={-{Latex[length=1.9mm]},draw=blue!70!black,line width=1.15pt},
  branch/.style={-{Latex[length=1.5mm]},draw=red!65!black,line width=.9pt},
  alternative/.style={-{Latex[length=1.9mm]},draw=violet!75!black,line width=1pt},
  region/.style={draw=cyan!60!blue,fill=cyan!7,rounded corners=4pt,line width=.8pt}
]
\foreach \shift/\tag/\id in {1.75/1/one,5.00/2/two,8.25/3/three,13.05/{3t}/last}{
  \begin{scope}[xshift=\shift cm]
    \path[region] (0,4.65) rectangle (2.05,6.28);
    \node[text=blue!65!black] at (1.025,5.99) {$Z_{\tag}$};
    \node[vertex] (L-\id) at (.30,5.17) {};
    \node[pink] (T-\id) at (1.025,5.54) {};
    \node[vertex] (R-\id) at (1.75,5.17) {};
    \draw[branch] (L-\id)--(T-\id);
    \draw[branch] (T-\id)--(R-\id);
    \node at (1.025,4.92) {$\cdots$};
  \end{scope}
}
\node[endpoint] (x) at (.42,5.17) {};
\draw[connector] (x)--node[above=3pt,text=blue!70!black] {$c_1$} (L-one);
\draw[connector] (R-one)--node[above=3pt,text=blue!70!black] {$c_2$} (L-two);
\draw[connector] (R-two)--node[above=3pt,text=blue!70!black] {$c_3$} (L-three);
\draw[draw=blue!70!black,line width=1.15pt] (R-three)--(10.82,5.17);
\node[font=\large] at (11.58,5.17) {$\cdots$};
\draw[connector] (12.24,5.17)--node[above=3pt,text=blue!70!black] {$c_{3t}$} (L-last);
\node[anchor=north] at (.42,4.96) {$x$};
\node[endpoint] at (R-last) {};
\node[anchor=north west] at (15.00,5.10) {$w$};

\draw[dashed,draw=violet!45!black,line width=.65pt]
  (1.75,4.65)--(.20,3.65) (3.80,4.65)--(15.45,3.65);
\draw[dashed,draw=violet!45!black,rounded corners=5pt,line width=.7pt]
  (.20,-.22) rectangle (15.45,3.65);
\path[region] (4.15,.05) rectangle (14.95,3.20);
\node[anchor=west,text=blue!65!black] at (4.43,2.88) {$Z_j(z,c_j)$};
\node[endpoint] (z) at (1.30,1.60) {};
\node[vertex,minimum size=6pt] (fork) at (4.82,1.60) {};
\node[pink,minimum size=6pt] (a) at (9.12,2.40) {};
\node[purple,minimum size=6pt] (aprime) at (9.12,.72) {};
\node[vertex,minimum size=6pt] (finish) at (13.55,1.60) {};
\draw[connector] (z)--node[above=4pt,text=blue!70!black] {$c_j$} (fork);
\draw[branch] (fork)--node[above=3pt,text=red!65!black] {$a$} (a);
\draw[branch] (a)--node[above=3pt,text=red!65!black] {$b$} (finish);
\draw[alternative] (fork)--node[below=3pt,text=violet!75!black] {$a'$} (aprime);
\draw[alternative] (aprime)--node[below=3pt,text=violet!75!black] {$b'$} (finish);
\node at (9.12,1.57) {$\vdots$};
\node[anchor=north] at (1.30,1.36) {$z$};
\node[anchor=north] at (4.82,1.36) {$z+c_j$};
\node[anchor=south] at (9.12,2.55) {$z+c_j+a$};
\node[anchor=north] at (9.12,.57) {$z+c_j+a'$};
\node[anchor=north] at (13.55,1.36) {$z+c_j+g$};
\end{tikzpicture}%
}
\caption{A chain of gadgets and the structure of one gadget.}
\label{fig:switching-gadget}
\end{figure}

Since $|\operatorname{col}(\mathcal P)|\le8t$ and $z\in Z_j(z,c)$ only if $c\in-D_j$, at least $|F|-8t-42\ge|F|/2$ choices of $c$ avoid $z$.
Moreover, if $Z_j(z,c)\cap Z_j(z,c')\ne\varnothing$ then $c'-c\in D_j-D_j$.
Thus each $Z_j(z,c)$ intersects at most $42^2$
such sets. Greedily choose
$C_{z,j}\subseteq F\setminus\operatorname{col}(\mathcal P)$
with
$
|C_{z,j}|\ge\frac{|F|}{2\cdot42^2}
>\frac{|F|}{4096}
$
such that the sets $Z_j(z,c)$, $c\in C_{z,j}$,
are pairwise disjoint and avoid $z$.
Fix these sets for all $(z,j)$ before revealing $V$.

Let
$
Y_{z,j}
=|\{c\in C_{z,j}:Z_j(z,c)\subseteq V\}|.
$
By disjointness, $Y_{z,j}$ is a sum of independent
Bernoulli variables, and
$
\mathbb EY_{z,j}
=|C_{z,j}|p^{|D_j|}
\ge\frac{|F|p^{42}}{4096}
\ge256\max\{3t,\log N\}.
$
By Lemma~\ref{lem:chernoff}\textup{(i)},
$
\Prob(Y_{z,j}<384t)
\le e^{-\mathbb EY_{z,j}/8}
\le N^{-32}.
$
Since $3t<N$, a union bound shows that
$Y_{z,j}\ge384t$ for all $(z,j)$ with probability
at least $1-N^{-30}$. Fix such a set $V$.

Fix $x\in U$ and let $x_0=x$. For
$j=1,\ldots,3t$, choose
$c_j\in C_{x_{j-1},j}$ with
$Z_j(x_{j-1},c_j)\subseteq V$.
At step $j$, at most $42(j-1)$ previously used
vertices other than $x_{j-1}$ and $j-1$
connecting colors exclude at most
$
42(j-1)+(j-1)=43(j-1)<129t
$
of the $384t$ choices. Hence $c_j$ can be chosen
so that $Z_j(x_{j-1},c_j)$ avoids all previously
used vertices and $c_1,\ldots,c_{j-1}$.
Add all the corresponding walks and let
$x_j=x_{j-1}+c_j+g\in U$.

Write $w=x_{3t}\in V$,
$
F_R=\operatorname{col}(\mathcal P)
\mathbin{\dot\cup}\{c_1,\ldots,c_{3t}\}
\subseteq F$,
and $|F_R|=9t+2\ell$. These choices are independent of $\mathcal P'$.
Now fix any $\ell$-element $\mathcal P'\subseteq\mathcal P_{\mathrm{flex}}$. Choose $p_j=(a_j,b_j)\in\mathcal P_j$ as above and follow the path with colors $c_j,a_j,b_j$ at step $j$. The resulting $x$--$w$ path $P_{\mathcal P'}$ is rainbow, has all vertices other than $x$ in $V$,
and satisfies
$
\operatorname{col}(P_{\mathcal P'})
=F_R\setminus\operatorname{col}(\mathcal P').
$ 
The rainbow path is illustrated in Figure~\ref{fig:switching-chain}.
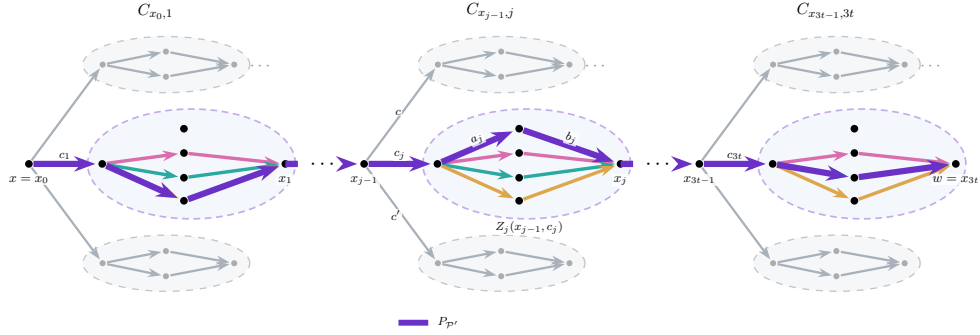
\begin{figure}[t]
\centering
\begingroup
\definecolor{pathpurple}{HTML}{6B30C6}
\definecolor{branchpink}{HTML}{DA6EAC}
\definecolor{branchteal}{HTML}{2AA9A1}
\definecolor{branchamber}{HTML}{DDA64A}
\definecolor{softblue}{HTML}{EAF2FB}
\definecolor{shadegray}{HTML}{F5F7F9}
\definecolor{fadedgray}{HTML}{99A3AD}
\tikzset{
 vtx/.style={circle,fill=black,inner sep=0pt,minimum size=5.2pt},
 dimvtx/.style={circle,fill=gray!65,inner sep=0pt,minimum size=4.1pt},
 mainpath/.style={draw=pathpurple,line width=2.75pt,-{Stealth[length=3.3mm,width=2.4mm]},shorten >=3pt,shorten <=2pt},
 alternative/.style={line width=1.5pt,-{Stealth[length=2.45mm,width=1.65mm]},shorten >=2.5pt,shorten <=2pt},
 lightpath/.style={draw=fadedgray!85,line width=1pt,-{Stealth[length=2.1mm,width=1.4mm]},shorten >=2pt,shorten <=2pt},
 mainset/.style={fill=softblue,fill opacity=0.49,draw=pathpurple!40,densely dashed,line width=.75pt},
 candidateset/.style={fill=shadegray,draw=gray!45,dashed,line width=.6pt},
 lab/.style={font=\footnotesize,fill=white,fill opacity=.78,text opacity=1,inner sep=1.2pt},
 smalllab/.style={font=\scriptsize,fill=white,fill opacity=.80,text opacity=1,inner sep=.75pt},
 }
\newcommand{\switchstage}[6]{%
 \begin{scope}[shift={(#1,0)}]
 \path[mainset] (3.88,0) ellipse [x radius=2.48,y radius=1.31];
 \path[candidateset] (3.28,2.37) ellipse [x radius=1.99,y radius=.67];
 \path[candidateset] (3.28,-2.37) ellipse [x radius=1.99,y radius=.67];
 \coordinate (#2-L) at (0,0);
 \coordinate (#2-s) at (1.75,0);
 \coordinate (#2-U) at (3.70,0.84);
 \coordinate (#2-Q) at (3.70,.26);
 \coordinate (#2-R) at (3.70,-.33);
 \coordinate (#2-B) at (3.70,-.88);
 \coordinate (#2-Rt) at (6.12,0);
 \coordinate (#2-Us) at (1.76,2.37);
 \coordinate (#2-Ua) at (3.27,2.68);
 \coordinate (#2-Ub) at (3.27,2.08);
 \coordinate (#2-Ut) at (4.90,2.37);
 \draw[lightpath] (#2-L) -- (#2-Us);
 \foreach \E in {#2-Ua,#2-Ub}{\draw[lightpath] (#2-Us) -- (\E); \draw[lightpath] (\E) -- (#2-Ut);}
 \foreach \V in {#2-Us,#2-Ua,#2-Ub,#2-Ut}{\node[dimvtx] at (\V) {};}
 \coordinate (#2-Ds) at (1.76,-2.37);
 \coordinate (#2-Da) at (3.27,-2.08);
 \coordinate (#2-Db) at (3.27,-2.68);
 \coordinate (#2-Dt) at (4.90,-2.37);
 \draw[lightpath] (#2-L) -- (#2-Ds);
 \foreach \E in {#2-Da,#2-Db}{\draw[lightpath] (#2-Ds) -- (\E); \draw[lightpath] (\E) -- (#2-Dt);}
 \foreach \V in {#2-Ds,#2-Da,#2-Db,#2-Dt}{\node[dimvtx] at (\V) {};}
 \draw[alternative,draw=branchpink] (#2-s) -- (#2-Q);
 \draw[alternative,draw=branchpink] (#2-Q) -- (#2-Rt);
 \draw[alternative,draw=branchteal] (#2-s) -- (#2-R);
 \draw[alternative,draw=branchteal] (#2-R) -- (#2-Rt);
 \draw[alternative,draw=branchamber] (#2-s) -- (#2-B);
 \draw[alternative,draw=branchamber] (#2-B) -- (#2-Rt);
 \draw[mainpath] (#2-L) -- (#2-s);
 \draw[mainpath] (#2-s) -- (#2-#6);
 \draw[mainpath] (#2-#6) -- (#2-Rt);
 \foreach \V in {#2-L,#2-s,#2-U,#2-Q,#2-R,#2-B,#2-Rt}{\node[vtx] at (\V) {};}
 \node[lab,anchor=north] at (0,-.24) {$#3$};
 \node[lab,anchor=north] at (6.12,-.24) {$#4$};
 \node[font=\normalsize] at (3.02,3.64) {$#5$};
 \node[font=\large,text=black!40] at (5.55,2.33) {$\cdots$};
 \end{scope}%
}

\resizebox{0.78\linewidth}{!}{%
\begin{tikzpicture}[font=\small,scale=0.65,transform shape]
 \switchstage{0}{one}{x=x_0}{x_1}{C_{x_0,1}}{B}
 \switchstage{8.0}{mid}{x_{j-1}}{x_j}{C_{x_{j-1},j}}{U}
 \switchstage{16.0}{last}{x_{3t-1}}{w=x_{3t}}{C_{x_{3t-1},3t}}{R}
 \draw[draw=pathpurple,line width=2.8pt] (one-Rt) -- (6.42,0);
 \node[font=\Large,fill=white,inner sep=0pt] at (7.05,0) {$\cdots$};
 \draw[mainpath] (7.67,0) -- (mid-L);
 \draw[draw=pathpurple,line width=2.8pt] (mid-Rt) -- (14.42,0);
 \node[font=\Large,fill=white,inner sep=0pt] at (15.05,0) {$\cdots$};
 \draw[mainpath] (15.67,0) -- (last-L);
 \node[smalllab,above] at (.86,.08) {$c_1$};
 \node[smalllab,above] at (16.85,.08) {$c_{3t}$};
 \node[smalllab,above] at (8.85,.07) {$c_j$};
 \node[smalllab,above,rotate=23] at (10.76,.45) {$a_j$};
 \node[smalllab,above,rotate=-18] at (12.89,.45) {$b_j$};
 \node[smalllab,anchor=east] at (8.91,1.23) {$c$};
 \node[smalllab,anchor=east] at (8.90,-1.22) {$c'$};
 \node[smalllab,below] at (11.94,-1.34) {$Z_j(x_{j-1},c_j)$};
 \draw[draw=pathpurple,line width=2.6pt] (8.82,-3.79) -- (9.53,-3.79);
 \node[font=\scriptsize,anchor=west] at (9.62,-3.79) {$P_{\mathcal P'}$};
 \end{tikzpicture}%
}
\endgroup
\caption{Constructing the rainbow path $P_{\mathcal P'}$.}
\label{fig:switching-chain}
\end{figure}
\end{proof}

The next lemma constructs a terminal absorber that uses all but one of the remaining colors, together with the colors of a suitable family of flexible $g$-pairs.

\begin{lemma}\label{lem:terminal-absorber}
Let $G$ be a finite abelian group, let $H\le G$, and let $U$ be a coset of $H$. Fix $g\in H\setminus\{0\}$, $e\in G$, and an integer $\ell\ge1$. Let $\mathcal P_{\mathrm{flex}}$ be a family of at least $\ell$ pairwise color-disjoint $g$-pairs in $H$. Let $C$ and $V_T$ be subsets of $G$ with $\operatorname{col}(\mathcal P_{\mathrm{flex}})\subseteq C\subseteq(H\cup\{e\})\setminus\{0\}$.

Suppose that for every $y\in U\cup(U+e)$ and distinct $a,b\in C$ with $a+b\ne0$, there are at least $5\ell+10$ pairs $(a_i,b_i)\in\mathcal P_{\mathrm{flex}}$ such that the sequences $y,y+a_i,y+a_i+a,y+a_i+a+b$ are rainbow directed paths, whose vertex sets other than $y$ are pairwise disjoint subsets of $V_T$.

Then, for every $x\in U$ and $L\subseteq C\setminus\operatorname{col}(\mathcal P_{\mathrm{flex}})$ with $|L|\le\ell$, there exist a subfamily $\mathcal P'\subseteq\mathcal P_{\mathrm{flex}}$ of size $\ell$ and a rainbow directed path $P$ starting at $x$ such that $V(P)\setminus\{x\}\subseteq V_T$ and $\operatorname{col}(P)=(L\cup\operatorname{col}(\mathcal P'))\setminus\{d\}$ for some $d\in L\cup\operatorname{col}(\mathcal P')$.
\end{lemma}

\begin{proof}
We construct a rainbow directed path $P$ starting at $x$, together with a family $\mathcal P'\subseteq\mathcal P_{\mathrm{flex}}$ of selected $g$-pairs and the set $R$ of colors in $L\cup\operatorname{col}(\mathcal P')$ not used by $P$. Our goal is to achieve
\[
|\mathcal P'|=\ell,\qquad |R|=1,\qquad
V(P)\setminus\{x\}\subseteq V_T.
\]
Throughout the construction, we have
\[
\operatorname{col}(P)\mathbin{\dot\cup}R
=L\mathbin{\dot\cup}\operatorname{col}(\mathcal P').
\]
Initially, $P=(x)$, $\mathcal P'=\varnothing$, and $R=L$.

We use two types of extensions (see Figure~\ref{fig:terminal-absorber-extensions}), each selecting a new pair $(a_i,b_i)\in\mathcal P_{\mathrm{flex}}$ and adding it to $\mathcal P'$.

\emph{Three-edge extension.}
If $a,b\in R$ are distinct and $a+b\ne0$, extend $P$ from its endpoint $y$ along
$
y\xrightarrow{a_i}y+a_i
\xrightarrow{a}y+a_i+a
\xrightarrow{b}y+a_i+a+b,
$
and replace $a,b$ in $R$ by $b_i$. This decreases $|R|$ by one.

\emph{Two-edge extension.}
For $a\in R$, choose $c\in(C\cap H)\setminus\{a,-a\}$ and append the first two edges of a path with colors $a_i,a,c$. Replace $a$ in $R$ by $b_i$, leaving $|R|$ unchanged. Such a $c$ exists since the two colors of a $g$-pair cannot both belong to $\{a,-a\}$, as $g\ne0$.

If $L=\varnothing$, we also use an initial one-edge extension: choose any $(a,b)\in\mathcal P_{\mathrm{flex}}$, append the first edge of a path with colors $a_i,a,b$, and select $(a_i,b_i)$, leaving $R=\{b_i\}$.

We verify that these operations are possible at every step. Let $y$ be the current endpoint of $P$. Since $x\in U$ and all used colors belong to $H\cup\{e\}$, with $e$ used at most once, we have $y\in W=U\cup(U+e)$.

Suppose $k<\ell$ pairs have been selected. Each operation adds at most three edges, so $|V(P)\setminus\{y\}|\le3k$. For any distinct $a,b\in C$ with $a+b\ne0$, the hypothesis gives at least $5\ell+10$ paths from $y$, pairwise vertex-disjoint outside $y$. Previously selected pairs and vertices of $P$ exclude at most
$
k+3k=4k<5\ell+10
$
of these paths. Hence some unselected pair $(a_i,b_i)$ yields a path meeting $P$ only at $y$. Moreover,
$
\{a_i,b_i\}\cap(\operatorname{col}(P)\cup R)
=\varnothing
$
by color-disjointness and
$L\cap\operatorname{col}(\mathcal P_{\mathrm{flex}})=\varnothing$.
Thus, after each extension, $P$ is rainbow, $V(P)\setminus\{x\}\subseteq V_T$, and $\operatorname{col}(P)\mathbin{\dot\cup}R=L\mathbin{\dot\cup}\operatorname{col}(\mathcal P')$.

Suppose first that $L\ne\varnothing$. Whenever $|R|\ge2$, a three-edge extension decreases $|R|$, unless $R=\{a,-a\}$. In that case, a two-edge extension replaces $a$ by $b_i\notin\{a,-a\}$, after which a three-edge extension leaves $|R|=1$. Thus at most $|L|\le\ell$ pairs are needed to reach $|R|=1$. If $L=\varnothing$, the initial one-edge extension already gives $|R|=1$.

Finally, repeat two-edge extensions until $|\mathcal P'|=\ell$, keeping $|R|=1$. Writing $R=\{d\}$, we obtain
$
\operatorname{col}(P)
=(L\cup\operatorname{col}(\mathcal P'))\setminus\{d\}$,
and
$V(P)\setminus\{x\}\subseteq V_T$,
as required.
\begin{figure}[t]
\centering
\resizebox{\linewidth}{!}{%
\begin{tikzpicture}[
    >=Latex,
    font=\scriptsize,
    dot/.style={circle,draw=black!65,fill=white,minimum size=3.2pt,inner sep=0pt},
    pair/.style={-{Latex[length=1.7mm]},line width=.75pt,draw=blue!70!black},
    used/.style={-{Latex[length=1.7mm]},line width=.75pt,draw=orange!85!black},
    unused/.style={-{Latex[length=1.7mm]},dashed,line width=.6pt,draw=black!40}
]
\begin{scope}
\node[anchor=west] at (0,1.35) {\textup{(a)} Three-edge extension};
\foreach \x/\name in {0.38/q0,2.35/q1,4.40/q2,6.62/q3}
    {\node[dot] (\name) at (\x,.78) {};}
\draw[pair] (q0)--node[above] {$a_i$} (q1);
\draw[used] (q1)--node[above] {$a$} (q2);
\draw[used] (q2)--node[above] {$b$} (q3);
\node[below=3pt of q0] {$y$};
\node[below=3pt of q1] {$y+a_i$};
\node[below=3pt of q2] {$y+a_i+a$};
\node[below=3pt of q3] {$y+a_i+a+b$};
\node at (3.50,.07) {$R\longmapsto(R\setminus\{a,b\})\cup\{b_i\}$};
\end{scope}
\begin{scope}[xshift=7.8cm]
\node[anchor=west] at (0,1.35) {\textup{(b)} Two-edge extension};
\foreach \x/\name in {0.38/r0,2.35/r1,4.40/r2,6.62/r3}
    {\node[dot] (\name) at (\x,.78) {};}
\draw[pair] (r0)--node[above] {$a_i$} (r1);
\draw[used] (r1)--node[above] {$a$} (r2);
\draw[unused] (r2)--node[above] {$c$} (r3);
\node[below=3pt of r0] {$y$};
\node[below=3pt of r1] {$y+a_i$};
\node[below=3pt of r2] {$y+a_i+a$};
\node[below=3pt of r3] {$y+a_i+a+c$};
\node at (3.50,.07) {$R\longmapsto(R\setminus\{a\})\cup\{b_i\}$};
\end{scope}
\end{tikzpicture}%
}
\caption{The two extensions in the terminal absorber. Both select a new $g$-pair $(a_i,b_i)$ and leave $b_i$ unused; the dashed edge of color $c$ in (b) is not traversed.}
\label{fig:terminal-absorber-extensions}
\end{figure}
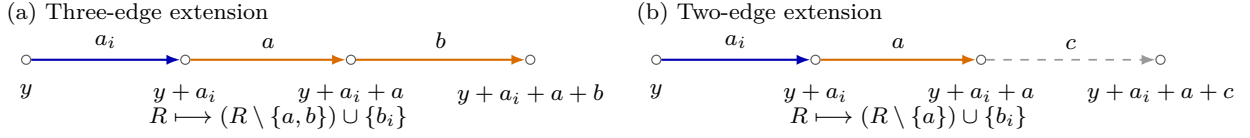
\end{proof}

\subsection{Constructing the rainbow path}
We adapt the dense-case rainbow-path construction of Bedert, Bucić, Kravitz, Montgomery and Müyesser~\cite[Theorem~9.11]{BBKMM}, which builds on the expansion-based approach of Bucić, Frederickson, Müyesser, Pokrovskiy and Yepremyan~\cite{BucicEtAl}. The main difference is that our dense color set lies in a subgroup $Q\le G$, while a small number of colors may lie outside $Q$. We use a Cayley graph diameter argument to pass to an expanding subgroup and adapt the switching and absorption steps to the relevant cosets. 

We first recall two rainbow-path results from~\cite{BBKMM}.
\begin{theorem}[{\cite[Theorem~9.10]{BBKMM}}]
\label{thm:three-quarter-path}
For every sufficiently large integer $N$, every group $G$ of order $N$ and every subset $S\subseteq G\setminus \{e\}$ with $|S|\ge 3N/4$, the Cayley digraph $\Cay(G,S)$ contains a rainbow path of length $|S|-1$.
\end{theorem}
\begin{lemma}[{\cite[Lemma~5.7]{BBKMM}}]
\label{lem:bbkmm-flexible-path}
Let $G$ be a finite abelian group of order $N$, and let $8N^{-1/9500}\le\zeta,\mu,\varepsilon,q\le1$. Suppose that
\begin{itemize}
\item[(i)] $S\subseteq G\setminus\{0\}$ satisfies $|S|\ge\varepsilon N$, and $\Cay(G,S)$ has no $\zeta$-sparse cut;
\item[(ii)] $J\subseteq G$ satisfies $|J|\le2^{-28}q^3\mu^3\varepsilon^2\zeta^2N$;
\item[(iii)] $M\subseteq G$ is a $q$-random subset, where $q\ge(1+\mu)|S|/N$;
\item[(iv)] $S'\subseteq S$ is a $q'$-random subset, where $0\le q'\le1-\mu q/4$.
\end{itemize}
Assume that $M$ and $S'$ are chosen independently. Then, with probability at least $1-7/N$, for every $S_F\subseteq S'$ and every two distinct vertices $x,y\in G$, there is a rainbow directed $x$--$y$ path in $\Cay(G,S\setminus S_F)$ using all but at most $\mu qN$ colors of $S\setminus S_F$, with all internal vertices in $M\setminus J$.
\end{lemma}
Now we are ready to prove Theorem~\ref{lem:dense-core-sparse-outside}.

\begin{proof}[Proof of Theorem~\ref{lem:dense-core-sparse-outside}]
\setcounter{denseproperty}{0}
Let $S=A\cup E$. By Lemma~\ref{lem:dense-expanding-subgroup}, applied with $\varepsilon=h^{-\gamma}$, there is a subgroup $H\le Q$ such that, writing $A_0=A\cap H$, we have $|A_0|\ge(1-h^{-\gamma})|A|\ge|A|/2$. Moreover, $\Cay(H,A_0)$ has no $\zeta$-sparse cut, where $\zeta=h^{-3\gamma}/128$. Write $J=S\setminus H$ and $N=|H|$. Since $|E|\le h^{1-2\gamma}\le h^{-\gamma}|A|$, we have $|J|\le2h^{-\gamma}|A|$. If $J=\varnothing$ and $|A_0|>3N/4$, the result follows from Theorem~\ref{thm:three-quarter-path} and Lemma~\ref{lem:cayley-walk}.

We next choose subsets $F,B$ of $A_0$ by assigning each color independently to $F$ and to $B$, with probabilities $1/4$ and $1/4$, respectively. Then we assign the remaining colors depending on the density of $A_0$ in $H$:
\begin{itemize}
\item \textbf{Type~I:} $|A_0|\le 3N/4$. Let $A_1=A_0$ and $A_2'=\varnothing$. By Lemma~\ref{lem:chernoff}\textup{(i)},
\begin{equation}\label{8.3-1}
\Prob\left(\min\{|F|,|B|\}\ge\frac{|A|}{16}\right)=1-o(1).
\end{equation}

\item \textbf{Type~II:} $|A_0|>3N/4$. Independently assign each color of $A_0\setminus(F\cup B)$ to $A_1'$ or $A_2'$, each with probability $1/2$, and let $A_i=A_i'\cup F\cup B$ for $i\in[2]$. By Lemma~\ref{lem:chernoff}\textup{(i)} and \textup{(ii)},
\begin{equation}\label{8.3-2}
\Prob\left(\min\{|F|,|B|\}\ge\frac{|A|}{16},\;
\frac{17N}{32}\le |A_i|\le\frac{25N}{32}
\text{ for all }i\in[2]\right)=1-o(1).
\end{equation}
\end{itemize}
In either type, $A_0=A_1\mathbin{\dot\cup}A_2'$. Conditional on $A_1$ in Type~I, the set $F\cup B$ is a $1/2$-random subset of $A_1$; conditional on either $A_i$ in Type~II, it is a $2/3$-random subset of $A_i$.

If Type~II holds and $\sigma(J)\in H$, choose $e\in J$; otherwise let $e=0$. This is possible since $J\ne\varnothing$ in Type~II. Set $U=\sigma(J\setminus\{e\})+H$, so that $U\ne H$ in Type~II. Let $p=1/32$, $q=1-3p$, and $\mu=h^{-2\gamma}/2^{60}$. Independently of the color partition, assign each vertex of $G$ to $V_A,V_R,V_M,V_T$ with probabilities $p,p,q,p$, respectively.

We will construct a rainbow directed path $P_1P_AP_RP_2P_T$, with $P_1$ trivial in Type~I. We first establish five properties simultaneously with probability $1-o(1)$. Properties (U1) and (U2) will provide the long paths $P_2$ in $U$ and $P_1$ in $H$, respectively. For (U1) and, in Type~II, (U2), we condition on $A_1$ and $A_2$, respectively, leaving $F\cup B$ random. For (U3)--(U5), we fix a color partition satisfying the bounds in \eqref{8.3-1} or \eqref{8.3-2}, and use only the randomness of the vertex partition.

Fix $u_0\in U$ and let $M_1=(V_M\cap U)-u_0$. For Lemma~\ref{lem:bbkmm-flexible-path}, take
$
(G,S,J,S',S_F,M,\varepsilon)_
{\text{Lemma~\ref{lem:bbkmm-flexible-path}}}
=(H,A_1,\varnothing,F\cup B,D,M_1,\varepsilon_0),
$
where $\varepsilon_0=h^{-\gamma}/4$ and $D\subseteq F\cup B$. Let us check the assumptions of Lemma~\ref{lem:bbkmm-flexible-path}. 
We take $q'=1/2$ in Type~I and $q'=2/3$ in Type~II.
Conditional on $A_1$, the set $F\cup B$ is a $q'$-random subset
of $A_1$, while $M_1$ is a $q$-random subset of $H$ independent
of $F\cup B$. Since $A_1=A_0$ in Type I, $\Cay(H,A_0)$ has no $\zeta$-sparse cut. In Type~II, the following argument applies to both $A_1$ and $A_2$. Every vertex of $\Cay(H,A_i)$ has in-degree and
out-degree $|A_i|$, so $e(X,H\setminus X)=e(H\setminus X,X)$ for every $X\subseteq H$. It therefore suffices to consider $|X|\le N/2$. For such $X$,
$
e(X,H\setminus X)
\ge |X|(|A_i|-|X|)
\ge \frac{|X|N}{32}
\ge \zeta|X|(N-|X|).
$
Also, the bounds in \eqref{8.3-1} or \eqref{8.3-2} ensure that (i) holds, while $q'\le2/3<1-\mu q/4$. Finally,
$\min\{\zeta,\mu,\varepsilon_0,q\}\ge8N^{-1/9500}$ for sufficiently
large $h$, since $N\ge|A_0|\ge h^{1-\gamma}/2$ and
$3\gamma<(1-\gamma)/9500$.

Conditional on $A_1$, Lemma~\ref{lem:bbkmm-flexible-path}, followed by translation by $u_0$, implies the following with probability at least $1-7/N$:

\medskip
\refstepcounter{denseproperty}\noindent\textbf{\thedenseproperty}\label{U1}
For every $D\subseteq F\cup B$ and distinct vertices $x,y\in U$,
there is a rainbow directed $x$--$y$ path in
$\Cay(G,A_1\setminus D)$ whose vertices lie in $U$, whose internal vertices lie in $V_M$, and which uses all but at most $\mu qN$ colors of $A_1\setminus D$.
\medskip

In Type~II, take $A_2$ in place of $A_1$ and $V_M\cap H$ in
place of $M_1$, with all other parameters unchanged. Conditional
on $A_2$, the set $F\cup B$ is a $2/3$-random subset of $A_2$,
independent of the $q$-random set $V_M\cap H$. The same argument
gives with probability at least $1-7/N$:

\medskip
\refstepcounter{denseproperty}\noindent\textbf{\thedenseproperty}\label{U2}
For every $D\subseteq F\cup B$ and distinct vertices $x,y\in H$,
there is a rainbow directed $x$--$y$ path in
$\Cay(H,A_2\setminus D)$ whose internal vertices lie in $V_M$,
and which uses all but at most $\mu qN$ colors of
$A_2\setminus D$.
\medskip

We next construct the short path $P_A$ from $H$ to $U$, using every color of $J\setminus\{e\}$ and auxiliary colors from $B$. Write $J\setminus\{e\}=\{j_1,\ldots,j_m\}$ and let $c_i=j_1+\cdots+j_i+H$ for $i\in[m]$, with $c_0=H$. If $m=0$, the one-vertex path suffices. Suppose $m\ge1$. For $i\in[m]$ and $z\in c_{i-1}$, let $\mathcal B_i(z)=\{b\in B:z+b,z+b+j_i\in V_A\}$. Since $j_i\notin H$, the sets $\{z+b,z+b+j_i\}$, $b\in B$, are pairwise disjoint. Hence $|\mathcal B_i(z)|\sim\operatorname{Bin}(|B|,p^2)$. Since $5m\le |B|p^2/2$ for sufficiently large $h$, by Lemma~\ref{lem:chernoff}\textup{(i)}, we have $\Prob(|\mathcal B_i(z)|<5m)\le\Prob(|\mathcal B_i(z)|\le |B|p^2/2)\le\exp(-|B|p^2/8)$. There are at most $mN=O(h^2)$ choices of $(i,z)$, and $|B|p^2\gg\log h$. Thus, with probability $1-o(1)$, every $\mathcal B_i(z)$ contains at least $5m$ colors.

Starting from any $x\in H$, we choose $b_i\in\mathcal B_i(z_{i-1})$ successively and append the edges $z_{i-1}\xrightarrow{b_i}z_{i-1}+b_i\xrightarrow{j_i}z_i$, where $z_i=z_{i-1}+b_i+j_i$. At step $i$, previously used colors and vertices exclude at most $3(i-1)<5m$ choices of $b_i$. Thus we can choose $b_i\in\mathcal B_i(z_{i-1})$ so that the path remains rainbow and has no repeated vertices. Hence the construction yields the following property with probability $1-o(1)$.

\medskip
\refstepcounter{denseproperty}\noindent\textbf{\thedenseproperty}\label{U3}
For every $x\in H$, there is a rainbow directed path $P_A$ starting at $x$ and ending in $U$, with $V(P_A)\setminus\{x\}\subseteq V_A$ and $\operatorname{col}(P_A)=(J\setminus\{e\})\mathbin{\dot\cup}B_A$ for some $B_A\subseteq B$.
\medskip

We next prepare the switching path $P_R$ using colors from $F$. Its endpoints will remain fixed while selected color pairs are left for the terminal path $P_T$. Let $t=\lfloor h^{1-2\gamma}/2^{20}\rfloor$ and $\ell=\lfloor tp^3/2^8\rfloor$. Of the unordered pairs of distinct colors in $F$, at most $|F|/2$ have sum zero. Averaging the remaining pairs over $H\setminus\{0\}$, we find $g\in H\setminus\{0\}$ for which at least
$
\frac{|F|(|F|-2)}{2(N-1)}
\ge\frac{|F|^2}{6N}
\ge\frac{h^{1-2\gamma}}{1536}
\ge3t+\ell
$
unordered pairs have sum $g$. These pairs are color-disjoint, since a color belongs to at most one pair with prescribed sum $g$. Choose $3t+\ell$ of these pairs, with an arbitrary ordering of the two colors in each pair, to form $\mathcal P$. Let $\mathcal P_{\mathrm{flex}}\subseteq\mathcal P$ be any subfamily of size $t+\ell$. Since both families are chosen using only $F$, they are independent of $V_R$ and $V_T$.

Since $|F|p^{42}\ge2^{20}\max\{3t,\log N\}$ for large $h$, Lemma~\ref{lem:switching-chain} applies in the coset $U$ with the $p$-random set $V_R\cap U$. Consequently, with probability at least $1-1/N$ we have the following property.

\medskip
\refstepcounter{denseproperty}\noindent\textbf{\thedenseproperty}\label{U4}
For every $z\in U$, there are $w\in V_R\cap U$ and $F_R\subseteq F$ such that $\operatorname{col}(\mathcal P)\subseteq F_R$, $|F_R|=9t+2\ell$, and, for every $\ell$-element subfamily $\mathcal P'\subseteq\mathcal P_{\mathrm{flex}}$, there is a rainbow directed $z$--$w$ path $P_{\mathcal P'}$ with $V(P_{\mathcal P'})\setminus\{z\}\subseteq V_R$ and $\operatorname{col}(P_{\mathcal P'})=F_R\setminus\operatorname{col}(\mathcal P')$.
\medskip

Finally, we establish \ref{U5}. Let $C=A_0\cup(\{e\}\setminus\{0\})$. Then $\operatorname{col}(\mathcal P_{\mathrm{flex}})\subseteq C\subseteq(H\cup\{e\})\setminus\{0\}$. Fix $y\in U\cup(U+e)$ and distinct $a,b\in C$ with $a+b\ne0$. Write $\mathcal P_{\mathrm{flex}}=\{(a_i,b_i):i\in[t+\ell]\}$. For each $i$, consider the three-edge walk $y\xrightarrow{a_i}y+a_i\xrightarrow{a}y+a_i+a\xrightarrow{b}y+a_i+a+b$. It is a rainbow path unless $a_i\in\{0,-a,-a-b,a,b\}$, which excludes at most five indices.

The three vertices other than $y$ form $Z_i=y+a_i+D$, where $D=\{0,a,a+b\}$. If $Z_i\cap Z_j\ne\varnothing$, then $a_j-a_i\in D-D$. Since $|D-D|\le9$, a greedy choice yields at least $(t-5)/10$ pairwise disjoint sets $Z_i$ corresponding to rainbow paths. These paths survive in $V_T$ independently with probability $p^3$, so the expected number of surviving paths is at least $(t-5)p^3/10$. Since $5\ell+10\le(t-5)p^3/20$ for sufficiently large $h$, by Lemma~\ref{lem:chernoff}\textup{(i)}, the probability that fewer than $5\ell+10$ paths survive is at most $\exp(-(t-5)p^3/80)$. Since $tp^3\gg\log h$, a union bound over at most $2N(N+1)^2$ choices of $(y,a,b)$ shows that, with probability $1-o(1)$, the hypothesis of Lemma~\ref{lem:terminal-absorber} holds. By that lemma, we obtain the following property.

\medskip
\refstepcounter{denseproperty}\noindent\textbf{\thedenseproperty}\label{U5}
For every $y\in U$ and $L\subseteq C\setminus\operatorname{col}(\mathcal P_{\mathrm{flex}})$ with $|L|\le\ell$, there exist an $\ell$-element subfamily $\mathcal P'\subseteq\mathcal P_{\mathrm{flex}}$ and a rainbow directed path $P_T$ starting at $y$ such that $V(P_T)\setminus\{y\}\subseteq V_T$ and $\operatorname{col}(P_T)=(L\cup\operatorname{col}(\mathcal P'))\setminus\{d\}$ for some $d\in L\cup\operatorname{col}(\mathcal P')$.
\medskip

The random variables $|V_A\cap H|$, $|V_M\cap U|$, and $|V_M\cap H|$ have binomial distributions with means $pN$, $qN$, and $qN$, respectively. For sufficiently large $h$, by Lemma~\ref{lem:chernoff}\textup{(i)} and a union bound,
$\Prob(V_A\cap H=\varnothing\text{ or }V_M\cap U=\varnothing\text{ or }|V_M\cap H|<2)\le e^{-pN/8}+2e^{-qN/8}=o(1)$.
We may therefore fix an outcome satisfying the size estimates, \ref{U1}, \ref{U3}--\ref{U5}, and, in Type~II, \ref{U2}, together with $V_A\cap H\ne\varnothing$, $V_M\cap U\ne\varnothing$, and $|V_M\cap H|\ge2$.

In Type~I, choose $x\in V_A\cap H$ and let $P_1=(x)$. In Type~II, choose distinct $z,x\in V_M\cap H$. By \ref{U2}, with $D=F\cup B$, there is a rainbow directed
$z$--$x$ path $P_1$ using all but at most $\mu qN$
colors of $A_2'$.

By \ref{U3}, construct $P_A$ from $x$, with auxiliary colors
$B_A\subseteq B$, and let $z_A\in U$ be its endpoint.
Applying \ref{U4} at $z_A$, choose $w\in V_R\cap U$ and
$F_R\subseteq F$ as in that property. We do not yet choose the corresponding switching path. Instead, choose $y\in V_M\cap U$. Since $w\ne y$, \ref{U1}, with $D=F_R\cup B_A$, gives a rainbow directed $w$--$y$ path $P_2$ using all but at most $\mu qN$ colors of $A_1\setminus(F_R\cup B_A)$.

Let $L$ consist of the colors of $S$ not belonging to $\operatorname{col}(P_1)\cup\operatorname{col}(P_A)\cup\operatorname{col}(P_2)\cup F_R$. Then $L\subseteq C\setminus\operatorname{col}(\mathcal P_{\mathrm{flex}})$ and $|L|\le2\mu qN+1\le\ell$, where the last inequality follows from $2\mu qN\le h^{1-2\gamma}/2^{59}$ and $\ell=\lfloor h^{1-2\gamma}/2^{43}\rfloor$.

By \ref{U5}, choose $\mathcal P'$, $P_T$, and $d$ for this $L$ and starting vertex $y$. We now choose the corresponding path $P_R$ supplied by \ref{U4}, from $z_A$ to $w$, with $\operatorname{col}(P_R)=F_R\setminus\operatorname{col}(\mathcal P')$.

The vertex partition and, in Type~II, $U\ne H$ ensure that $P_1P_AP_RP_2P_T$ is a directed path as illustrated in Figure~\ref{fig:rainbow-five-paths}. It uses every color of $S\setminus\{d\}$ exactly once. Thus $S$ is rearrangeable by Lemma~\ref{lem:cayley-walk}.
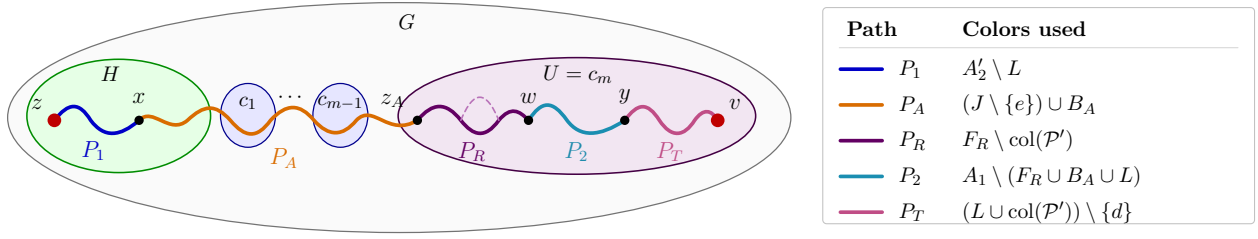
\begin{figure}[t]
\centering
\resizebox{\linewidth}{!}{%
\begin{tikzpicture}[
  x=1cm,y=1cm, line cap=round,line join=round,
  lab/.style={font=\large},
  cos/.style={font=\normalsize},
  joint/.style={circle,fill=black,inner sep=1.7pt},
  endpt/.style={circle,fill=red!75!black,inner sep=2.5pt},
  segment/.style={line width=1.9pt}
]
\filldraw[fill=black!2,draw=black!55,line width=.65pt]
  (7.25,1.87) ellipse [x radius=7.12,y radius=2.09];
\node[anchor=west,font=\normalsize] at (7.10,3.62) {$G$};
\filldraw[fill=green!10,draw=green!55!black,line width=.8pt]
  (2.15,1.90) ellipse [x radius=1.67,y radius=1.03];
\node[cos] at (2.00,2.66) {$H$};
\filldraw[fill=violet!10,draw=violet!53!black,line width=.8pt]
  (10.48,1.90) ellipse [x radius=3.25,y radius=1.06];
\node[cos] at (10.50,2.66) {$U=c_m$};
\filldraw[fill=blue!10,draw=blue!55!black,line width=.65pt]
  (4.50,1.91) ellipse [x radius=.50,y radius=.58];
\filldraw[fill=blue!10,draw=blue!55!black,line width=.65pt]
  (6.18,1.91) ellipse [x radius=.50,y radius=.58];
\node[cos] at (4.50,2.14) {$c_1$};
\node[cos] at (6.18,2.14) {$c_{m-1}$};
\node[cos] at (5.30,2.24) {$\cdots$};
\coordinate (a) at (.98,1.82);
\coordinate (b) at (2.52,1.82);
\coordinate (c) at (7.58,1.82);
\coordinate (d) at (9.60,1.82);
\coordinate (e) at (11.35,1.82);
\coordinate (f) at (13.04,1.82);
\draw[segment,blue!78!black]
  (a) .. controls (1.24,2.15) and (1.43,2.16) .. (1.64,1.82)
      .. controls (1.89,1.49) and (2.17,1.52) .. (b);
\draw[segment,orange!85!black]
  (b) .. controls (2.77,2.08) and (3.07,1.54) .. (3.40,1.82)
      .. controls (3.70,2.10) and (3.90,2.13) .. (4.18,1.82)
      .. controls (4.45,1.50) and (4.66,1.48) .. (4.91,1.82)
      .. controls (5.20,2.13) and (5.44,2.13) .. (5.70,1.82)
      .. controls (5.97,1.52) and (6.22,1.52) .. (6.47,1.82)
      .. controls (6.78,2.14) and (7.05,1.50) .. (c);
\draw[segment,violet!78!black]
  (c) .. controls (7.86,2.15) and (8.10,2.16) .. (8.37,1.82)
      .. controls (8.60,1.53) and (8.81,1.48) .. (9.06,1.82)
      .. controls (9.28,2.08) and (9.43,2.00) .. (d);
\draw[violet!55,line width=.85pt,densely dashed]
  (8.37,1.82) .. controls (8.60,2.38) and (8.83,2.38) .. (9.06,1.82);
\draw[segment,cyan!66!black]
  (d) .. controls (9.84,2.18) and (10.03,2.20) .. (10.30,1.82)
      .. controls (10.54,1.50) and (10.89,1.52) .. (e);
\draw[segment,magenta!77!black]
  (e) .. controls (11.59,2.14) and (11.76,2.16) .. (11.98,1.82)
      .. controls (12.20,1.49) and (12.43,1.52) .. (12.62,1.82)
      .. controls (12.78,2.04) and (12.92,2.02) .. (f);
\foreach \q in {b,c,d,e}{\node[joint] at (\q) {};}
\node[endpt] at (a) {};
\node[endpt] at (f) {};
\node[font=\large,above left=3pt] at (a) {$z$};
\node[font=\large,above=5pt] at (b) {$x$};
\node[font=\large,above left=4pt] at (c) {$z_A$};
\node[font=\large,above=5pt] at (d) {$w$};
\node[font=\large,above=5pt] at (e) {$y$};
\node[font=\large,above right=3pt] at (f) {$v$};
\node[lab,blue!78!black] at (1.69,1.27) {$P_1$};
\node[lab,orange!85!black] at (5.15,1.13) {$P_A$};
\node[lab,violet!78!black] at (8.57,1.25) {$P_R$};
\node[lab,cyan!63!black] at (10.47,1.25) {$P_2$};
\node[lab,magenta!77!black] at (12.17,1.25) {$P_T$};

\draw[black!24,rounded corners=2pt,line width=.55pt,fill=white]
  (14.95,-.06) rectangle (22.80,3.87);
\node[anchor=west,font=\bfseries\normalsize] at (15.25,3.49) {Path};
\node[anchor=west,font=\bfseries\normalsize] at (17.35,3.49) {Colors used};
\draw[black!17,line width=.4pt] (15.25,3.14)--(22.50,3.14);
\draw[blue!78!black,line width=1.9pt] (15.27,2.75)--(16.03,2.75);
\node[anchor=west,font=\normalsize] at (16.20,2.75) {$P_1$};
\node[anchor=west,font=\normalsize] at (17.35,2.75) {$A_2'\setminus L$};
\draw[orange!85!black,line width=1.9pt] (15.27,2.10)--(16.03,2.10);
\node[anchor=west,font=\normalsize] at (16.20,2.10) {$P_A$};
\node[anchor=west,font=\normalsize] at (17.35,2.10) {$(J\setminus\{e\})\cup B_A$};
\draw[violet!78!black,line width=1.9pt] (15.27,1.45)--(16.03,1.45);
\node[anchor=west,font=\normalsize] at (16.20,1.45) {$P_R$};
\node[anchor=west,font=\normalsize] at (17.35,1.45) {$F_R\setminus\operatorname{col}(\mathcal P')$};
\draw[cyan!66!black,line width=1.9pt] (15.27,.80)--(16.03,.80);
\node[anchor=west,font=\normalsize] at (16.20,.80) {$P_2$};
\node[anchor=west,font=\normalsize] at (17.35,.80) {$A_1\setminus(F_R\cup B_A\cup L)$};
\draw[magenta!77!black,line width=1.9pt] (15.27,.15)--(16.03,.15);
\node[anchor=west,font=\normalsize] at (16.20,.15) {$P_T$};
\node[anchor=west,font=\normalsize] at (17.35,.15) {$(L\cup\operatorname{col}(\mathcal P'))\setminus\{d\}$};
\end{tikzpicture}%
}
\caption{The rainbow path $P_1P_AP_RP_2P_T$, with the color sets of its five segments listed on the right.}
\label{fig:rainbow-five-paths}
\end{figure}
\end{proof}

\subsection{Proof of Theorem~\ref{thm:asymptotic}}

\begin{proof}[Proof of Theorem~\ref{thm:asymptotic}]
Let $S\subseteq G\setminus\{0\}$ have sufficiently large size $n$. By Theorem~\ref{thm:main}, we may assume that $\diss(S)\le C_0(\log n)^2$ for an absolute constant $C_0$. Set $\gamma=2^{-20}$, $\delta=\gamma/10$, and $\alpha=\delta/100$.

By Lemma~\ref{lem:concentrated-core}, there is a coset $a+H$ such that $A=S\cap(a+H)$ satisfies $|S\setminus A|\le n^{2\alpha}$ and $|A\setminus(x+K)|\ge\lceil n^\alpha\rceil$ for every proper subgroup $K<H$ and every $x\in G$. Write $h=|H|$ and $q=\ord_{G/H}(a+H)$.

If $h\ge n^{1+\delta}$, then $S$ is rearrangeable by Lemma~\ref{lem:large-subgroup-core}. If $q\ge n^\delta$, the same follows from Theorem~\ref{thm:large-order-absorption}, applied with $\varepsilon=2\alpha<\delta$. We may therefore assume that $h<n^{1+\delta}$ and $q<n^\delta$.

Let $Q=H+\langle a\rangle$. Then $a+H\subseteq Q$ and $|Q|=hq<n^{1+2\delta}$. For $A_Q=S\cap Q$ and $E_Q=S\setminus Q$, we have $|A_Q|\ge n-n^{2\alpha}$ and $|E_Q|\le n^{2\alpha}$, so $|Q|\ge n/2$ for sufficiently large $n$. Since $(1+2\delta)(1-\gamma)<1$ and $2\alpha<1-2\gamma$, we obtain $|Q|^{1-\gamma}\le n^{(1+2\delta)(1-\gamma)}\le n-n^{2\alpha}\le|A_Q|$ and $|E_Q|\le n^{2\alpha}\le(n/2)^{1-2\gamma}\le|Q|^{1-2\gamma}$. By Theorem~\ref{lem:dense-core-sparse-outside}, $S=A_Q\cup E_Q$ is rearrangeable.
\end{proof}

\section*{AI declaration}
During this work, we used ChatGPT 5.4, 5.5, 5.6 and 6 Pro for mathematical discussion, numerical optimization, assistance with drafting and checking proofs, and language editing. All figures in this paper were based on ideas provided by the authors with the help of AI.

Early in 2026, inspired by the use of dissociated sets in Bedert and Kravitz~\cite{BedertKravitz}, we began studying the relationship between dissociated dimension and Graham's rearrangement conjecture and obtained Theorem~\ref{thm:main}. The main ideas and arguments were developed by the authors. ChatGPT 5.4 and 5.5 Pro assisted with numerical optimization of the parameters and constants. Subsequent discussions with AI led to an asymptotically sharp refinement of Theorem~\ref{thm:main}, which is not presented in this paper.

Section~\ref{sec:large-subgroup} was also developed primarily
by the authors. For Sections~\ref{sec:large-order-coset}
and~\ref{sec:dense-core}, discussions with ChatGPT 5.6
and 6 Pro helped identify adaptations of existing methods
and develop preliminary proof drafts. In Section~\ref{sec:large-order-coset}, these concerned the local-repair arguments
of~\cite{CostaEtAlKneser,PhamSauermann}, leading to
Lemma~\ref{lem:large-order-forbidden} and
Theorem~\ref{thm:large-order-absorption}.
In Section~\ref{sec:dense-core}, they concerned the
expansion approach of~\cite{BucicEtAl} and the absorption
argument of~\cite[Section~9]{BBKMM}, contributing to
Lemmas~\ref{lem:dense-expanding-subgroup},
\ref{lem:switching-chain}
and~\ref{lem:terminal-absorber}, and their applications
in Theorem~\ref{lem:dense-core-sparse-outside}.
The authors studied the original papers, determined
the necessary modifications, and revised and checked
the resulting arguments.

The authors directed all stages of the research and made the final decisions concerning the statements, arguments, and presentation. They take full responsibility for the content of this paper, including the correctness of all results and proofs.

\appendix
\section{Proof of Lemma \ref{lem:nested-to-local-remote}}\label{app:local-repair}
In this section, we prove Lemma~\ref{lem:nested-to-local-remote}.

We use the following formulation of \cite[Definition~3.7 and Corollary~3.6]{CostaEtAlKneser}.

\begin{definition}\label{def:local-remote}
Let $D\ge1$ and $K\ge0$ be integers, and let $\Theta>0$. The random ordering has the \emph{$(D,K,\Theta)$ local--remote chain property} with constant $B$ if, for every $1\le h\le D$, $c\in[n]$, and $F\subseteq[1,c]$ with $|F|\le K$,
\begin{equation}\label{eq:right-local-remote}
\sum_{\substack{c<b_1<\cdots<b_h\le n\\b_h-c<n-|F|}}
\Prob\left(\sum_{j=c+1}^{b_\ell}X_j=Z_\ell(\mathbf b)\text{ for all }\ell\in[h]\,\middle|\,\mathcal F_F\right)\le B\Theta^h
\end{equation}
for arbitrary $G$-valued $\mathcal F_F$-measurable prescribed sums $Z_\ell(\mathbf b)$, where $\mathbf b=(b_1,\ldots,b_h)$. We also require the reversed estimate: if $F\subseteq[c,n]$ and $|F|\le K$, then
\begin{equation}\label{eq:left-local-remote}
\sum_{\substack{1\le a_1<\cdots<a_h<c\\c-a_1<n-|F|}}
\Prob\left(\sum_{j=a_\ell}^{c-1}X_j=Z_\ell(\mathbf a)\text{ for all }\ell\in[h]\,\middle|\,\mathcal F_F\right)\le B\Theta^h
\end{equation}
for arbitrary $\mathcal F_F$-measurable prescribed sums, where $\mathbf a=(a_1,\ldots,a_h)$. The constant $B$ is uniform over these choices.
\end{definition}

\begin{lemma}[Costa, Della Fiore, Feng and Liu~{\cite[Lemma~4.1]{CostaEtAlKneser}}]
\label{lem:costa-local-repair}
Fix $\beta>0$. There exists an integer $D=D(\beta)\ge2$
with the following property. Let $A\subseteq G\setminus\{0\}$
be an $n$-element subset of a finite abelian group, let
$\sigma:[n]\to A$ be a uniformly random bijection, and let
$0<\Theta\le1$. Suppose that $\sigma$ has the
$(D,100D^2,\Theta)$-local--remote chain property of
Definition \ref{def:local-remote}, with constant $C_D$.
If $\Theta\le n^{-\beta}$, then there exists an integer
$
n_0=n_0(\beta,C_D)
$
such that $A$ admits a valid ordering whenever $n\ge n_0$.
\end{lemma}

The proof of \cite[Lemma~4.1]{CostaEtAlKneser} allows any fixed integer $D>3/\beta$.
Now we prove Lemma~\ref{lem:nested-to-local-remote}.
\begin{proof}[Proof of Lemma~\ref{lem:nested-to-local-remote}]
Fix $1\le h\le D$. We verify the two estimates in Definition~\ref{def:local-remote} with $K=100D^2$ and $\Theta=n^{-\beta}$.

Fix $c\in[n]$ and $F\subseteq[1,c]$ with $|F|\le100D^2$, and condition on $\mathcal F_F$. For each tuple $c<b_1<\cdots<b_h\le n$ satisfying $b_h-c<|S_F|$, define $R_\ell=\{X_j:c<j\le b_\ell\}$ for $\ell\in[h]$. The entries at the unexposed positions form a uniformly random ordering of $S_F$. Since $[c+1,b_h]\cap F=\varnothing$, the sets $R_1\subsetneq\cdots\subsetneq R_h\subsetneq S_F$ form a uniformly random strict nested chain with $|R_\ell|=b_\ell-c$.

Distinct endpoint tuples correspond to distinct size tuples $(m_1,\ldots,m_h)$, where $m_\ell=b_\ell-c$. Since the targets $Z_\ell(\mathbf b)$ are fixed after conditioning on $\mathcal F_F$, we have
\begin{align*}
&\sum_{\substack{c<b_1<\cdots<b_h\le n\\b_h-c<|S_F|}}
\Prob\left(\sum_{j=c+1}^{b_\ell}X_j=Z_\ell(\mathbf b)\ (\ell\in[h])\,\middle|\,\mathcal F_F\right)\\
&\quad\le\sum_{1\le m_1<\cdots<m_h<|S_F|}\sup_{z_1,\ldots,z_h\in G}
\Prob\left(\sigma(R_i)=z_i\ (i\in[h])\,\middle|\,\mathcal F_F\right)
\le Bn^{-\beta h},
\end{align*}
where, for each size tuple on the right, the chain has $|R_i|=m_i$. The last inequality follows from \eqref{eq:nested-chain-sum}, proving \eqref{eq:right-local-remote}.

For $F\subseteq[c,n]$ with $|F|\le100D^2$, consider $1\le a_1<\cdots<a_h<c$ satisfying $c-a_1<|S_F|$. Conditional on $\mathcal F_F$, the sets $R_i=\{X_j:a_{h+1-i}\le j<c\}$ form a uniformly random strict nested chain in $S_F$ with $|R_i|=c-a_{h+1-i}$. Applying the same argument with the targets $Z_{h+1-i}(\mathbf a)$ proves \eqref{eq:left-local-remote}.

Thus the random ordering has the $(D,100D^2,n^{-\beta})$ local--remote chain property with constant $B$. By Lemma~\ref{lem:costa-local-repair} and the choice of $D$ above, $S$ is rearrangeable for sufficiently large $n$.
\end{proof}
\end{document}